\let\ReallyTheEnd=\end
\def\end{\vfil\eject}
%


\input psfig

\newif\ifboxfigure      
\boxfigurefalse

\def\BoxIt#1#2{
	\vbox{\hrule
	\hbox{\vrule\kern#2\vbox{\kern#2#1\kern#2}\kern#2\vrule}
		   \hrule}}

\def\insertRaster #1 pixels #2  by #3 scaled #4 {
			\medskip
			 \hbox to \hsize{%

			 \hss
			 \RasterBox {#1} {#2} {#3} {#4}
			 \hss
			 }%
}

\def\RasterBox #1 #2 #3 #4{


\dimen5=65pt
\divide\dimen5 by 72

\dimen0=#2\dimen5
\divide\dimen0 by 1000
\dimen1=#3\dimen5
\divide\dimen1 by 1000
\dimen2=#3\dimen5
\divide\dimen2 by 1000
\dimen3=#2\dimen5
\divide\dimen3 by 1000

\setbox4=\hbox to #4\dimen0{
 \vbox to #4\dimen1{
 \vss
 \psfig{figure=#1,height=#4\dimen2,width=#4\dimen3}
 }
 \hss
 }
 \ifboxfigure\BoxIt{\box4}{0pt}
 \else\box4
 \fi
 }


\input psfig
\def\QP{\narrower\smallskip\noindent}
\def\ref{\hangindent=1pc \hangafter=1 \noindent}
\def\QED{  \rlap{$\sqcup$}$\sqcap$ \smallskip}
\def\mod{{\rm mod\;}}
\def\[{$\,}
\def\]{\,$}
\def\C{{\bf C}}
\def\R{{\bf R}}
\def\Z{{\bf Z}}

\def\={\;=\; }

\font\bit=cmssi12 at 12truept
\font\tenmsy=msym10

\textfont8=\tenmsy
\mathchardef\ssm="7872
\mathsurround = 1pt
\abovedisplayskip=6pt
\belowdisplayskip=6pt
\parskip=2pt

\def\Rat{{\rm Rat}}
\def\Map{{\rm Map}}
\def\Per{{\rm Per}}
\def\Aut{{\rm Aut}}
\def\result{{\rm result}}
\def\CP{{\bf CP}}
\def\RP{{\bf RP}}
\def\sl{{\cal S}}
\def\QP{\medskip\leftskip=.4in\rightskip=.4in\noindent}
\def\subarr{\hookrightarrow}
\def\darrow{\rightarrow\!\!\!\!\!\rightarrow}
\def\M{{\cal M}_2}
\def\fm{^{\rm fm}}
\def\cm{^{\rm cm}}
\def\tm{^{\rm tm}}
\def\PSL{{\rm PSL}}
\def\SO{{\rm SO}}
\def\PSU{{\rm PSU}}

\centerline{\bf Remarks on Quadratic Rational Maps.}\smallskip

\centerline{J. Milnor} \smallskip

\centerline{Stony Brook, August 1992}
\bigskip\medskip

\centerline{\bf Contents.}\medskip

\halign{\hskip .7in # \hfil&\hskip .2in \hfil #\cr
\S1. Introduction. & 1\cr
\noalign{\kern 2pt}

{\bf Geometry:} &\cr\noalign{\kern 1pt}

\S2. The Space \[\Rat_2\]
of Quadratic Rational Maps. & 2\cr\noalign{\kern 1pt}

\S3. The Space \[{\cal M}_2\] of Holomorphic Conjugacy Classes. & 3\cr\noalign
{\kern 1pt}

\S4. The compactification \[\widehat\M\cong\CP^2\]. & 8\cr\noalign{\kern 1pt}


\S5. Maps with Symmetries. & 12\cr\noalign{\kern 1pt}

\S6. Maps with marked critical points or fixed points. &14\cr\noalign
{\kern 2pt}

{\bf Dynamics:} & \cr\noalign{\kern 1pt}

\S7. Hyperbolic Julia Sets and Hyperbolic Components in Moduli Space. &
 20\cr\noalign{\kern 1pt}

\S8. The ``Escape Locus'': Totally Disconnected Julia Sets. & 24\cr\noalign
{\kern 1pt}

\S9. Complex 1-Dimensional Slices.& 30\cr\noalign{\kern 1pt}

\S10. Real Quadratic Maps. & 39\cr
\noalign{\kern 6pt}

Appendix A. Resultant and Discriminant.& 44\cr\noalign{\kern 1pt}

Appendix B. The Space \[\Rat_d\] of Degree \[d\] Rational Maps.& 45\cr
\noalign{\kern 1pt}

Appendix C. Normal forms and relations between conjugacy invariants.&48\cr
\noalign{\kern 1pt}

Appendix D. Geometry of periodic orbits. & 52\cr\noalign{\kern 1pt}

Appendix E. Totally Disconnected Julia Sets in Degree \[d\].& 54\cr
\noalign{\kern 1pt}

Appendix F. A ``Sierpinski carpet'' as Julia set~~ ({\bit written with
Tan Lei}).& 57\cr\noalign{\kern 2pt}

References.&61\cr}

\bigskip\medskip

\centerline{\bf \S1. Introduction.}\medskip

This will be an expository description of quadratic rational maps.\smallskip

Sections 2 through 6 are concerned
with the geometry and topology of such maps. The space
\[\Rat_2\] of all quadratic rational maps from the
Riemann sphere to itself is a smooth complex 5-manifold, having
the homotopy type of an \[{\rm SO}(3)$-bundle
over the real projective plane. However the ``moduli space'' \[\M\],
consisting of all holomorphic conjugacy classes of maps in \[\Rat_2\],
has a much simpler structure, and
is biholomorphic to the coordinate space \[\C^2\]. (More precisely,
\[\M\] can be described as an orbifold whose underlying space is
isomorphic to \[\C^2\].) The locus \[\Per_n(\mu)\]
consisting of conjugacy classes with a 
periodic point of period \[n\] and multiplier \[\mu\] is an algebraic curve in
\[\M\cong\C^2\]. For the special cases \[n=1\] and \[n=2\] this curve
is a straight line. The moduli space \[\M\] has a natural compactification
\[\widehat\M\], isomorphic to the projective plane
\[\CP^2\]. We also consider quadratic maps together with a marking of the
critical points, or of the fixed points. As an example, the moduli space
\[\M\cm\] for maps with marked critical points
is an orbifold with one essentially
singular point, and has the homotopy type of a 2-sphere.

Sections 7--10 survey of some topics from the
dynamics of quadratic rational maps. There are few proofs.
Those maps which are hyperbolic on their Julia set give rise
to ``hyperbolic components'' in moduli space, as studied by
Rees [R3]. If we work in the compactified
moduli space \[\widehat\M\], then every hyperbolic component is a topological
4-cell with a preferred center point. However if we work in \[\M\cong\C^2\]
then
there is one exceptional component which has a more complicated topology,
namely the ``escape component'',\break
consisting of maps with totally disconnected Julia set.
Section 9 attempts to explore and picture moduli
space by means of complex one-dimensional slices. (Compare Rees [R4], [R5].)
Section 10 describes the theory of real quadratic rational maps.\smallskip

For convenience in exposition, some technical details have
been relegated to appendices:
Appendix A outlines some classical algebra. Appendix B describes the
topology of the space of rational maps of degree \[d\]. Appendix C outlines
several
convenient normal forms for quadratic rational maps, and computes relations
between various invariants.\break Appendix D describes some geometry associated
with the curves \[\Per_n(\mu)\subset\M\].\break
Appendix E describes
totally disconnected Julia sets containing no critical points.\break Finally,
Appendix F, written in collaboration with Tan Lei, describes an example
of a connected quadratic Julia set for which no two components of the
complement have a common boundary point.
\bigskip

\centerline{\bf\S2. The Space \[\Rat_2\] of Quadratic Rational Maps.}
\medskip

This section will set the stage by giving a brief description of the space
of all quadratic rational maps.\smallskip

It will be convenient to identify the compactified plane \[\hat\C=\C\cup
\infty\] with the unit sphere \[S^2\] via stereographic projection, and to call
either one the {\bit Riemann sphere\/}. Let \[\Rat_d\] be the space consisting
of all holomorphic maps of degree \[d\] from \[S^2\] to itself.
Information about \[\Rat_d\] may be found in
Appendix B. (Compare Segal [Se].) For example, \[\Rat_d\]
is a smooth connected complex manifold of dimension \[2d+1\], and the
fundamental group \[\pi_1(\Rat_d)\] is cyclic of order \[2d\] for \[d\ge 1\].
For \[d=1\], note that the space \[\Rat_1\] can be identified
with the group \[\PSL(2,\C)\]
consisting of all M\"obius transformations from the Riemann sphere to itself.
\smallskip
 
Now let us specialize to the case \[d=2\]. Each map \[f\] in the space
\[\Rat_2\], consisting of all {\bit quadratic rational maps\/},
can be expressed as a ratio
$$	f(z)~=~{p(z)\over q(z)}~=~{a_0z^2+a_1z+a_2\over b_0z^2+b_1z+b_2} $$
with degree \[d={\rm Max}\big({\rm deg}(p)\,,\,{\rm deg}(q)\big)\] equal to
two. It follows easily that \[\Rat_2\]
can be identified with the Zariski
open subset of complex projective 5-space consisting of all points\break
\[(a_0:a_1:a_2:b_0:b_1:b_2)\] in \[\CP^5\;\] for which the {\bit resultant}
$$	\result(p,q)~=~{\rm det}\left[\matrix{a_0&a_1&a_2&0\cr
0&a_0&a_1&a_2\cr b_0&b_1&b_2&0\cr 0& b_0&b_1&b_2\cr}\right] $$
is non-zero. (See Appendix A.)
The topology of this space can be described roughly as follows.

{\QP{\bf Theorem 2.1.} \it The space \[\Rat_2\] contains a compact
non-orientable manifold \[M^5\] as deformation retract. This manifold
\[M^5\] can be described as the unique non-trivial principal
\[{\rm SO(3)}$-bundle over the projective plane \[\RP^2\].\medskip}

The proof can be outlined as follows. (For details, see Appendix B.)
Every quadratic rational map has two distinct critical points
\[\omega_1\ne \omega_2\] in the Riemann sphere \[S^2\], and two distinct
critical values \[f(\omega_1)\ne f(\omega_2)\]. Let \[M^5\subset\Rat_2\]
be the sub-space consisting of quadratic rational maps such that:

(a) the two critical points \[\omega_1\,,\;\omega_2\] are ``antipodal''
in the sense that \[\omega_2=-1/\bar\omega_1~,\]

(b) the two critical values \[f(\omega_1)\,,\;f(\omega_2)\] are also
antipodal, and

(c) every point on the ``equator'' midway between the critical points maps
\hfil\break\indent\hskip .05in
to a point on the equator midway between the critical values.

\noindent It is not hard to check that \[M^5\] is indeed a smooth manifold,
embedded in \[\Rat_2\] as a\break deformation retract, and that the
continuous map \[f\mapsto\{\omega_1\,,\,\omega_2\}\] from \[M^5\] to the real
projective plane is the projection
map of a principal fibration, with fiber equal to
the group \[\SO(3)\cong\PSU(2)\subset\PSL(2,\C)\]
consisting of all rotations of the \[2$-sphere.\break Using results of Graeme
Segal, we will show in Appendix B that the fundamental group\break \[\pi_1(M^5)\cong
\pi_1(\Rat_2)\] is cyclic of order \[4\], and conclude that this bundle must
be non-trivial.\smallskip

As one immediate consequence of Theorem 1, we see that \[M^5\] (or \[\Rat_2\])
has the rational homology of a 3-sphere. However,
the 2-fold orientable covering manifold of \[M^5\]
is homeomorphic to the product \[{\rm SO}(3)\times S^2\] (hence the universal
covering of \[M^5\] is homeomorphic to \[S^3\times S^2\]).
In \S6 we will discuss the corresponding 2-sheeted covering manifold
of \[\Rat_2\]. This covering manifold can be identified with the space of
{\bit critically marked\/} quadratic rational maps, denoted by \[\Rat_2\cm\].
Its elements can be described as ordered triples
\[(f\,,\, \omega_1\,,\,\omega_2)\] where \[f\in\Rat_2\] and where \[\omega_1
\ne\omega_2\] are the two critical points of \[f\].
\bigskip\bigskip

\centerline{\bf \S3. The Space \[{\cal M}_2\]
of Holomorphic Conjugacy Classes.}\medskip

The group \[\Rat_1\cong{\rm PSL}_2\C\] of M\"obius transformations acts
on the space \[\Rat_2\] of quadratic rational maps by conjugation,
$$	g\,\in\,\Rat_1\qquad{\rm and}\qquad f\,\in\,\Rat_2\qquad{\rm yield}
	\qquad g\circ f\circ g^{-1}\,\in\,\Rat_2~.$$
Two maps in \[\Rat_2\] are said to be {\bit holomorphically conjugate\/}
if they belong to the same orbit.

{\bf Definition:}
The quotient space of \[\Rat_2\] under this action will be denoted
by \[{\cal M}_2\], and called the {\bit moduli space\/} of holomorphic
conjugacy classes \[\langle f\rangle\] of quadratic rational maps \[f\].

This action of \[{\rm PSL}_2\C\] is not free. For example the M\"obius
transformation \[g(z)=-z\] acts trivially on any odd function, such as \[f(z)
=a(z+z^{-1})\]. 
Hence we might expect the quotient space \[{\cal M}_2\] to have singularities.
In fact however, we will see that it has the simplest possible description,
and can be identified with the complex affine space \[\C^2\].
(On the other hand, since it is defined as a non-trivial quotient space,
\[\M\] does have a natural orbifold structure which reflects
the complications of the group action.
Compare \S5.)\smallskip

In order to describe this affine structure, let us study
fixed points. Every map\break \[f\in\Rat_2\]
has three not necessarily distinct fixed points \[z_1\,,\,z_2\,,\,z_3\in
S^2\]. Let \[\mu_i\] be the {\bit multiplier\/} of \[f\] at \[z_i\]
(that is the first derivative, suitably interpreted in the special case
when \[z_i\] is the point at infinity), and let
$$	\sigma_1=\mu_1+\mu_2+\mu_3\,,\quad\sigma_2=\mu_1\mu_2+
\mu_1\mu_3+\mu_2\mu_3\,,\quad\sigma_3=\mu_1\mu_2\mu_3 $$
be the elementary symmetric functions of these multipliers. (Note that \[\mu_i
=1\] if and only if \[z_i\] is a multiple fixed point, so that \[z_i=z_j\]
for some \[j\ne i\].)

{\QP{\bf Lemma 3.1.} \it These three multipliers determine \[f\] up to
holomorphic\break conjugacy, and are subject only to the restriction that
$$	\mu_1\mu_2\mu_3-(\mu_1+\mu_2+\mu_3)+2\;=\;0\;, \eqno (1)$$
or in other words
\vskip -.3in
$$	\sigma_3~=~\sigma_1\;-\;2~.\eqno (1') $$
Hence the moduli space \[{\cal M}_2\] is canonically isomorphic to \[\C^2\],
with coordinates \[\sigma_1\] and \[\sigma_2\].\medskip}

We will sometimes use the notation \[\langle f\rangle=\langle \mu_1\,,\,
\mu_2\,,\,\mu_3\rangle\] for the conjugacy class of a map \[f\] having
fixed points of multiplier \[\mu_1\,,\,\mu_2\] and \[\mu_3\].
If \[\mu_1\mu_2\ne 1\], then we can solve equation (1) for
$$	\mu_3\;=\;{2-\mu_1-\mu_2\over 1-\mu_1\mu_2}\;. \eqno (2) $$
On the other hand, if \[\mu_1\mu_2=1\] then it follows easily
from (1) that \[\mu_1=\mu_2=1\] so that \[z_1=z_2\] is a double fixed point.
In this case \[\mu_3\] can be arbitrary.\smallskip

{\bf Proof of Equation (1).} First suppose that
the \[\mu_i\] are all different from \[1$, so that there is no double fixed
point. Then the classical formula
$$	\sum 1/(1-\mu_i)\;=\;1  $$
is proved by integrating \[dz/(z-f(z))\]. (See
for example in [M2, \S9].) Clearing denominators, we obtain (1).
On the other hand, if \[\mu_1=1\] then \[z_1\] is a double fixed point, with
say \[z_1=z_2\] and \[\mu_1=\mu_2=1\]. The equation (1) is then true
for any value of \[\mu_3\].\QED

{\bf Proof that the holomorphic conjugacy class is
determined by \[\{\mu_1\,,\,\mu_2\,,\,\mu_3\}\].}
First consider a map \[f\] which has at least two
distinct fixed points. After
conjugating by a M\"obius transformation, we may assume that these two fixed
points are at zero and infinity. It follows easily that \[f\] has the form
$$	f(z)\;=\;z\,{a z+b\over c z+d}\;, $$
where \[a\ne 0\;,\quad d\ne 0\], and \[ad-bc\ne 0\]
since \[f\] has degree two. After multiplying numerator and denominator
by a constant, we may assume that \[d=1\].
If we replace \[f(z)\] by \[f(kz)/k\], the effect will be to multiply both
\[a\] and \[c\] by \[k\]. Thus there is a {\bit unique~}
choice of \[k\]
which has the effect of replacing \[a\] by \[1\].
This yields the normal form
$$	f(z)\;=\; z\,	{z+b\over c z+1}\qquad{\rm with}\qquad
1-b\,c\ne 0\;, \eqno (3) $$
where \[b=\mu_1\] and \[c=\mu_2\] are evidently
equal to the multipliers at zero and infinity respectively. Thus \[f\] is
uniquely determined, up to holomorphic conjugacy, by the multipliers
\[\mu_1\] and \[\mu_2\] associated with any two distinct fixed points.
Here the determinant \[1-\mu_1\mu_2\] cannot vanish, but there are no
other restrictions on \[\mu_1\] and \[\mu_2\]. The multiplier at the third
fixed point is then determined by Equation (2). For further information,
see Appendix C.

Now suppose there is only one fixed point.
After a M\"obius conjugation, we may assume that
this fixed point \[z_1=z_2=z_3\] is the point at infinity, and that \[f^{-1}
(\infty)=\{0,\infty\}\]. This implies that \[f\] has the form \[f(z)=p(z)/z\]
for some quadratic polynomial \[p(z)\]. Here the difference \[f(z)-z=
(p(z)-z^2)/z\] can have no zeros in the finite plane, hence \[p(z)-z^2\]
must be constant so that \[f(z)=z+c/z\], with critical points \[\pm\sqrt c\].
If we normalize so that the
critical points of \[f\] are \[\pm 1\], then \[c=1\] and
$$	f(z)~=~z+z^{-1}~. \eqno (4) $$
In this case, the multipliers at the unique fixed point are given by
\[\mu_1=\mu_2=\mu_3=1\], and again the conjugacy class is uniquely
determined by these multipliers.

Evidently we can realize any triple \[\{\mu_1\,,\,\mu_2\,,\,\mu_3\}\] which
satisfies Equation (1). Finally, note that the unordered
collection \[\{\mu_i\}\] of multipliers is determined by the three
elementary symmetric functions \[\sigma_n=\sigma_n(\mu_1\,,\,\mu_2\,,\,
\mu_3)\]. Since equation $(1')$ shows that
\[\sigma_3\] is determined by \[\sigma_1\], this completes the proof of 3.1.
\QED\smallskip

{\bf Remark 3.2: Cubic polynomial maps.} There is a strong analogy between
the theory of quadratic rational maps and of cubic polynomial maps. 
(Compare [M3], [M5].) In both
cases there are three fixed points and two critical points. Furthermore,
in both cases the moduli space of holomorphic conjugacy classes has
dimension two, and can be identified with \[\C^2\], with the elementary
symmetric functions of the multipliers at the fixed points, subject to a
single linear relation, as coordinates. In the cubic polynomial case, this
linear relation takes the form \[\sigma_2-2\sigma_1+3=0\].\medskip

{\bf Remark 3.3: Affine structure.}
Since the complex manifold \[\M\cong \C^2\] has many holomorphic
automorphisms, it is
not immediately clear that the affine structure imposed by taking the
\[\sigma_i\] as affine coordinates has any preferred status. However, the
following three lemmas show that this affine structure does indeed have very
special properties. (For
a different coordinate system, which would impose a different and less
useful affine structure, see 6.3.)\smallskip

{\bf Definition.} For each \[\eta\in\C\]
let \[\Per_1(\eta)\subset\M\] be the set of all conjugacy classes \[\langle f
\rangle\] of maps \[f\] which have a fixed point with multiplier equal to
\[\eta\]. (See Figure 1 for a plot of the \[\Per_1(\eta)\] in the real case.)

\midinsert
\centerline{\psfig{figure=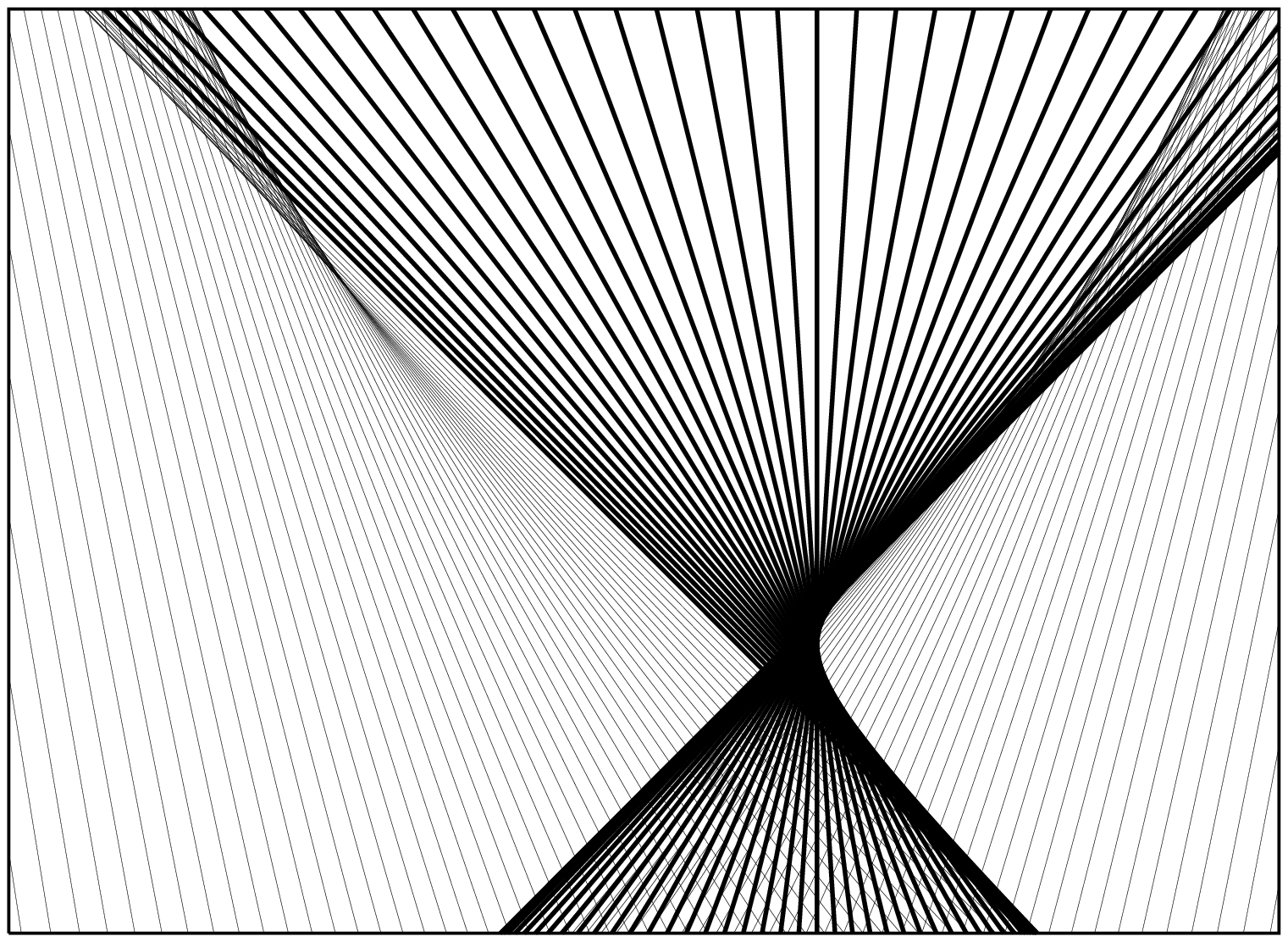,height=3.5in}}
{\narrower\smallskip\noindent
\bit Figure 1. The (real) lines \[\Per_1(\eta)\] in the real \[(\sigma_1\,,\,
\sigma_2)$-plane. (Compare \S10 and Figures 15, 16.)
The region \[(\sigma_1\,,\,\sigma_2)\in[-12,10]\times[-10,22]\] is shown.
(Horizontal scale exaggerated. Those lines with \[-1\le\eta\le 1\],
corresponding to attracting or parabolic real fixed points, have been
emphasized.) The envelope of this family \[\{\Per_1(\eta)\}\]
consists of the symmetry locus \[\sl\] of \S5,
together with the line \[\Per_1(1)\]. This envelope cuts the real
plane into three regions for which the
real quadratic map \[f\] has three distinct real fixed points, and two regions
for which \[f\] has only one real fixed point.\par}
\vskip -.1in
\endinsert

\eject
{\QP{\bf Lemma 3.4.} \it 
For each \[\eta\in\C\] this locus \[\Per_1(\eta)\subset\M\]
is a straight line with respect to the coordinates \[\sigma_1\,,\;\sigma_2\],
with slope \[~d\sigma_2/d\sigma_1=\eta+\eta^{-1}\]. For \[\eta\ne 0\] it is
given by the equation
$$	\sigma_2~=~(\eta+\eta^{-1})\,\sigma_1-(\eta^2+2\,\eta^{-1})~, \eqno (5) $$
while \[\Per_1(0)\] is the vertical line \[\sigma_1=+2\].\smallskip}

{\bf Proof.} The multipliers at the three fixed points are the roots of
the equation
$$	\eta^3-\sigma_1\,\eta^2+\sigma_2\,\eta-\sigma_3~=~0~. $$
Substituting \[\sigma_3=\sigma_1-2\] and solving for \[\sigma_2\],
we obtain the required equation.\QED

{\bf Definition.} More generally, for any integer \[n\ge 1\] and any number
\[\eta\ne 1\] in \[\C\],
let \[\Per_n(\eta)\] be the set of \[\langle f\rangle\in\M\]
having a periodic point of period \[n\] and multiplier \[\eta\]. (For the
special case \[\eta=1\], the definition needs more care. Compare Appendix D.
One possibility would be to simply define \[\Per_n(1)\] as the limit of
\[\Per_n(\eta)\] as \[\eta\to 1\;,~\eta\ne 1\].)
The following result will be proved in 4.2.

{\QP{\bf Lemma 3.5.} \it Each \[\Per_n(\eta)\] is an algebraic curve
in \[\M\] with degree equal to the number of period \[n\] hyperbolic components
in the Mandelbrot set.\smallskip}
\eject



Thus for \[n=1\] and \[n=2\] the curve \[\Per_n(\eta)\] is a straight line,
but for \[n=3\] it is a cubic curve. For period \[n=2\] we have the
following simple description.

{\QP{\bf Lemma 3.6.} \it The 
curves \[\Per_2(\eta)\] are parallel straight
lines of slope \[-2\], given by the equation
$$	2\sigma_1+\sigma_2~=~\eta~. $$
\smallskip}

As noted above, the case \[\eta=+1\] is exceptional. In fact
the proof will show that
there is no quadratic rational map having a period \[2\] orbit with
multiplier equal to \[+1\].\smallskip

{\bf Proof of 3.6.} Note that the fixed points of the \[4$-th degree map
\[f^{\circ 2}\] consist of the fixed points of \[f\] together with the period
\[2\] orbits (if any) of \[f\]. First consider
a map \[f\in\Rat_2\] with fixed point multipliers \[\{\mu_i\}\] such that no
\[\mu_i\] is
equal to \[\pm 1\]. Then we will show that the five fixed points of
\[f^{\circ 2}\] are all distinct. In fact, three of
these are the three distinct fixed points of \[f\]. These have
multipliers \[\mu_1^2\,,\,\mu_2^2\,,\,\mu_3^2\ne +1\] when considered as fixed
points of \[f^{\circ 2}\].
The remaining two must constitute a period two orbit for \[f\]. Neither
of these points can coincide with a fixed point of \[f\], since such a
multiple fixed point of \[f^{\circ 2}\]
would have to have multiplier \[+1\], and the two
cannot coincide with each other, since they would then constitute an extra
fixed point for \[f\]. It follows that the multiplier \[\eta\]
for this period two orbit cannot be \[+1\]. Hence
the rational fixed point formula for \[f^{\circ 2}\] takes the form
$$	{1\over 1-\mu_1^2}+{1\over 1-\mu_1^2}+{1\over 1-\mu_1^2}+{1\over
 1-\eta}+{1\over 1-\eta}~~=~~1~. $$
(Compare [M2] or the proof of 3.1.)
We can solve this equation, so as to express \[\eta\] as a certain
rational function of the elementary symmetric functions \[\sigma_i
=\sigma_i(\mu_1\,,\,\mu_2\,,\,\mu_3)\].
In fact, making use of the relation $(1')$ and carrying out the division
(preferably by computer),
we find the required formula \[\eta=2\sigma_1+\sigma_2\].

Now suppose that \[\langle f\rangle\] belongs to the locus \[\Per_1(-1)\],
or in other words suppose that \[f\] has
a fixed point of multiplier \[\mu_i=-1\]. Then 
\[f^{\circ 2}\] has a multiple fixed point, with multiplier
\[\mu_i^2=+1\]. If \[f\] has \[m\] distinct fixed points (where \[m=2\]
or \[m=3\]), then it follows easily that \[f^{\circ 2}\] has at most
\[m+1\] distinct fixed points. {\it Hence \[f\] cannot have
any period two orbit.\/} In fact, as \[\eta\to 1\], the unique
period two orbit for \[f\]
degenerates to the fixed point \[z_i\] of multiplier \[-1\].
(It follows from formula (1) that there cannot be two fixed points of
multiplier \[-1\].)
Note that the equation for the locus \[\Per_1(-1)\], as given by 3.4,
coincides precisely with the locus \[\eta=2\sigma_1+\sigma_2=+1\].
Thus it is convenient to define
$$	\Per_2(1)~=~\Per_1(-1)~=~\{\;\langle f\rangle \in\M\;:\;
 2\sigma_1+\sigma_2=1\;\}~. $$
(Compare Figure 6.)
Finally, suppose that \[f\] has a double fixed point \[~z_1=z_2~\] with
multiplier \[\mu_1=\mu_2=1\], and that the third fixed point \[z_3\] has
multiplier \[\mu_3\ne-1\]. Then a straightforward argument by continuity
shows that the formula \[\eta=2\sigma_1+\sigma_2\] for the multiplier
of the period two orbit remains true.
In this case, a brief computation shows that
the multiplier \[\eta=2\sigma_1+\sigma_2\] for the period two orbit is
equal to \[5+4\,\mu_3\ne+1\].\QED


\vfil\eject

\centerline{\bf \S4. The Compactification \[\widehat\M\cong\CP^2\].}
\smallskip

The coordinate plane \[\C^2\]
embeds naturally in the projective plane \[\CP^2\]. Since \[\M\]
is isomorphic to \[\C^2\]
with coordinates \[\sigma_1\] and \[\sigma_2\], there is a corresponding
compactification \[\widehat\M\cong\CP^2\], consisting of \[\M\] together with
a 2-sphere of {\bit ideal points\/} at infinity. Elements of this 2-sphere
can be thought of very roughly as limits of quadratic rational maps
which degenerate towards a fractional linear or constant map. However
caution is needed, since such a limit cannot be uniform over the entire
Riemann sphere.
\smallskip

In terms of the multipliers \[\{\mu_i\}\] at the fixed points, this 2-sphere
at infinity\break can be
described as follows. If at least one of the elementary
symmetric functions\break
\[\sigma_i=\sigma_i(\mu_1\,,\,\mu_2\,,\,\mu_3)\] tends
to infinity, then at least one of the \[\mu_i\] must tend to infinity.
If only \[\mu_3\], for example, tends to infinity, then it follows from
formula (2) that the product \[\mu_1\,\mu_2\] must tend to \[+1\].
On the other hand, if two of the \[\mu_i\] tend to infinity,
then using (2) we see that the third must tend to zero.
{\it Thus the collection of ideal points in \[\widehat\M\] can be identified
with the set of
unordered triples of the form \[\langle\mu\,,\,\mu^{-1}\,,\,\infty
\rangle\] with \[\mu\in\hat\C=\C\cup\infty\].\/} It seems appropriate
to use the notation \[\widehat\Per_1(\infty)\subset\widehat\M\]
for this 2-sphere of points at infinity.
A useful parameter on \[\widehat\Per_1(\infty)\]
is the sum \[\mu+\mu^{-1}\in\hat\C\], which can be identified
with the limiting ratio
$$	{\sigma_2\over\sigma_3}~=~{1\over\mu_1 }+{1\over\mu_2 }
	+{1\over\mu_3 }~. $$

If we exclude the special case \[\mu=1\], then the dynamics
of a representative map for a point in \[\M\] which is ``close'' to the
ideal point \[\langle\mu\,,\,\mu^{-1}\,,\,\infty\rangle\] can be described
as follows. (For \[\mu=1\], a useful description would
be more complicated, involving the theory of \'Ecale cylinders [La].)
As in \S3 (3) or Appendix C (22), use the normal form
$$	f(z)~=~z(z+\mu_1)/(\mu_2z+1)\qquad{\rm with}\qquad \mu_1\,\mu_2\ne 1
~. \eqno (6) $$
First suppose that \[\mu\ne 0\,,\,\infty\], and let
\[\mu_1\approx\mu\;,~\mu_2\approx\mu^{-1}\], hence
\[\mu_1\mu_2\approx 1\]. It turns out that, over
most of the Riemann sphere, this map \[f\] is uniformly close to the linear
map \[z\mapsto z/\mu_2\], or equivalently
\[z\mapsto\mu\,z\]. However, the behavior
is quite different in a small neighborhood of the point \[z=-1/\mu_2\]:
This neighborhood, which includes both critical points, maps over
the entire Riemann sphere. The case \[\mu=0\] is similar.
Any \[\langle f\rangle\in\M\] which is close to the ideal point
\[\langle 0,\infty,\infty\rangle\] has a convenient
representative which is uniformly
close to a constant throughout most of the sphere.
We will make these statements more precise, as
part of the proof of the following result.

{\QP{\bf Lemma 4.1.} \it For any period \[n\ge 2\] and for any multiplier
\[\eta\in\C\], the only possible
limit points of the curve \[\Per_n(\eta)\subset\M\]
on the 2-sphere at infinity are ideal points of the form
\[\langle\mu\,,\,\mu^{-1}\,,\,\infty\rangle\] where \[\mu\] is a \[q$-th
root of unity, with \[q\le n\].\smallskip}

In particular, the limiting ratio \[\sigma_2/\sigma_3=\mu+\mu^{-1}\]
is necessarily a point in the real interval \[[-2\,,\,2]\]. For example,
if \[\mu=e^{2\pi im/n}\] then \[\sigma_2/\sigma_3=2\cos(2\pi m/n)\].
(Compare Figures 16, 17.) It is conjectured that the case \[q=1\] cannot
occur, and
that the set of all limit points of \[\Per_n(\eta)\] is precisely the set of
\[\langle\mu\,,\,\mu^{-1}\,,\,\infty\rangle\] such that \[\mu\] is a
\[q$-th root  of unity with \[1<q\le n\].
\smallskip
\eject

{\bf Proof of 4.1.} First suppose that \[\mu\ne 0\,,\,1\,,\,\infty\].
As in the discussion above, we use the normal form
(6). Using the notations
$$\delta~=~1-\mu_1\,\mu_2~,\qquad \ell(z)~=~\mu_2z+1~, $$
let us assume that the determinant \[\delta\] is
very close to zero.
Note that the linear function \[\ell(z)\] will be close to zero if
and only if \[z\] is close to \[-1/\mu_2\]. A brief computation shows that
$$	\mu_2{f(z)\over z}~=~\mu_2{z+\mu_1\over \mu_2\,z+1}~=~1-{\delta
\over\ell(z)}~~,$$
and that
$$	\mu_2\,f'(z)~=~1-{\delta\over \ell(z)^2}~. $$
Let us partition the \[z$-plane into three non-overlapping
regions \[D\,,\,A\,,\,C\],
according as the number \[|\ell(z)|^3\] is less than \[|\delta|^2\], or between
\[|\delta|^2\] and \[|\delta|\], or greater than \[|\delta|\].
Thus \[D\] is a very small disk centered at the pole
\[-1/\mu_2\,\], \[A\] is a small annulus surrounding this disk, and the
complementary region \[C\]
is everything else, including zero and infinity. Note that:
$$	|\mu_2\,{f(z)\over z} -1|~\le~ |\delta|^{1/3}\qquad{\rm for}\qquad z\in A\cup
C~,$$
and
$$	|\mu_2\,f'(z)-1|~\le~|\delta|^{1/3}\qquad{\rm for}\qquad z\in C~, $$
but
$$	|\mu_2\,f'(z)-1|~\ge~|\delta|^{-1/3}\qquad{\rm for}\qquad z\in D~. $$
Thus \[f(z)/z\] is uniformly close to the constant \[1/\mu_2\approx\mu\]
everywhere outside of the small disk \[D\]. The derivative \[f'(z)\] is very
large throughout \[D\], and is uniformly
close to \[1/\mu_2\approx \mu\] and hence bounded away from zero throughout
the outside region \[C\]. It follows that
both critical points must belong to the annulus \[A\].

Now consider a periodic point of period \[n\ge 2\]. If its orbit is
disjoint from the disk \[D\], then \[1=f^{\circ n}(z)/z\approx
\mu^n\]. If the orbit touches both \[D\] and
\[A\], then we may assume that \[z\in A\], and that \[f^{\circ q}(z)\in D\]
for some \[1\le q<n\] which we take to be
minimal. In this case, it follows that \[\mu^q\approx
1\]. Finally, if an orbit touches \[D\] but not \[A\], then its multiplier
must tend to infinity as \[\delta\to 0\]. Thus, in the limit as
\[\delta\to 0\] and \[\mu_1\to\mu\], we can have an orbit
of period \[n\ge 2\] with bounded multiplier
only if \[\mu\] is a \[q$-th root of unity with \[q\le n\]. (Note that
this argument allows the possibility that \[\mu=1\].)\break
This completes the proof for the case \[\mu\ne 0\,,\,\infty\].

To handle the case \[\mu=\infty\] (or equivalently \[\mu=0\]) we need a
slightly different argument. Again use the normal form \[(6)\], but
now assume that the multiplier \[\mu_1\] at the origin is very large in
absolute value, and that the multiplier \[\mu_2\] at infinity is very
close to zero. It then follows from formula (2) that the multiplier \[\mu_3\]
at the third fixed point \[z_3=(\mu_1-1)/(\mu_2-1)\approx -\mu_1\] is also
very large in absolute value. We write \[\mu_1\,,\,\mu_3\approx\infty\] but
\[\mu_2\approx 0\]. For \[z\] in the disk \[|z|<2\], it then follows 
from the computation
$$	f'(z)~=~ {\mu_2\,z^2+2z+\mu_1\over (\mu_2\,z+1)^2} $$
that the derivative \[f'(z)\] is uniformly close to \[\mu_1\]. Hence this
disk maps diffeomorphically onto a region \[U\], which is approximately the
disk of radius \[2|\mu_1|\] enclosing both finite fixed points. Let \[D\subset
U\] be the disk centered at the midpoint of the two finite fixed points,
with radius equal to the distance between them. Since the two
finite fixed points play a symmetric role, it follows that the pre-image
\[f^{-1}(D)\subset D\]
splits up as a neighborhood \[N_1\] of \[z_1=0\], throughout which
\[f'\approx\mu_1\],
and a neighborhood \[N_3\] of \[z_3\] throughout which \[f'\approx\mu_3\].
It is now easy to check that the Julia set \[J(f)\]
is a Cantor set, contained in the union \[N_1\cup N_3\], and that every
orbit outside of the Julia set converges to the fixed point at infinity.
(Compare \S8 below.) Thus \[f'(z)\] is approximately equal to either
\[\mu_1\] or \[\mu_3\] throughout the
Julia set; hence the multiplier of any orbit of period \[\ge 2\] tends
to infinity as \[\mu_1\,,\,\mu_3\to\infty\].\QED

{\QP{\bf Theorem 4.2.} \it For any \[\mu\in\C\],
the degree of the curve \[\Per_n(\mu)\subset\M\]
is equal \[\nu_2(n)/2\], where the numbers \[\nu_2(n)\]
are defined inductively by the formula
$$	2^n~=~\sum_{m|n}\nu_2(m)~, $$
to be summed over all positive integers \[m\] which divide \[n\]. Equivalently,
this degree is equal to the number of hyperbolic components of period \[n\]
in the Mandelbrot set. 
\smallskip}

(Compare 3.5 and Appendix D.) Here is a table giving some examples.
\medskip

\centerline{\vbox{\hrule\hbox{\vrule\vbox{
{\halign{\quad # & \quad# & # & # & # & # & # & # & #\quad\cr
\noalign{\smallskip}
$n$ & 1 & 2 & 3 & 4 & 5 & 6 & 7 & 8\cr
degree & 1 & 1 & 3 & 6 & 15 & 27 & 63 & 120\cr
\noalign{\smallskip}}
}}\vrule}\hrule}}\smallskip


{\bf Proof of 4.2.} By 3.4, it suffices to consider the case \[n>1\].
Since the definition of \[\Per_n(\mu)\] is purely
algebraic, it is not difficult to check that it is an algebraic curve in
\[\M\cong\C^2\]. (See Appendix D.)
In fact it is convenient to consider its closure \[\widehat\Per_n(\mu)\]
in the
projective space \[\widehat\M\cong\CP^2\]. By definition, the degree of a
curve in \[\CP^2\] is equal to its number of intersections with any straight
line, counted with multiplicity.
As line in \[\M\], we choose the closure of the
locus \[\Per_1(0)\], with equation \[\sigma_1
=2\]\break (or equivalently \[\sigma_3=\mu_1\mu_2\mu_3=0\]). This line can be
identified with the set of all quadratic polynomial maps \[f(z)=z^2+c\],
having a fixed point of multiplier zero at infinity. (Here
\[\sigma_2=4c\].) The closure \[\widehat\Per_1(0)\] within the compactified
space \[\widehat\M\] contains just one
point at infinity \[\langle 0,\infty,\infty\rangle\].
According to 4.1, the curve \[\widehat\Per_n(\mu)\] does
not contain this point.
Thus it suffices to consider intersections in the finite plane. First consider
the case \[\mu=0\]. The points of the intersection \[\Per_1(0)\cap\Per_n(0)\]
can be described as the conjugacy classes of maps \[f_c(z)= z^2+c\] for
which the finite critical point \[0\] has period exactly \[n\]
under \[f_c\]. By definition, these are exactly the center points of the
various period \[n\] components of the Mandelbrot set. Furthermore,
it follows easily from [DH1, \S3] that
the multiplicity of such a value \[c\] as solution to the equation
\[f_c^{\circ n}(0)=0\] is always \[1\]. In other words, the intersection
is always transverse, with intersection multiplicity \[1\]. More
generally, whenever \[|\mu|<1\]
a similar argument shows that the curves \[\Per_1(0)\]
and \[\Per_n(\mu)\] intersect transversally, with exactly one
intersection point in each period \[n\] component of the Mandelbrot set.

To identify this number of intersection points with \[\nu_2(n)/2\],
note that the equation \[f_c^{\circ n}(0)=0\] has degree \[2^n/2\]. In other
words, the number of centers in the Mandelbrot set with period dividing \[n\]
is equal to \[2^n/2\]. After discarding all of those centers corresponding
to proper divisors of \[n\], we obtain the required number, namely
\[\nu_2(n)/2\].

For the general case, it is convenient to introduce algebraic curves
\[Q_n^*\supset Q_n\to P_n\] as follows. Let \[Q_n^*\subset\C^2\] be the set
of all pairs \[(c\,,\,z)\] satisfying the polynomial equation \[f_c^{\circ n}
(z)=z\]. (Thus the point \[z\] must be periodic with period \[m\]
dividing \[n\] under the map \[f_c\].) Let \[Q_n\subset Q_n^*\] be the union
of those irreducible components of the curve \[Q_n^*\] for which a generic
point \[(c_0\,,\,z_0)\] has the property that \[z_0\] has period exactly \[n\]
under \[f_{c_0}\]. (According to Bousch [Bou], there is exactly
one such irreducible component; in other words, \[Q_n\] is irreducible.)
Evidently we can write
$$	Q_n^*~=~\bigcup_{m|n} Q_m~, $$
taking the union over all divisors \[m\] of \[n\]. The cyclic group of order
\[n\] operates on \[Q_n\] by the transformation \[(c\,,\,z)\mapsto
(c\,,\,f_c(z))\]. Let \[P_n\] be the quotient variety of \[Q_n\] under this
action. Thus a point of \[P_n\] can be described as a pair \[(c\,,\,\{z_i\})\]
consisting of a parameter value \[c\] and a periodic orbit \[z_1\mapsto z_2
\mapsto\cdots\] under \[f_c\] which (at least in the generic case) has period
exactly \[n\].

For each fixed value of \[z\], note that the defining equation \[f_c^{\circ n}
(z)=z\] has degree \[2^n/2\] in \[c\]. In other words, the projection map
\[(c,z)\mapsto z\] from \[Q_n^*\] to the \[z$-plane has degree \[2^n/2\].
If we restrict to the subvariety \[Q_n\], it follows easily that the
corresponding projection map \[(c,z)\mapsto z\]
has degree \[\nu_2(n)/2\]. If \[z=z_1\mapsto z_2\mapsto
\cdots\] is the orbit of \[z\] under \[f_c\], then it follows
that each projection \[(c,z)\mapsto z_i\] from \[Q_n\] to \[\C\] also
has degree \[\nu_2(n)/2\]. Now
note that the multiplier \[\eta=(2z_1)(2z_2)\cdots(2z_n)\] of such a periodic
orbit is, up to a constant factor, just the product \[z_1\,z_2\cdots z_n\].
It follows easily that the projection \[(c,z)\mapsto\eta\] from \[Q_n\]
to the \[\eta$-plane has degree \[\sum_1^n\nu_2(n)/2=n\,\nu_2(n)/2\].
For example this can be proved by considering 2-dimensional cohomology
with compact support for the composition
$$	Q_n\ssm\eta^{-1}(0)~\to~\prod_1^n\C\ssm\{0\}~\buildrel{\rm product}
	\over\longrightarrow~\C\ssm\{0\} $$
of proper maps. Finally, since the
projection \[Q_n\to P_n\] has degree \[n\], this implies that the
projection \[(c\,,\,\{z_i\})\mapsto\eta\] from \[P_n\] to the \[\eta$-plane
has the required degree \[\nu_2(n)/2\]. In other words, for generic
choice of \[\eta\] there are \[\nu_2(n)/2\] corresponding points in the curve
\[P_n\], which map to
\[\nu_2(n)/2\] distinct points of the \[c$-plane. Of course for particular
values of \[\eta\] there may be coincidences, but this will not affect the
count with multiplicity. Now using the argument above,
it follows that the curve \[\Per_n(\eta)\subset\M\]
has degree \[\nu_2(n)/2\]. This proves 4.2 and 3.5.\QED


{\bf Remark 4.3. The curves \[\Per_3(\eta)\].} As an example (without proofs),
let us look at the special case \[n=3\]. It follows
from 4.2 that each \[\Per_3(\eta)\] is a curve of degree three.
For most values of \[\eta\], this curve is non-singular of genus one, and has
two ends corresponding to the two intersection
points of its projective completion with the line at infinity (namely
a double intersection point at \[\langle \omega,\bar\omega,\infty
\rangle\] where \[\omega\] is a primitive cube root of unity, and
a single intersection point at \[\langle -1,-1,\infty\rangle\]). However,
there are three special values of \[\eta\] which behave differently.
For \[\eta=0\], the curve
\[\Per_3(0)\] has genus zero, with a transverse self-intersection point
corresponding to the map for which both critical points lie in a single period
\[3\] orbit (Figures 2, 9). Similarly, for
\[\eta=-8\] the curve \[\Per_3(-8)\] has genus zero,
with a single transverse self-intersection point corresponding to the map
\[z\mapsto 1/z^2\], which has
two distinct period \[3\] orbits with multiplier \[-8\].
Finally, for \[\eta=1\] the cubic curve \[\Per_3(\eta)\]
degenerates into a union of three straight lines:
$$	\Per_3(1)~=~\Per_1(\omega)\cup\Per_1(\bar\omega)\cup\Per_2(-3)~,
\eqno (7) $$
where \[\omega\] is a primitive cube root of unity. (Compare Figures 8, 10.)
The first two lines, with slope \[-1\],
correspond to maps for which one period \[3\] orbit degenerates to a fixed
point of mutiplier \[\omega\] or \[\bar\omega\], while the third straight
line, with equation \[\sigma_2=-2\sigma_1-3\],
corresponds to maps for which the two period \[3\] orbits coincide. 
For some reason, which I do not understand, this locus is precisely equal
to the line \[\Per_2(-3)\].
This third line is visible as part of the boundary of a hyperbolic
component in Figure 16. 
Thus the curve \[\Per_3(1)\] has two finite self-intersections (corresponding
to the map
\[z\mapsto\omega(z+z^{-1})\] and its complex conjugate),
and one self-intersection at infinity.

\bigskip \smallskip

\centerline{\bf \S5. Maps with Symmetries.}\medskip

By an {\bit automorphism\/} of a quadratic rational map \[f\], we will
mean a M\"obius transformation \[g\] which commutes with \[f\], so
that \[g\circ f\circ g^{-1}=f\]. The collection of all automorphisms
of \[f\] forms a finite group
$$	\Aut(f)~\subset~\Rat_1~\cong~\PSL(2,\C)~, $$
which measures the extent to which the action of \[\Rat_1\] on \[\Rat_2\]
by conjugation fails to be free at \[f\].

{\QP{\bf Theorem 5.1.} \it A quadratic rational map possesses a
non-trivial automorphism if and only if it is conjugate to a
map in the unique normal form
$$      f(z)~=~k\,(z+z^{-1})~, \eqno (8) $$
with \[k\in\C\ssm\{0\}\].
For \[f\] in this normal form, if \[k\ne-1/2\] the group \[\Aut(f)\]
is cyclic of order two,
consisting of the maps \[z\mapsto\pm z\]. However, for
\[k=-1/2\] the group \[\Aut(f)\] is non-abelian of order six\footnote
{$^1$}{\rm Compare  [DM], [Mc]}.\medskip}

See Figure 12 for a picture of the \[k$-plane.\smallskip

{\bf Remark 5.2.} For a map in this normal form, note that the point at
infinity is a fixed point with multiplier \[\mu=1/k\]. There are two
other fixed points at \[z=\pm\sqrt{k/(1-k)}\], both with multiplier \[2k-1\].
Thus the fixed point multipliers of \[f\] are
$$	\{\mu_i\}~=~\{k^{-1}\,,\,2k-1\,,\,2k-1\}~.	\eqno (9) $$
There are two special values of \[k\] for which all three multipliers
are equal. These are the exceptional point \[k=-1/2\] of 5.1,
with \[\mu_1=\mu_2=\mu_3=-2\], and the point 
\[k=1\] with \[\mu_1=\mu_2=\mu_3=1\]. In the latter case, all three fixed
points coincide with the point at infinity, as discussed in the proof
of 3.1.\smallskip


{\bf Proof of 5.1.} First consider an automorphism which has order
two. Any element of order two in \[\PSL(2\,,\,\C)\] is conjugate to the map
\[g(z)=-z\], so it suffices to look at quadratic rational maps \[f\] which
commute with  \[z\mapsto-z\]. In other words, it suffices to look at
odd functions, \[f(-z)=-f(z)\].
Writing \[f(z)\] as a quotient \[p(z)/q(z)\] of two polynomials, we see easily
that \[f\] is odd if and only if one of these two polynomials is odd and the
other is even. If \[p(z)\] is even and \[q(z)\] is odd, then we can write
$$	f(z)~=~ {k\,z^2+\ell\over z}~=~k\,z+\ell\,z^{-1}\qquad{\rm with}
	\qquad k\ell\ne 0~. $$
Choosing \[\lambda=\pm\sqrt{\ell/k}\], we see that the
conjugate map \[f(\lambda\,z)/\lambda\] has the required
form \[z\mapsto k\,(z+z^{-1})\]. On the other hand, if \[p(z)\] is odd and
\[q(z)\] is even, then the conjugate map \[1/f(1/z)\] will have the form
{\it even$/$odd~}, so that the above argument applies.\smallskip

For \[f\] in this normal form (8), note that the two critical points
\[\pm 1\] are interchanged by the automorphism \[z\mapsto -z\]. In fact for any
quadratic rational map \[f\] and automorphism \[g\] it is clear that the set
of critical points \[\{\omega_1\,,\,\omega_2\}\] must be mapped into itself
by \[g\]. Hence the automorphism group \[\Aut(f)\] contains a subgroup
\[\Aut^0(f)\] of index \[\le 2\] consisting of automorphisms which fix
each of the two critical points. Suppose that this subgroup contains a
non-trivial automorphism \[g\]. After a M\"obius change of coordinates, we may
assume that the two critical points are zero and infinity. Thus the
non-trivial automorphism
\[g\] fixing these two points must have the form \[g(z)=\lambda\,z\] for some
\[\lambda\ne 0\,,\,1\]. The
equation \[f(\lambda\,z)=\lambda\,f(z)\] then implies that \[f(0)\in\{0\,,\,
\infty\}\], and similarly that\[f(\infty)\in\{0\,,\, \infty\}\]. If \[f\]
also fixed both critical points, then evidently \[f(z)=\alpha\,z^2\] for some
constant \[\alpha\ne 0\], and the equation \[f(\lambda\,z)=\lambda\,f(z)\]
would imply that \[\lambda=1\], contrary to our hypothesis. Since the two
critical values must be distinct, the only other possibility is that \[f\]
interchanges the two critical points. Thus \[f\] must have the form
\[f(z)=\alpha/z^2\], and
after a scale change we may assume that \[\alpha=1\] so that
$$	f(z)~=~1/z^2~.	\eqno (10) $$
A brief computation then shows that the group \[\Aut^0(f)\] of automorphisms
which fix zero and infinity is the cyclic group of order three, consisting
of all maps \[g(z)=\lambda\,z\] with \[\lambda^3=1\]. The full group of
automorphisms for this map (10) is generated by this subgroup, together
with the involution \[z\mapsto 1/z\]. Making use of the discussion above, or by
direct computation, we see that the map (10) is holomorphically
conjugate to the  special case \[z\mapsto -(z+z^{-1})/2\] of formula (8).
Further details of the proof are straightforward.\QED

In terms of the fixed points of \[f\], we can reformulate this result as
follows:\smallskip

{\bf Case 1.} If \[f\] has three distinct fixed points (or in other words
if \[\mu_i\ne 1\]), then \[\Aut(f)\] coincides with the group
consisting of all permutations of the fixed points which preserve the
multiplier. Thus \[\Aut(f)\] has order \[1\,,~2\] or \[6\] according as
the \[\mu_i\] are distinct, two are equal, or all three are equal.\smallskip

{\bf Case 2.} If \[f\] has only two distinct fixed points, then \[\Aut(f)\]
is trivial.\medskip

{\bf Case 3.} If \[f\] has only one fixed point, then \[\Aut(f)\] is cyclic
of order two.\smallskip

{\bf Definition.} Let \[\sl\subset\M\] be the {\bit symmetry locus\/},
consisting
of all conjugacy classes \[\langle f\rangle\] of quadratic maps which possess a
non-trivial automorphism. (For other characterizations of \[\sl\], see 5.4
and 6.4.) Using formula (9), we easily prove the following.

{\QP{\bf Corollary 5.3.} \it The symmetry locus \[\sl\] is a curve of degree
three and genus zero in \[\M\cong\C^2\]. It can be defined parametrically
by the equations
$$	\sigma_1=4\,k-2+k^{-1}\,,\qquad \sigma_2=4\,k^2-4\,k+5-2\,k^{-1}~, $$
as \[k\] varies over \[\C\ssm\{0\}\]. This curve is non-singular, except
for a cusp 
at the point \[\langle z\mapsto 1/z^2\rangle\], with \[k=-1/2\,,~
\sigma_1=-6\,,~\sigma_2=12\].
\medskip}

Compare Figure 15, which shows the intersection of \[\sl\] with the real
\[(\sigma_1\,,\,\sigma_2)$-plane, and see also Figure 1.\smallskip

{\bf Remark 5.4. Orbifold structure.} Since the action
\[~g,\,f\mapsto g\circ f\circ g^{-1}~\]
of the group \[\PSL(2,\C)\] on the space \[\Rat_2\] is proper and locally
free, it follows that the quotient space \[\M\]
has an associated orbifold structure. In fact, if
\[U\subset\Rat_2\] is a complex
\[2$-manifold transverse to the orbit \[\{g\circ f_0\circ g^{-1}\}\], then
the finite group \[\Aut(f_0)\] acts on \[U\] in a neighborhood \[U_0\]
of \[f_0\],
and the quotient of \[U_0\] by this action is precisely a coordinate
neighborhood of \[\langle f_0\rangle\] in the orbifold
\[\M\]. {\it Evidently, we can
describe the symmetry locus \[\sl\] as the set of points
of \[\M\] at which this natural orbifold structure is
non-trivial\/}. Note that this structure is particularly
non-trivial at the cusp point \[\langle z\mapsto 1/z^2\rangle\].

\bigskip\bigskip

\centerline{\bf \S6. Maps with Marked Critical Points or Fixed Points.
}\medskip

Recall from \S2 that a {\bit critically
marked\/} quadratic rational map \[(f\,,\,\omega_1\,,\,\omega_2)\] is
a map \[f\in\Rat_2\] together with an ordered list of its critical points.
The space \[\Rat_2^{\rm cm}\] of all critically marked quadratic rational
maps is
a smooth two-sheeted covering manifold of \[\Rat_2\]. The M\"obius group
\[\Rat_1\cong\PSL(2,\C)\] acts on \[\Rat_2^{\rm cm}\] by conjugation
$$	g\cdot(f\,,\,\omega_1\,,\,\omega_2)~=~\big(g\circ f\circ g^{-1}\,,\,
	g(\omega_1)\,,\,g(\omega_2)\big)~. $$
The quotient space of \[\Rat_2^{\rm cm}\] under this action will be denoted
by \[{\cal M}_2^{\rm cm}\], and called the {\bit critically marked moduli
space\/}. Following Rees [R3], we will show that
this moduli space is a smooth complex
manifold except at one singular point, corresponding to the special map
$$	f(z)~=~1/z^2~. \eqno (10) $$
To understand this space \[\M\cm\]
it is convenient to use the normal form \[(f\,,\,0\,,\,\infty)\] where
$$	f(z)\;=\; {\alpha z^2+\beta\over \gamma z^2+\delta} \qquad{\rm with}
\qquad \alpha\delta-\beta\gamma=1\; , \eqno (11) $$
so that the two marked critical points are zero and infinity respectively.
(Compare\break Appendix C.)
Note that maps in this form satisfy \[f(z)=f(z')\] if and only if \[z'=\pm z\].
It follows that the Julia set \[J(f)\] is invariant under the involution
\[z\leftrightarrow-z\].

This normal form
is unique except for the scale change which replaces \[f(z)\] by \[f(\lambda^2
z)/\lambda^2\]. This acts on the unimodular matrix of coefficients by the
transformation
$$	\left [\matrix{\alpha& \beta\cr \gamma& \delta}\right]\;\mapsto
	\left [\matrix{\alpha\lambda& \beta/\lambda^3\cr \gamma\lambda^3&
 \delta/\lambda}\right]\qquad{\rm for\;any}\quad\lambda\,\in\,\C\ssm\{0\}\;. 
 \eqno (12) $$
(Note that we can change the signs of all coefficients by taking \[\lambda=
-1\].) It is easy to check that
this action of the group \[\C\ssm\{0\}\] on the manifold
\[{\rm SL}(2,\C)\] of complex unimodular matrices
is free with a single exception: The cube roots of unity
act trivially on the orbit \[\alpha=\delta=0\], which corresponds to
the special mapping (10).\smallskip

We can introduce three expressions which are invariant under this action
(12) by the formulas
$$	A=\alpha\delta=1+\beta\gamma\;,\qquad B=\alpha^3\beta\;,\qquad
\quad C=\gamma\delta^3\;. \eqno (13) $$
It seems difficult to interpret these quantities geometrically, but they are
quite convenient to work with.

{\QP{\bf Lemma 6.1.}
\it The moduli space \[{\cal M}_2^{\rm cm}\] for critically marked quadratic
rational maps can be identified with
the hypersurface \[W\] consisting of all triples
\[(A,B,C)\in\C^3\] which satisfy
the equation
$$	A^3(A-1)\;=\; BC\; .\eqno (14) $$
This algebraic surface
is non-singular except at the point \[A=B=C=0\],\break corresponding to
the mapping $(10)$, where it has an essential singularity.
The deck transformation \[(f\,,\,\omega_1\,,\,\omega_2)\leftrightarrow
(f\,,\,\omega_2\,,\,\omega_1)\] of \[{\cal M}_2^{\rm cm}\] over
\[{\cal M}_2\], which interchanges the two critical points, 
corresponds to the map
$$(A,B,C)~\leftrightarrow~(A,C,B)~.$$
\smallskip}

{\bf Proof.} It is clear that these quantities \[A,\,B,\,C\]
are indeed invariant under (12), and that they satisfy the relation (14).
Conversely, given \[A,\,B,\,C\] satisfying (14) then we can find
\[\alpha,\,\beta,\,\gamma,\,\delta\] satisfying (13),
unique up to the action of \[\C\ssm\{0\}\], as follows.
If \[B\ne 0\] or \[A\ne 0\], then we can set \[\alpha=1\], and solve
uniquely for \[\beta=B\,,\; \delta=A\] and either \[\gamma=(A-1)/B\] or
\[\gamma=C/A^3\]. The case \[C\ne 0\] is similar, and the case \[A=B=C=0\]
reduces to the mapping (10).
Finally, note that the conjugacy which replaces \[f(z)\] by
$$ {1\over f(1/z)}\;=\;{\delta z+\gamma\over\beta z+\alpha} $$
interchanges the roles of the two critical points, and interchanges the
invariants \[B\,,\,C\].\QED

{\bf Remark 6.2: Singularities.} In
order to understand the singularity of the hypersurface
(14) at the origin, it is convenient to make the following.
{\bf Definition.} A surface in \[\C^3\]
has a singularity of {\bit type~} \[(p,q,r)\] if it can be reduced to the
form  \[z_1^p+z_2^q+z_3^r=0\] by a local holomorphic change of variable,
where \[p\,,\,q\,,\,r>1\]. Such a point is indeed always singular. In fact
the surface is locally homeomorphic to the cone over a 3-manifold with
non-trivial fundamental group. (See for example [M1].) In particular,
for a singularity of type \[(2,2,r)\] this fundamental group is cyclic
of order \[r\]. Similarly, we will say that a curve in \[\C^2\] has a
singularity of {\bit type~} \[(p,q)\] if it can be reduced to the form
\[z_1^p+z_2^q=0\].

Clearly the hypersurface (14) has a singularity of type \[(2,2,3)\] at
the origin. Hence a neighborhood of the origin is homeomorphic to the
cone over a 3-dimensional lens space which has fundamental group of order
3. We can resolve this singularity locally by passing to a 3-sheeted
covering space which is ramified at this single singular point.
(Compare 6.6.)\smallskip

{\bf Remark 6.3: \[{\cal M}_2\cong\C^2\].}
Lemma 6.1 provides a quite different proof that the
moduli space \[{\cal M}_2\] is isomorphic to \[\C^2\]. Evidently we can
obtain \[{\cal M}_2\] from the algebraic surface (14) by identifying
each triple \[(A,B,C)\] with \[(A,C,B)\]. Let us introduce the sum
\[\Sigma=B+C\], which is invariant under this involution. Given any
pair \[(A,\Sigma)
\in\C^2\], we can solve the equations \[A^3(A-1)=BC\] and \[\Sigma=B+C\]
uniquely for the unordered pair \[\{B,C\}\]. This the quotient surface
is isomorphic to \[\C^2\], with coordinates \[A\] and \[\Sigma\].
(However, these new coordinates are not compatible with the compactification
introduced in \S4.)
Of course this proof immediately raises a question: How are these new
coordinates \[(A,\,\Sigma)\] related to the coordinates \[(\sigma_1\,
,\,\sigma_2)\]
of \S3? This question will be answered in Appendix C.\smallskip

{\bf Remark 6.4: \[{\cal M}_2^{\rm cm}\] as 2-sheeted covering.} Evidently
the critically marked moduli space \[{\cal M}_2^{\rm cm}\cong W\]
can be considered as a 2-sheeted ramified covering space of \[{\cal M}_2
\cong\C^2\]. Evidently the covering map is ramified precisely over the
symmetry locus \[\sl\] of \S5. For there exists an automorphism of \[f\]
interchanging the two critical points if and only if \[\langle f\rangle
\in\sl\].
The ramification locus or symmetry locus corresponds to the set of \[(A,B,C)\]
in the hypersurface \[W\]
which satisfy \[B=C=\Sigma/2\], and hence are invariant
under the involution \[B\leftrightarrow C\]. In terms of the coordinates
\[(A\,,\,\Sigma)\] on \[{\cal M}_2\], this locus can be described
by the \[4$-th degree\footnote{$^1$}
{Thus \[\sl\] is a \[4$-th degree curve in terms of the coordinates
\[(A,\,\Sigma)\] for \[\M\cong\C^2\], but a cubic curve in terms
of the more natural coordinates \[(\sigma_1\,,\,\sigma_2)\].}
equation
$$	4A^3(A-1)~=~\Sigma^2	~. \eqno (15) $$
As noted already in \S5, this locus \[\sl\]
can be described geometrically as a curve of genus zero in \[\C^2\]
with a single cusp point at the origin.\smallskip

{\bf Remark 6.5: Homotopy type.} This critically marked moduli
space \[{\cal M}_2^{\rm cm}\] has the
homotopy type of the 2-sphere. In fact the correspondence
$$	\left[\matrix{\alpha&\beta\cr\gamma&\delta}\right]\;\mapsto\;
{\alpha^3\over\gamma}\;=\;{B\over A-1}\;=\;{A^3\over C} $$
is a smooth map from \[{\cal M}_2^{\rm cm}\] onto the Riemann sphere
\[S^2=\C\cup\infty\] with the property
that the inverse image of any point is isomorphic to \[\C\]. (This map is of
course not a fibration: local triviality fails about the singular point
\[A=B=C=0\], which maps to \[0\].) The topological 2-sphere
$$	0\le A\le 1\,,\quad |B|=\sqrt{A^3(1-A)}\, ,\quad C=-\bar B $$
is embedded in  \[{\cal M}_2^{\rm cm}\] as a deformation retract, and
maps homeomorphically onto \[\bar\C\]. (This 2-sphere provides a natural
example of a ``teardrop orbifold'', which is simply\-connected but has one
point with non-trivial orbifold structure.)
\eject

\centerline{\bf Marked Fixed Points.}\smallskip

Instead of numbering the critical points of a rational map \[f\],
we can equally well number the fixed points.
Let \[\Rat_2^{\rm fm}\] be the space of {\bit fixed point marked\/} quadratic
rational maps, that is ordered 4-tuples \[(f\,,\,z_1\,,\,z_2\,,\,z_3)\] where
\[z_1\,,\,z_2\,,\,z_3\in S^2\] are the fixed points
of \[f\]. Here a double (or triple) fixed point is to be listed twice
(respectively three times). Let \[{\cal M}_2^{\rm fm}\]
be the quotient of this space by the group \[\Rat_1\], acting by conjugation.
Indeed we can go further and mark {\it both\/} the fixed points and the
critical points, thus producing a space \[\Rat_2^{\rm tm}\] of {\bit totally
marked\/} rational maps \[(f,z_1,z_2,z_3,\omega_1,\omega_2)\], with
quotient moduli space \[{\cal M}_2^{\rm tm}\].

{\QP{\bf Lemma 6.6.} \it The space \[\Rat_2^{\rm fm}\] is a smooth complex
5-manifold, and \[\Rat_2^{\rm tm}\] is an unramified 2-fold covering
manifold of it. The action of \[\Rat_1\] on \[\Rat_2\tm\] by conjugation is
free, so that the quotient space \[{\cal M}_2^{\rm tm}\] under this action
is a smooth complex 2-manifold. However, the action of \[\Rat_1\] on
\[\Rat_2\fm\] has one non-free orbit, hence the moduli space
\[{\cal M}^{\rm fm}\] has
one singular point. It corresponds to the conjugacy
class of the map \[z\mapsto z+z^{-1}\] which has just one triple fixed point.
\medskip}

{\bf Proof.}
To see that the space \[\Rat_2^{\rm tm}\] of totally marked maps
is a smooth complex 5-manifold, we will
show that it can be identified with an open subset of the product
\[S^2\times  S^2\times S^2\times S^2\times S^2\]. More precisely, we
show that a point
$$	(f\,,\,z_1\,,\,z_2\,,\,z_3\,,\,\omega_1\,,\,\omega_2)
~\in~\Rat_2^{\rm tm} $$
is uniquely determined by its 5-tuple
\[(z_1\,,\,z_2\,,\,z_3\,,\,\omega_1\,,\,\omega_2)\] of fixed and critical
points, and that a given 5-tuple actually occurs if and only if:

(a) \[\omega_1\ne\omega_2\], and

(b) for each \[i\ne j\] the cross-ratio $${(z_i-\omega_1)\,(z_j-\omega_2)
\over (z_j-\omega_1)\,(z_i-\omega_2)}~\in~\C\cup\infty$$ 
is well defined and different from \[-1\].

As in (11), it is convenient to consider the special case \[\omega_1=0\,,
\;\omega_2=\infty\], so that \[f(z)=(\alpha z^2+\beta)/(\gamma z^2+\delta)\].
With this choice of critical points, Condition (b) simply says
that \[z_i\ne -z_j\]. In fact the two points \[z_i\]
and \[-z_i\] cannot be distinct and both fixed
since they have the same image under \[f\], and we cannot have \[z_i=z_j\in\{
0,\infty\}\], since a critical fixed point cannot also be a double fixed point.
The fixed points of \[f\] are the roots of the equation \[\gamma z^3-\alpha z^2
+\delta z-\beta=0\]. (Here, as usual, a fixed
point at infinity corresponds to a polynomial equation of reduced degree.)
Evidently the set of roots determines and is determined by
the point \[(\alpha:\beta:\gamma:\delta)\in\CP^3\], which is subject only
to the determinant
inequality \[\alpha\delta-\beta\gamma\ne 0\]. Expressing these
coefficients in terms of the fixed points, a brief computation show that
this determinant inequality is equivalent to Condition (b).

Thus \[\Rat_2^{\rm tm}\] is a smooth complex 5-manifold.
Since \[\Rat_2^{\rm fm}\] is the quotient space of \[\Rat_2^{\rm
tm}\] by the fixed point free holomorphic involution
$$	(z_1\,,\,z_2\,,\,z_3\,,\,\omega_1\,,\,\omega_2)~\leftrightarrow~
	(z_1\,,\,z_2\,,\,z_3\,,\,\omega_2\,,\,\omega_1)~, $$
it follows that
\[\Rat_2^{\rm fm}\] is also a smooth complex manifold. Similarly, since
the 5-tuple \[(z_1\,,\,z_2\,,\,z_3\,,\,\omega_1\,,\,\omega_2)\] must
contain at least three distinct points, the action of \[\Rat_1\] on \[\Rat_2^{
\rm tm}\] by conjugation is necessarily free, and it follows that
the quotient space \[{\cal M}_2^{\rm tm}\] is a smooth complex manifold.

On the other hand, the action of \[\Rat_1\] on \[\Rat_2^{\rm fm}\] is
not free: If \[f(z)=z+z^{-1}\] with just one triple fixed point
at infinity, then the involution \[z\mapsto -z\] in \[\Rat_1\] acts
trivially on the point \[(f\,,\,\infty\,,\,\infty\,,\,\infty)\in
\Rat_2^{\rm fm}\]. In fact, it follows immediately
from Lemma 3.1 that the moduli space \[\Rat_2
^{\rm fm}\] can be identified with the hypersurface consisting of all
points \[(\mu_1\,,\,\mu_2\,,\,\mu_3)\in\C^3\] satisfying the polynomial
equation
$$  \mu_1\mu_2\mu_3-(\mu_1+\mu_2+\mu_3)+2~=~0\;. $$
It is easy to check that this surface has a singular point
of type \[(2,2,2)\] at the point \[\mu_1=\mu_2=\mu_3=1\], corresponding to
a map with a triple fixed point. (Compare 6.2.)\QED

To summarize, we have a commutative diagram
$$	\matrix{\Rat_2^{\rm tm} & \buildrel 6\over\longrightarrow &
		 \Rat_2^{\rm cm}\cr
	{\scriptstyle 2}\downarrow~~ & & {\scriptstyle 2}\downarrow~~\cr
	\Rat_2^{\rm fm} & \buildrel 6\over\longrightarrow & \Rat_2\cr} $$
of holomorphic maps, where the two vertical maps are unramified 2-sheeted
coverings. The horizontal maps are 6-sheeted coverings,
ramified along the double fixed point
locus \[\prod(\mu_i-1)=0\], and with a more complicated ramification
along the triple fixed point orbit \[\mu_1=\mu_2=\mu_3=1\]. Similarly there
is a commutative diagram
$$      \matrix{{\cal M}_2^{\rm tm} & \buildrel 6\over\longrightarrow &
                 {\cal M}_2^{\rm cm}\cr
        {\scriptstyle 2}\downarrow~~ & & {\scriptstyle 2}\downarrow~~\cr
        {\cal M}_2^{\rm fm} & \buildrel 6\over\longrightarrow & {\cal M}_2\cr} $$
of holomorphic maps, where now all four maps have ramification points.
However, the projection \[{\cal M}_2^{\rm tm}\to{\cal M}_2^{\rm fm}\]
ramifies only at the singular point
\[\mu_1 =\mu_2=\mu_3=1\], so that \[{\cal M}_2^{\rm tm}\] can be considered
as a desingularization of \[{\cal M}_2^{\rm fm}\]. Similarly,
the projection\break \[{\cal M}_2^{\rm tm}\to
{\cal M}_2^{\rm cm}~~\] has an isolated ramification point at the unique
singular point\break
\[\mu_1\!=\!\mu_2\!=\!\mu_3=-2\],~ so that it can be considered at least
locally as a desingularization.\smallskip

{\bf Remark 6.7: Topology.} I know almost nothing about the homology or homotopy of the
spaces \[\Rat_2^{\rm tm}\] and \[\Rat_2^{\rm fm}\] and
\[{\cal M}_2^{\rm tm}\] and \[{\cal M}_2^{\rm fm}\], or about the fiber
bundle\break \[\Rat_1\subarr \Rat_2^{\rm tm}\darrow {\cal M}_2^{\rm tm}\].
Any information would be appreciated.
\smallskip

{\bf Remark 6.8: One marked fixed point.} Sometimes it is convenient to
consider maps with just one distinguished fixed point. (Compare
(23) in Appendix C.) The multiplier \[\mu\]
at this fixed point is then an invariant, and the product \[\tau\] of the
multipliers at the other two fixed points is also an invariant. The pair
\[(\mu\,\,,\tau)\in\C^2\] determines the map up to conjugacy. For we can
solve for \[\sigma_1=\sigma_3+2=\mu\tau+2\], and the sum of the multipliers at
the other two fixed points is equal to \[\sigma_1-\mu=\mu\tau+2-\mu\].
We can easily solve for \[\sigma_2=\mu(\mu\tau+2-\mu)+\tau\]. Thus \[\tau\]
is an affine parameter along the line \[\Per_1(\mu)\subset\M\].
\eject

{\bf Remark 6.9. Marked cubic polynomials.}
There is a completely analogous concept of markings for cubic
polynomial maps. (Compare 3.2.) A cubic polynomial
map is uniquely determined by its
fixed points \[z_i\] and its critical points \[\omega_j\], which are subject
only to the equality
\[~	(z_1+z_2+z_3)/3\,=\,(\omega_1+\omega_2)/2	~\]
of barycenters, and the inequality
$$	z_1\,z_2+z_1\,z_3+z_2\,z_3~\ne~3\omega_1\,\omega_2 ~.$$
(The equality \[z_1\,z_2+z_1\,z_3+z_2\,z_3\;=\;3\,\omega_1\,\omega_2\],
together with the equality of barycenters, would characterize triples
\[\{z_1\,,\,z_2\,,\,z_3\}\] which have a common image under
every cubic map with critical points \[\{\omega_1\,,\,\omega_2\}\].) From
this description, it is not difficult to check that the space of all
totally marked cubic polynomial maps is a manifold having the homotopy
type of the 
non-orientable 2-sphere bundle over a circle. The corresponding
moduli space, consisting of conjugacy classes of
totally marked cubic polynomial maps, is a complex manifold which can be
identified with the complement of a quadratic curve in \[\CP^2\]. This
moduli space has the homotopy type of \[\RP^2\].\smallskip

These descriptions seem rather complicated. In fact, in the
polynomial case it is
usually much more convenient to work with {\bit monic centered polynomials\/},
or sometimes with {\bit critically marked monic centered polynomials\/},
rather than getting into the complications of a fixed point marking.
However, for quadratic
rational maps there does not seem to be any correspondingly
convenient normal form.

\end

\def\Rat{{\rm Rat}}
\def\Map{{\rm Map}}
\def\Per{{\rm Per}}
\def\result{{\rm result}}
\def\CP{{\bf CP}}
\def\RP{{\bf RP}}
\def\sl{{\cal S}}
\def\QP{\medskip\leftskip=.4in\rightskip=.4in\noindent}
\def\subarr{\hookrightarrow}
\def\darrow{\rightarrow\!\!\!\!\!\rightarrow}
\def\M{{\cal M}_2}
\def\fm{^{\rm fm}}
\def\cm{^{\rm cm}}
\def\tm{^{\rm tm}} 
\ifnum\pageno<2 \pageno=20 \fi

\centerline{\bf \S7. Hyperbolic Julia Sets and Hyperbolic Components
in Moduli Space.}\medskip

This will be a brief description of results due to Rees and Tan Lei.
Recall that a rational map is {\bit hyperbolic\/} (ie., hyperbolic on its
Julia set) if and only
if the orbit of every critical point converges to some attracting periodic
orbit. The hyperbolic maps form an open subset of moduli space, and the
connected components of this open set are called {\bit hyperbolic
components\/}. Rees works with the critically marked moduli space \[{\cal M}_2
^{\rm cm}\]. (Compare \S6.) However we can work equally well with the unmarked
moduli space \[\M\] of \S3, or with the moduli space \[\M\fm\] or \[\M\tm\]
of \S6. She shows that the
hyperbolic components can be divided
into four classes, as follows. (The names are my own.)

{\bf Type B: Bitransitive.} Each of the two critical points belongs to the
immediate basin\footnote{$^1$}{The
{\bit immediate basin\/} of an attracting periodic point \[z_0=f^{\circ n}
(z_0)\] is the component of \[z_0\] in the open set consisting of all
points whose orbit under \[f^{\circ n}\] converges to \[z_0\], or
equivalently the component of \[z_0\] in the Fatou set \[\hat\C\ssm J(f)\].}
of some attracting periodic point, where these two periodic points
are distinct but belong to the same orbit. Evidently the period
must be two or more.\smallskip

{\bf Type C: Capture\footnote{$^2$}{\rm I am told that the word ``capture''
is used in a quite different sense, as a construction for passing from
quadratic polynomial maps to quadratic rational maps, in the thesis
of B. Wittner [W]. Unfortunately, I have not been able to get a copy of
this important work.}
.} Only one critical point belongs to the immediate
basin of a periodic point, but the orbit of the other critical point
eventually falls into this immediate basin.
Again the period must be two or more. (Compare 8.2 below.)\smallskip

{\bf Type D: Disjoint attractors.} The two critical points belong to the
attracting basins for two disjoint attracting periodic orbits. One particularly
interesting and important class of examples are those obtained by {\bit
mating~} two quadratic polynomial maps. (Compare the discussion below.)
\smallskip

{\bf Type E: Escape.} Both critical orbits converge to the same attracting
fixed point. It follows that the Julia set is totally disconnected. To fix
our ideas, we will take this fixed point to be the
point at infinity, and say that both critical orbits ``escape'' to infinity.
There is just one such hyperbolic component. It will be discussed in \S8.
\smallskip

For an analogous discussion of polynomial
maps, see [M3], [M4].\smallskip

\pageinsert
\insertRaster 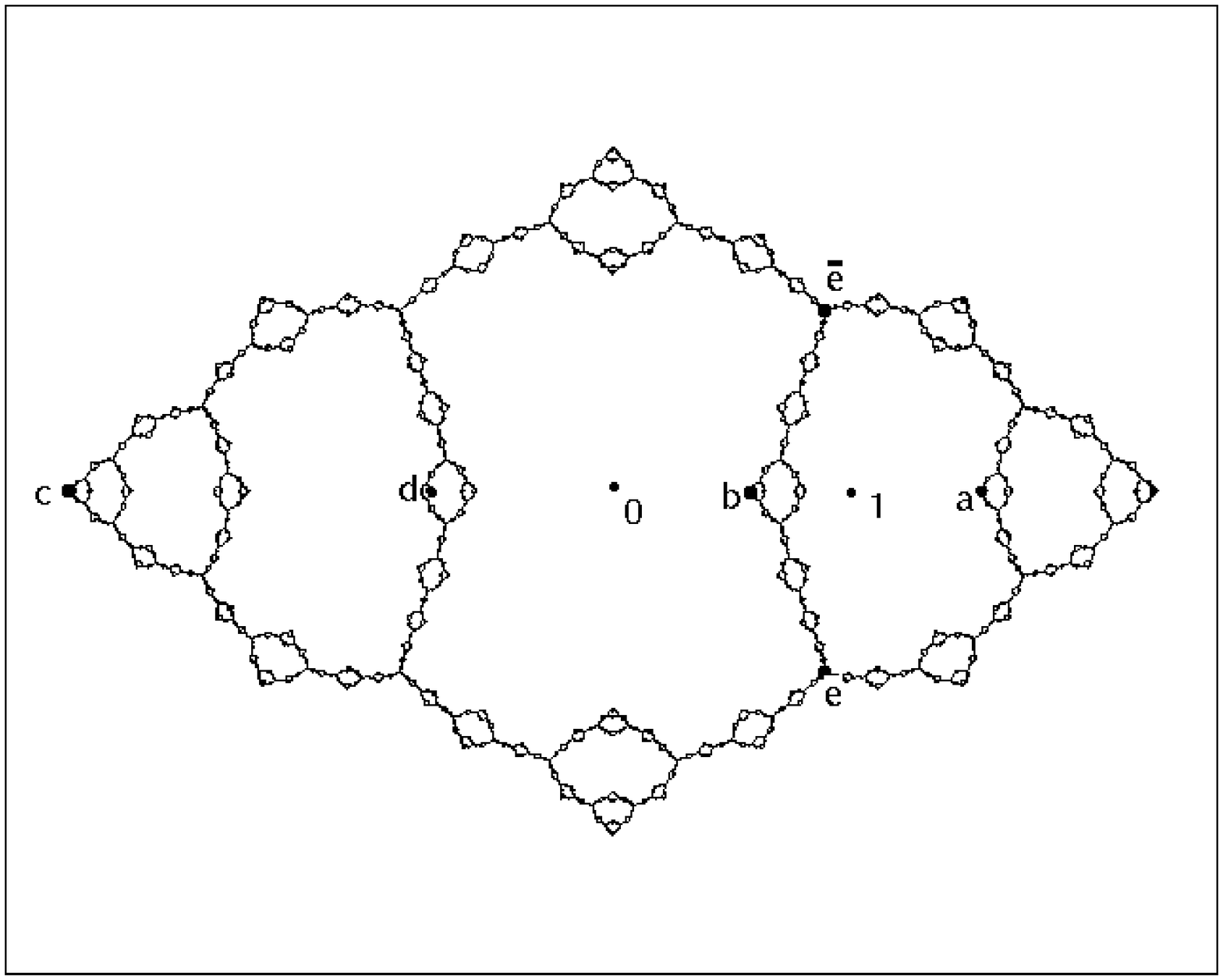 pixels 640 by 512 scaled 450
\smallskip
{\QP\bit Figure 2. Type B: Julia set for \[z\mapsto 1-1/z^2\], with both
critical points in the period three orbit \[0\mapsto\infty\mapsto 1\mapsto 0\].
The fixed points
are \[d\,,\,e\,,\,\bar e\], and the other period three orbit is \[a\mapsto b
\mapsto c\mapsto a\].\medskip}

\insertRaster 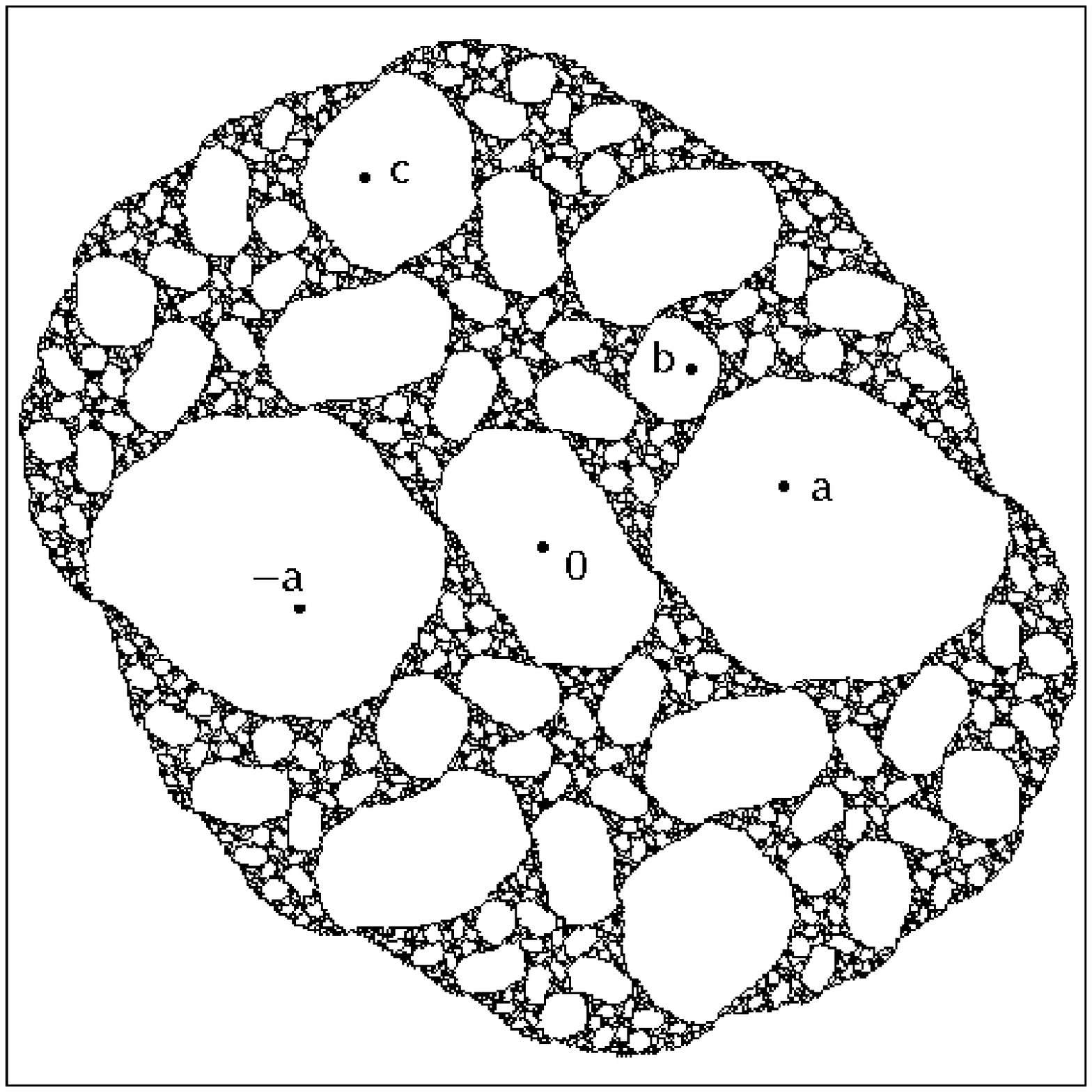 pixels 480 by 480 scaled 450
\smallskip
{\QP\bit Figure 3. Type C: Julia set for the map \[z\mapsto a+1/(z^2-a^2)\]
with\break \[a\approx.7467+.2286\,i\]. Here the periodic critical orbit \[~~\infty
\leftrightarrow a~~\] captures the\break orbit \[0\mapsto b\mapsto c
\mapsto -a\mapsto\infty\] of the other critical point.\par}
\endinsert

There are infinitely many hyperbolic components of each of the types B, C
and D, and they share many similarities. For every hyperbolic map \[f\]
of type B, C or D, the Julia set \[J(f)\] is connected and locally connected.
(Compare 7.1 and 8.2 below.) It follows that each component
of its complement \[\hat\C\ssm J(f)\] is conformally isomorphic to the unit disk.
Following McMullen,
Rees shows that each hyperbolic component of type B, C or D contains
precisely one map \[f_0\], called its ``center point'', which
is {\bit post-critically
finite\/}. (That is, the orbit of each critical point under \[f_0\]
is periodic or eventually periodic.) Furthermore, with one trivial exception,
she shows that each hyperbolic component of type B, C or D is a topological
4-cell. The unique exception is the hyperbolic component \[H_*\cm\subset
\M\cm\] of type B which
is centered at the unique singular point \[\langle z\mapsto 1/z^2\rangle\]
of the space \[\M\cm\].
Even in this exceptional case, we can get rid of the singularity
and obtain a topological cell simply by working in one of the other versions
of moduli space. In fact the corresponding hyperbolic component \[H_*\subset\M\]
is a topological cell, and its 6-sheeted ramified covering \[H_*\fm\subset\M\fm\]
is also a topological cell. The corresponding set in \[\M\tm\] is a disjoint
union of two copies of \[H_*\fm\].

\midinsert\vskip -.2in
\insertRaster 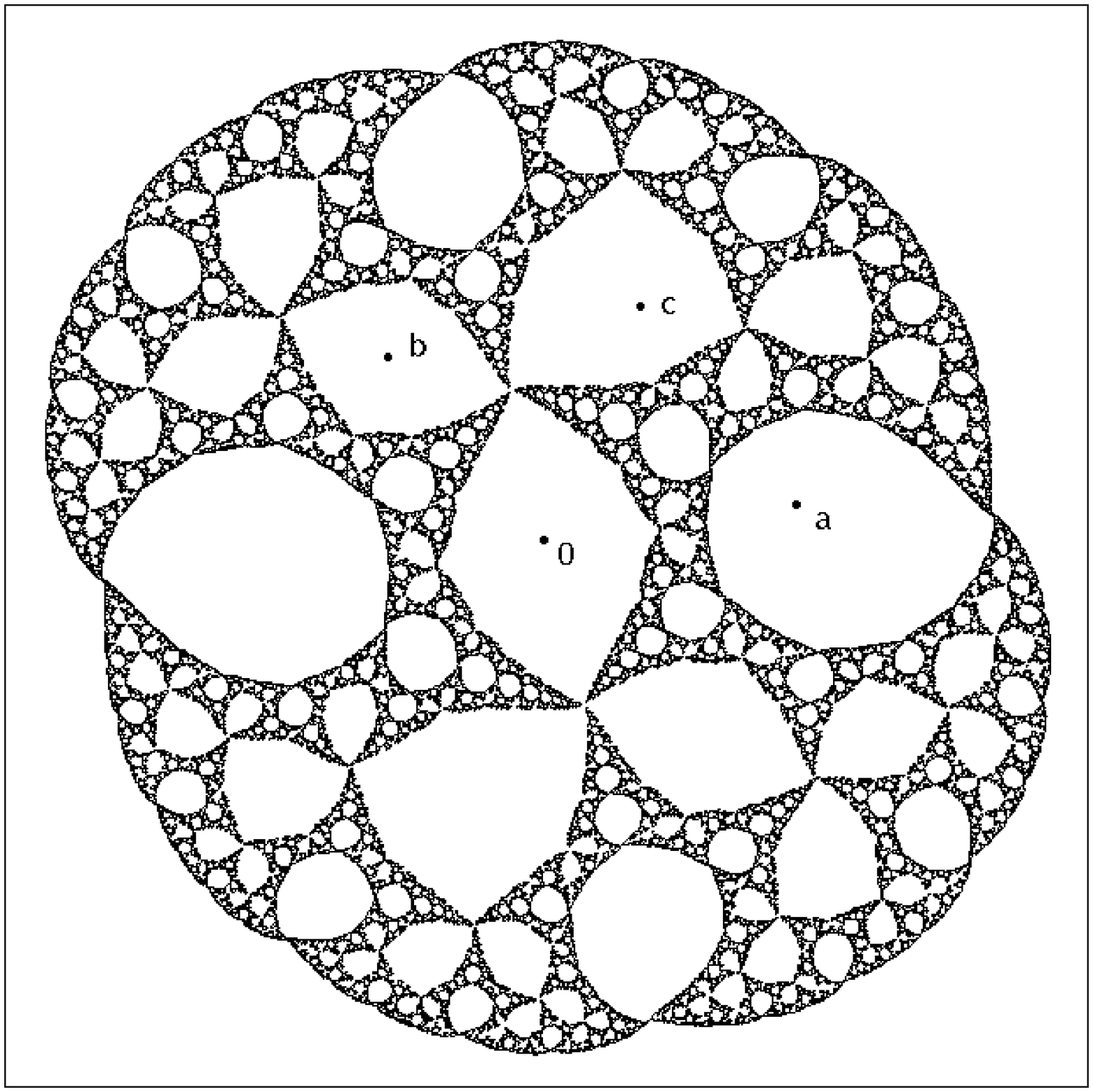 pixels 640 by 640 scaled 450
\vskip -.1in
{\QP\bit Figure 4. Type D: Julia set for the map \[z\mapsto a+1/(z^2-a^2)\]
with\break \[a\approx .8200+.1446\,i\]. (See [Bie].)
Here there are two disjoint periodic critical\break
orbits \[\infty\leftrightarrow
a\] and \[0\mapsto b\mapsto c\mapsto 0\]. This example is constructed by
mating.\par} 
\vskip -.1in
\endinsert

Some representative examples of Julia sets for post-critically finite
hyperbolic maps of types B, C, D are shown in Figures 2, 3, 4. (It is probably
impossible to distinguish these three types just by looking at the Julia set.)
In each case, the critical points have been placed at zero and infinity,
so that the symmetry of the Julia set is evident.

{\bf Mating.} Suppose that we are given two quadratic maps,
say \[f_a(z)=z^2+a\] and \[f_b(z)=z^2+b\], both with connected Julia set.
It is often possible to paste together the filled Julia sets for
\[f_a\] and \[f_b\] so as to obtain a Riemann sphere which decomposes into
two halves, one half with the dynamics of \[f_a\] and the other
with the dynamics of
\[f_b\]. (Compare [Ta], [Sh2], [W].)
{\it According to Tan Lei, in the post-critically finite case, this
is possible if and only if \[f_a\] and \[f_b\] do not belong to complex
conjugate limbs of the Mandelbrot set.\/} Evidently,
her construction yields many
examples of quadratic rational maps with two distinct superattractive cycles.
However, not every component of type D can be obtained in this way. Wittner
has described a real quadratic map with attracting cycles of period 3 and 4
which cannot be obtained by mating. (Compare Appendix F.)
\eject

{\bf Quadratic Mating Conjecture:} It is conjectured
that the mating \[f_a\amalg f_b\] can be well defined as an element
of \[\M\cm\] which depends continuously on \[a\]
and \[b\], whenever \[f_a\] and \[f_b\] do not belong to complex conjugate
limbs of the Mandelbrot set \[M\]. (For this purpose, it is convenient to
say that a map in  \[M\] belongs
to the \[t$-{\bit limb\/} if it either has a fixed point of
multiplier \[e^{2\pi it}\], or is separated from the central region
of \[M\] by the map which has
such a fixed point, where \[t\] can be {\it any\/}
angle in \[\R/\Z\].)

One way of trying to construct such a mating can be described as follows.
Let \[U(f_a\,,\,\epsilon)\] be the neighborhood of the filled Julia set
\[K(f_a)\] consisting of all points for which the
Green's function (canonical potential function)
\[G_a\] of \[K(f_a)\] takes values \[G_a(z)<\epsilon\].
Let \[M=M(f_a\,,\,f_b\,,\,\epsilon)\] be the compact Riemann surface obtained
from the disjoint union \[U(f_a\,,\,\epsilon)\sqcup U(f_b\,,\,\epsilon)\]
by identifying the open subset \[U(f_a\,,\,\epsilon)\ssm K(f_a)\] with
\[U(f_b\,,\,\epsilon)\ssm K(f_b)\] under the correspondence
\[(G,t)\leftrightarrow (\epsilon-G,-t)\], where \[G\] is the Green's function
and \[t\] is the external angle. By the Uniformization Theorem, \[M\] is
conformally isomorphic to the Riemann sphere \[\hat\C\]. More explicitly,
there is a
unique conformal isomorphism which takes the critical points of \[f_a\] and
\[f_b\] to \[+1\] and \[-1\] respectively, and takes the point with coordinates
\[G=\epsilon/2\,,~t=0\] to the point at infinity in
\[\hat\C\]. Now \[f_a\] and \[f_b\] fit together to yield a holomorphic
map from \[M(f_a\,,\,f_b\,,\,\epsilon)\] to \[M(f_a\,,\,f_b\,,\,2\epsilon)\].
Identifying each of these manifolds with \[\hat\C\]. This yields a holomorphic
map from \[\hat\C\] to itself which has critical points at \[\pm 1\] and a
fixed point at infinity. Thus it is a quadratic rational map of the form
\[\;	\amalg(f_a\,,\,f_b\,,\,\epsilon)\,:\,w~\mapsto~
 (w+w^{-1}+c)/\lambda~~, \]
where \[\lambda\] is the multiplier at infinity. (Compare (23) in Appendix C.)
If the coefficients \[c\]
and \[\lambda\] tend to well defined finite limits as \[\epsilon\to 0\], 
then the limiting rational map may be called the {\bit mating\/} \[~
f_a\amalg f_b\]. (For other, more standard,
forms of the definition, see [Sh2], [Ta], [Bie].)

{\bf Common Fatou Boundary Points.} In the examples of Julia sets which are
illustrated above, there are many pairs of Fatou components which have a
common boundary point, or (in Figure 2) even a Cantor set of common boundary
points. However, this need not be the case. Appendix F, written in
collaboration with Tan Lei, describes an example of a hyperbolic map
for which no two Fatou components have a common boundary point. The
corresponding Julia set is a ``Sierpinski carpet''.

{\bf A Compactness Question.} How can one decide whether
some given hyperbolic component has compact closure within the moduli
space \[\M\cong\C^2\],
or whether it is unbounded? Any hyperbolic component with an attracting
fixed point is certainly unbounded,
and Figures 7, 9, 10, 16, 17 show many other examples of unbounded components.
On the other hand, consider a component obtained by mating
two hyperbolic components of period \[\ge 2\] of the Mandelbrot set.
If the Quadratic Mating Conjecture is true, then it is easy to show that
such a component has compact closure, canonically homeomorphic to
\[\bar D\times\bar D\] (at least in the critically marked case).
McMullen [Mc] has conjectured that whenever a degree \[d\] Julia set \[J(f)\]
is a Sierpinski carpet
(Appendix F) the corresponding hyperbolic component in \[{\cal M}_d\] has
compact closure.\smallskip

To conclude this section, we prove the following.

{\QP{\bf Lemma 7.1.} \it If the Julia set of a hyperbolic map is connected,
then it is locally connected.\smallskip}
\eject

{\bf Proof outline} (suggested by Douady). It suffices to consider the
post-critically finite case. For according to [Mc] every hyperbolic map
with connected Julia set can be
deformed through hyperbolic maps to a
post-critically finite one, and according to [MSS] or [Ly] the topology of
the Julia set, and the isotopy class of its embedding into \[\hat\C\]
does not change under such a deformation. In the case of a
periodic Fatou component \[U\], every internal ray
from the super-attracting point
lands at a well defined boundary point, which depends continuously on the
initial angle. (Compare the argument in
[DH2, p. 24].) Since the continuous image of a locally
connected space is locally connected, it follows that the boundary of \[U\]
is not only connected but also locally connected. Since every Fatou component
is eventually periodic ([Su]), it follows that: {\it
the boundary of every Fatou component is connected and locally connected.\/}
Let \[V\subset\hat\C\] be the complement of the post-critical set for \[f\].
We will use the Poincar\'e metric on \[V\], and
on its universal covering manifold \[\widetilde V\]. Note that \[f^{-1}\]
lifts to a well defined contracting map \[\widetilde{f^{-1}}\] on
\[\widetilde V\]. Each Fatou component \[U\] which contains no post-critical
point lifts homeomorphically to a
subset \[U'\subset\widetilde V\], and \[\widetilde{f^{-1}}\] maps \[U'\]
to a set of strictly smaller diameter in the Poincar\'e metric. In fact
a straightforward compactness argument shows that the diameter shrinks by
a factor bounded away from \[1\]. {\it
For any \[\epsilon>0\], it follows that there are only finitely
many Fatou components of diameter \[>\epsilon\].\/} From this,
it follows easily that \[J\] is locally connected. In fact, if two points of
\[J\] are close to each other, then we can find a connected subset of \[J\]
of small diameter containing both points as follows. Draw a straight line
segment between them, and replace its intersection with each Fatou component
\[U\] with the entire boundary of \[U\] if \[U\] is small, or with a suitable
small connected subset of \[\partial U\] otherwise.\QED

\bigskip

\centerline{\bf \S8. The ``Escape Locus'': Totally Disconnected Julia Sets.}\smallskip
First a lemma above rational maps of arbitrary degree, to be
proved in Appendix E.

{\QP{\bf Lemma 8.1.} \it If all of the critical values of a rational map
are contained in a single component of the Fatou set \[\hat\C\ssm J\], then
the Julia set \[J\] is totally disconnected. 
Every orbit
in the Fatou set converges either to an attracting fixed point, or to a
parabolic fixed point of multiplicity precisely equal to two.\medskip}

Here the multiplicity \[m\ge 1\] of a (finite) fixed point \[f(z_0)=z_0\]
is defined to be the degree of the first non-zero term in the Taylor expansion
about \[z_0\],
$$	f(z) -z~=~\alpha(z-z_0)^m+{\rm (higher~terms)~,\qquad with}\quad
	\alpha\ne 0~. $$
In the quadratic
case, there is an explicit dichotomy, and an effective criterion.
See Yin [Y] and Makienko [Mak], as well as [Pr1], [R3].

{\QP{\bf Lemma 8.2} \it
The Julia set \[J\] of a quadratic rational map is
either connected, or totally disconnected and
homeomorphic, with its dynamics,
to the one-sided shift on two symbols. It is totally disconnected
if and only if either:

\ref
~~$(a)$ both critical orbits converge to a common attracting fixed point, or

\ref
~~$\,(b)$ both critical orbits converge to a common fixed point of multiplicity
two but neither critical orbit actually lands on this point.
\medskip}

The proof will be given at the end of this section.
\medskip

{\bf 8.3. Quadratic Examples.} For the map \[\;f(z)=
z+2+z^{-1}\,,\;\] one critical orbit \[1\mapsto 4\mapsto 6.25\mapsto
\cdots\]
converges to the parabolic fixed point at infinity, while the other critical
orbit \[-1\mapsto 0\mapsto\infty\] actually lands at this fixed point.
Thus the Julia set is connected. In
fact \[J(f)\] is the interval \[[-\infty,0]\]. For \[f(z)=
\pm(z+z^{-1})\], both critical orbits
\[~	\pm 1\,\mapsto\,\pm 2\,\mapsto\,\pm 2.5\,\mapsto\,\cdots ~ \]
converge to the parabolic point at infinity. Since the multiplicity
of the fixed point at infinity is either three (if the multiplier is \[+1\])
or one\break (if the multiplier is \[-1\]), the Julia set is again
connected. It equals the imaginary
axis plus the point at infinity. On the other hand, for \[z\mapsto z+c+
z^{-1}\] with \[c>2\], the Julia set is totally disconnected.\smallskip

{\bf 8.4. Cubic Examples.} The situation for maps of higher degree is more
complicated.  The rational map
\[z\mapsto 2+2z^3/(27(2-z))\] has connected Julia set, although the critical
points \[0,\,0,\,3,\,\infty\] are all contained in an orbit \[3\mapsto
0\mapsto 2\mapsto\infty\], which lands at a superattracting fixed point.
The polynomial map \[f(z)=z^3-12z+12\] has a Cantor
set, contained in the interval \[[-4\,,\,3]\], as Julia set, although
the critical point \[+2\] belongs to this Julia set. The
Julia set for the map \[z\mapsto {5\over 36}\,(z-3)^2(z+4)\]
is disconnected, but contains the connected interval \[[0\,,\,5]\]. For
a thorough analysis of such polynomial examples, see [BH].
\smallskip

{\bf The escape locus.}
There is just one hyperbolic component of type E in the moduli space \[\M\].
We will call it the {\bit hyperbolic
escape locus\/}, and denote it by \[E\subset\M\]. (Perhaps a better term
would be ``shift locus'' or ``totally disconnected locus'' ?)
If we work
with critically marked conjugacy classes, there is a corresponding
unique component\break
\[E\cm\subset\M\cm\]. This escape component is very different
from the other hyperbolic components.
For a map \[f\] of type E, the Julia set \[J(f)\] is a Cantor
set, and its complement is a connected open set of infinite connectivity.
Such a map can never be post-critically finite, so this hyperbolic
component \[E\] does not have any preferred center point.

Here is a well known collection of examples. Consider the one-parameter
family of quadratic polynomials \[~z~\mapsto ~z^2+c~,\]
which can be identified with the one-dimensional
slice \[\sigma_3=0\] through the moduli space \[\M\].
The intersection of this family with the escape locus \[E\]
is precisely the complement of the {\bit Mandelbrot set\/.}
This complement is conformally isomorphic to \[\C\ssm\bar D\],
with free cyclic fundamental group.
The corresponding
escape locus for polynomials of higher degree has been studied by Blanchard,
Devaney and Keen, who show that it has a very rich fundamental group.

Using the critically marked moduli space,
Rees shows that the escape locus \[E^{\rm cm}\] has a rather
complicated topology. In particular, her description implies
that \[E^{\rm cm}\] has the homotopy type of a Klein bottle, and hence
has non-abelian fundamental group. However, the unmarked escape locus
\[E\subset\M\] has a much simpler description.

{\QP{\bf Lemma 8.5.} \it The hyperbolic escape locus \[E\subset\M\] is
homeomorphic to the product \[D\times (\C\ssm\bar D)\],
where \[D\] is the open unit disk.\smallskip}

More explicitly, if \[\mu=\mu(f)\in D\] is the multiplier
at the unique attracting fixed point of \[f\], we will show that the
correspondence \[\langle f\rangle\mapsto \mu(f)\] yields
a fibration of \[E\] over \[D\], with fiber \[E_\mu\]
homeomorphic to \[\C\ssm\bar D\] (or equivalently to \[D\ssm\{0\}\]).
As in 6.8, it will be convenient to
work with the coordinates \[(\mu,\tau)\], where \[\mu\] is the multiplier
at the preferred (attracting) fixed point and \[\tau\] is the product of the
multipliers at the other two (repelling) fixed points. Thus each coordinate
line \[\mu={\rm constant}\in D\] in the \[(\mu,\tau)$-plane corresponds to
a straight line \[\Per_1(\mu)\subset\M\],
consisting of all conjugacy classes of maps which have
an attracting fixed point of multiplier \[\mu\]. (Compare \S3.
These various straight lines \[\Per_1(\mu)\]
intersect each other within \[\M\], but
not within the escape locus, since  for \[\langle
f\rangle\in E\] the map \[f\] has only one attracting fixed point.)

The proof will be based on Goldberg and Keen [GK].
For each fixed\[\mu\in D\], Goldberg
and Keen show that the complement \[M_\mu\] of the escape locus in
\[\Per_1(\mu)\cong\C\] is canonically homeomorphic to the Mandelbrot set
\[M=M_0\]. (This statement is quite
likely true even in the limiting case \[\mu=1\]; compare Figure 5.)
They prove this as an application of the theory of polynomial-like mappings
(Douady and Hubbard [DH3]).
For \[\langle f\rangle\in M_\mu\] the map \[f\] restricted to
a suitable neighborhood of its Julia set
is polynomial-like of degree two with connected Julia set, and hence is
hybrid-equivalent to a unique conjugacy class in \[M_0\]. Furthermore,
by [DH3, \S7] the resulting map from \[\bigcup_\mu M_\mu\subset\M\] to
\[M_0\] is continuous.

{\bf Caution:} This statement is formulated rather differently in [GK],
since the authors work with marked critical points
and hence study a 2-fold
branched covering space \[\Per_1(\mu)\cm\], which
contains two disjoint copies of the Mandelbrot set for \[\mu\ne 0\].
Note also that some notations (such as \[\Rat_2\]) which are used in [GK]
are incompatible with the notations used here.

It follows from the Riemann Mapping Theorem that the complement
$$E_\mu~=~\Per_1(\mu)\cap E~\cong~\C\ssm M_\mu $$
is conformally isomorphic to the region \[\C\ssm\bar D\].
In fact we can choose a canonical conformal isomorphism
$$	\psi_\mu:\C\ssm\bar D~\to~\C\ssm M_\mu~, $$
normalized so that the multiplier at infinity, \[~\lambda(\mu)\,=\,
1/\lim_{z\to\infty}\psi'_\mu(z)~\]
is real and positive. We must show that \[\psi_\mu\] depends
continuously on \[\mu\], using the topology of locally uniform convergence,
so that the correspondence
$$(\mu\,,\,z)~\mapsto ~(\mu\,,\,\psi_\mu(z))~\in~\mu\times  E_\mu $$
yields the required homeomorphism between \[D\times(\C\ssm\bar D)\]
and the escape locus \[E\].

Note that this multiplier \[\lambda(\mu)\] provides an invariant measure
of the ``size'' of the open set \[E_\mu=\psi_\mu(\C\ssm\bar D)\]. More
generally, if \[U_1\subset U_2\] are simply-connected neighborhoods of
infinity in \[\hat\C\] and if \[\lambda_1\] and \[\lambda_2\] are the
multipliers at infinity of the corresponding Riemann maps
\[\C\ssm\bar D\to U_i\], then it follows
easily from the Schwarz Lemma that \[\lambda_1\le\lambda_2\], with equality
if and only if \[U_1=U_2\].

If this correspondence \[\mu\mapsto\psi_\mu\]
were not continuous, then we could choose a sequence
of points \[\mu_i\in D\] converging to a limit \[\hat\mu\in D\] so that the
corresponding sequence of functions \[\psi_{\mu_i}\] on \[\C\ssm\bar D\]
did not converge locally uniformly to \[\psi_{\hat\mu}\]. However, since the
\[\psi_{\mu_i}\] belong to a normal family, we may assume after
passing to a subsequence that these functions do converge
locally uniformly to some univalent limit \[\hat\psi\ne\psi_{\hat\mu}\].
It is not hard to check that the image \[\hat\psi(\C\ssm\bar D)\] must be
a proper subset of the open set \[\psi_{\hat\mu}(\C\ssm\bar D)\]. Therefore,
as noted above, the multiplier \[\hat\lambda\] of \[\hat\psi\] at
infinity must be strictly less than the multiplier \[\lambda(\hat\mu)\] of
\[\psi_{\hat\mu}\], say \[\hat\lambda<\lambda(\hat\mu)/(1+\epsilon)\]
with \[\epsilon>0\]. Now the circle of radius \[1+\epsilon\] in \[\C\ssm\bar
D\] corresponds under \[\psi_{\hat\mu}\] to a loop which encloses a small
neighborhood of the compact set \[M_{\hat\mu}\]. For \[\mu_i\] sufficiently
close to \[\hat\mu\] the compact set \[M_{\mu_i}\] must be contained in this
neighborhood, and it follows easily that
$$	\lambda(\mu_i)~\ge~\lambda(\hat\mu)/(1+\epsilon)~. $$
Passing to the locally uniform limit as \[i\to\infty\], since
\[\lambda(\mu_i)\to\hat\lambda\], we obtain a contradiction.\QED

{\bf Remark 8.6: The escape locus in \[\widehat\M\].} If we work with
the compactified
moduli space \[\widehat\M\] of \S4, then the appropriate ``hyperbolic escape
locus'' \[E^*\subset\widehat\M\] is a topological 4-cell with a preferred
center point, just like all other hyperbolic
components. In fact, let \[E^*\] consist of \[E\] together with the
2-cell consisting of all ideal points \[\{\mu\,,\,\mu^{-1}\,,\,\infty\}\]
for which \[|\mu|<1\]. Then \[E^*\] fibers over the open disk with an open
disk as fiber; and each fiber contains one and only one ideal point.
By definition, the ``center'' of this filled in component is
the ideal point \[\langle 0\,,\,\infty\,,\,\infty\rangle\]. This center point
can perhaps be identified with the improper map \[\langle z\mapsto z^2+\infty
\rangle\], or with the limit of a sequence of quadratic maps ``tending''
to a constant map.
\medskip

Next let us study the escape locus \[E\cm\subset\M\cm\], with marked
critical points. We will prove the following. (For a more precise description
of \[E\cm\], see [R3].)

{\QP{\bf Lemma 8.7.} The critically marked escape locus
\[E\cm\] contains a Klein bottle as retract. Hence the fundamental group
of \[E\cm\] is non-abelian.\medskip}

This space \[E\cm\] can be described as
a 2-fold branched covering of \[E\], branched along the symmetry locus
\[\sl\cap E\], which consists of all pairs \[(\mu\,,\,\tau)\] with \[\mu\in
 D\ssm\{ 0\}\] and \[\tau=({2-\mu\over\mu})^2\]. However, it seems somewhat
easier to take a different approach. We will first
work with fixed point markings, and describe the corresponding escape
locus \[E^{\rm fm}\subset{\cal M}_2^{\rm fm}\]. Recall that a point in
\[\M\fm\] is specified by a 3-tuple \[(\mu_1\,,\,\mu_2\,,\,\mu_3)\] satisfying
the cubic equation (1). For a map in the escape locus, one of the three
fixed points must be attracting and the other two must be repelling. Thus
\[E\fm\] actually splits into three distinct components, depending on which
of the three marked fixed points
is attracting. To fix our ideas, suppose that \[|\mu_1|\,,\,
|\mu_2|>1>|\mu_3|\]. We will take the pair \[(\mu_1\,,\,\mu_2)\in(\C\ssm\bar D)
\times(\C\ssm\bar D)\] as independent parameters, solving for \[\mu_3\] by
equation (2). Of course not every such pair determines a map in the escape
locus, or even a map with \[|\mu_3|<1\]. We first prove the following.

{\QP{\bf Lemma 8.8.} \it If \[\;|\mu_1|>6\;\] and \[\;|\mu_2|>6\;\],
then the associated map \[f\] belongs to the escape
locus.\smallskip}

This estimate is probably far from sharp. (The largest \[|\mu_1|\] and
\[|\mu_2|\] which I know, for \[f\] outside the escape locus, is given
by \[\mu_1=\mu_2=-3\], for \[f(z)=-(z+z^{-1})\].)\smallskip

{\bf Proof of 8.8.} We will use the normal form
$$	f(z)~=~{z\,(z+\mu_1)\over \mu_2\,z+1}~=~{\mu_1+z~~\over ~\mu_2+
	z^{-1}} $$
of equation (3). (Compare Appendix C.)
If \[|\mu_1|>6\] and \[|\mu_2|>6\], then a brief computation shows
that \[f\] maps the
annulus \[\;2/|\mu_2|\,<\,|z|\,<\,|\mu_1|/2\;\] into a compact subset of
itself. Furthermore,
the polynomial \[\mu_2\,z^2+2z+\mu_1\], whose
roots are the critical points of \[f\], is strictly non-zero outside of this
annulus, hence both critical points are contained in the annulus.
Using the Poincar\'e metric for this annulus,
we see that both critical orbits must converge
to a common attracting fixed point. The conclusion then follows by 8.2.\QED

{\bf Proof of 8.7.}
We can now easily construct a mapping from our preferred
component of \[E\fm\] to the torus by the correspondence
$$	(\mu_1\,,\,\mu_2)~\mapsto~(\arg(\mu_1)\,,\,\arg(\mu_2))~. $$
The subset \[|\mu_1|=|\mu_2|=7\] maps homeomorphically to the torus,
and hence is embedded in \[E\fm\] as a retract.\smallskip

Now let us also mark the critical points, and hence pass
to the 2-sheeted covering manifold \[E\tm\]. (The single
ramification point for the map \[\M\tm\to\M\fm\]
does not belong to the escape locus, and hence causes no difficulty.)
It is not difficult to show that a choice of critical point for \[f\]
is equivalent to
a choice of sign for \[\pm\sqrt{1-\mu_1\mu_2}\]. (Compare the discussion
following (22) in Appendix C.) Since the ratio
\[(1-\mu_1\mu_2)/(\mu_1\mu_2)\] always lies in the left half-plane when
\[|\mu_1\mu_2|>1\], this is equivalent to
making a choice of sign for the geometric mean \[\eta=\pm\sqrt{\mu_1\mu_2}\].
In other words, our component of \[E\tm\] can be identified with
an open subset of the manifold
consisting of all \[(\mu_1\,,\,\mu_2\,,\,\eta)\in(\C\ssm\bar D)^3\] for
which \[\eta^2=\mu_1\mu_2\]. In particular, we obtain an explicit retraction
onto the torus consisting of all
$$ (\mu_1\,,\,\mu_2\,,\,\eta)~\in~\C^3\quad{\rm for~
which}\quad|\mu_1|= |\mu_2|=|\eta|=7\quad {\rm and}\quad\eta^2=\mu_1\mu_2~.
\eqno (16) $$
In order to obtain the hyperbolic component \[E\cm\], as described by Rees,
we must pass to a quotient space by identifying under the involution of
\[E\tm\] which interchanges the role of the two repelling fixed points.
A brief computation shows that this corresponds to the fixed point free
involution
$$	(\mu_1\,,\,\mu_2\,,\,\eta)~\leftrightarrow
	~(\mu_2\,,\,\mu_1\,,\,-\eta)~. $$
The quotient space \[E\cm\] is still a smooth manifold. If we collapse
the torus (18) under this
involution, we obtain a Klein bottle. Thus this argument proves
that \[E\cm\] contains a Klein bottle as retract, and hence that \[\pi_1(E
\cm)\] retracts onto the fundamental group of a Klein bottle.\QED

(If we try to carry out the same argument for the unmarked escape locus
\[E\subset\M\], then we must simply work with the torus \[|\mu_1|=|\mu_2|=7\]
and the orientation reversing involution \[(\mu_1\,,\,\mu_2)\mapsto(\mu_2\,,\,
\mu_1)\]. In this case, the quotient space is a M\"obius band, which has the
homotopy type of a circle. Thus we could conclude only that \[E\] contains a
circle as retract.)
\eject

Recall that \[\M\cm\] has the homotopy type of a 2-sphere.
It is conjectured that the inclusion map \[E\cm\to\M\cm\] is homologically
non-trivial, or more explicitly that it induces an isomorphism of 2-dimensional
homology with mod 2 coefficients.
\bigskip

{\bf Proof of 8.2.}
To conclude this section, let us prove that every quadratic Julia set is either
connected or totally disconnected.
Note that \[J\] is connected if and only if every
component of the Fatou set \[\hat\C\ssm J\] is simply-connected. According to
Sullivan, every component of \[\hat\C\ssm J\] is eventually periodic, and every
periodic component is either a Siegel disk, a Herman ring, or an immediate
basin for some attracting or parabolic point. According to Shishikura [Sh1],
there are no Herman rings in the quadratic case.
We will make frequent use of the fact the a ramified
covering of a simply-connected region which has only one ramification
point must again be simply-connected.

First suppose that every periodic Fatou component is simply-connected.
Any cycle of Fatou components must either contain a critical point (in the
attracting or parabolic cases), or have a critical orbit which is dense
in its boundary
(in the Siegel case). Thus there is at most one critical point left over.
It follows inductively that every Fatou component is simply-connected.

Now suppose that some periodic Fatou component is not simply-connected.
Evidently it must be an immediate basin for some attracting or
parabolic periodic point. Such an immediate basin can be reconstructed
from a simply-connected neighborhood of the fixed point in the attracting
case, or from a simply-connected petal in the parabolic case, by taking a
direct limit of successive ramified coverings, ramified only at the critical
values. The result will again be simply-connected, unless both critical
points belong to the same connected component. But in that case, it follows
from Lemma 8.1 that \[J(f)\] is totally disconnected. (Compare [R3] and
Appendix E.)

If both critical orbits converge to the same attracting or multiplicity two
fixed point, then at least one critical point must belong to the immediate
basin \[U\]. Hence the map \[f\] restricted to \[U\] is two-to-one. Since
\[f\] is quadratic, and since \[U\] is fixed under \[f\],
this implies that \[U\] is fully invariant: \[U
=f^{-1}(U)\]. Hence both critical
points belong to \[U\], and again we can apply 8.1.

Finally, we must prove that a totally disconnected \[J\] is homeomorphic
to a one-sided two-shift. In the hyperbolic case, this is proved in [GK].
The key step is the construction of an embedded disk \[\Delta\] which
contains both critical values and is forward invariant,\break
\[f(\Delta)\subset
\Delta\]. The two components of \[\hat\C\ssm f^{-1}(\Delta)\] then cover the
Julia set, and form the required Bernoulli partition. In the parabolic
case, we modify this argument by constructing a simply connected open
petal \[P\] which contains both critical values and satisfies \[f(P)\subset
P\]. Again the two components of \[\hat\C\ssm f^{-1}(P)\] cover the Julia
set and form the required Bernoulli partition. The proof that a point in the
Julia set is uniquely determined by its symbol sequence with respect to this
partition is now more delicate. However, since a completely analogous argument
is carried out in Appendix E, details will be left to the reader.\QED

\vfil\eject

\centerline{\bf\S9. Complex 1-Dimensional Slices.}\medskip

The moduli space \[\M\]  can be thought of as a kind of table of contents,
with each point \[\langle f\rangle\in\M\] corresponding to a different form
of dynamic
behavior. Since this space is a complex 2-manifold, it is difficult to
visualize all of it it directly. Hence it may be helpful to try to
describe what kinds of behavior occur for \[\langle f\rangle\] belonging to
some 1-dimensional slice through \[M\].\smallskip

As a first example of a dynamically
interesting 1-dimensional slice, suppose that we choose some fixed integer
\[n\ge 1\] and some fixed complex number \[\mu\] and
consider the curve \[\Per_n(\mu)\subset\M\] consisting
of all conjugacy classes of maps which possess a periodic
point of period \[n\] and multiplier \[\mu\]. (See 3.4, 4.2.)
The case \[\mu=0\] is of particular interest. (Compare [R4], [R5],
[M3], [M5], and see Figures 7, 9.) Evidently the center point of
every hyperbolic component must belong to at least one of these curves
\[\Per_n(0)\]. For any
\[\mu\] in the closed disk \[|\mu|\le 1\], the fact that \[\langle f\rangle\]
belongs to \[\Per_n(\mu)\] imposes very strong restrictions on the dynamics
of \[f\]. The cases where \[\mu\] is a root of unity are noteworthy, since
these curves contain faces where two or more hyperbolic components of \[\M\]
come together along a common boundary. (Figures 5, 6, 8.)\smallskip


We can also consider curves for which one critical orbit is eventually
periodic, say \[f^{\circ t}
(\omega_1)=f^{\circ t+n}(\omega_1)\]. See Figure 11 for the case \[t=2\,,~
n=1\]. These curves also contain
common boundary faces between two different hyperbolic components.
The symmetry locus of \S5 is another curve of interest. (Figure 12.)
\smallskip

As a final interesting family of curves, for each integer \[t\ge 1\] we
can consider all \[f\] such that the \[t$-th forward image of one critical
point is equal to the other critical point, \[f^{\circ t}(\omega_1)=\omega_2\].
See Figure 13 for the case \[t=1\]. Note that the center
point of every hyperbolic component
of type B or C must belong to at least one of these curves.\medskip

Here is a more detailed discussion of some sample curves.\smallskip

\midinsert
\vskip -.1in
\insertRaster 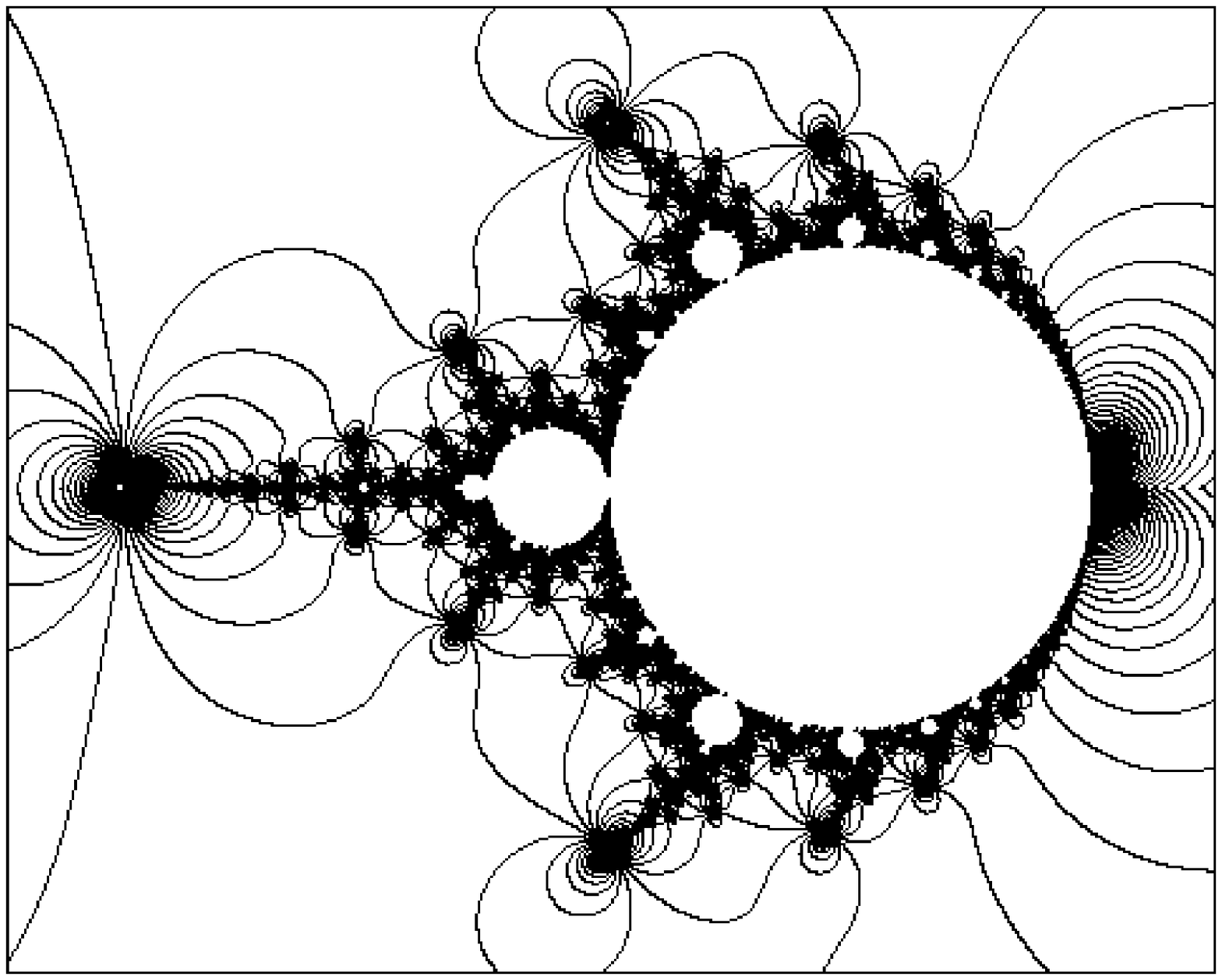 pixels 480 by 384 scaled 450
{\QP\bit Figure 5. The bifurcation locus
in \[\Per_1(1)\] seems to be a homeomorphic
copy of the Mandelbrot set boundary,
with the cusp point straightened out. Here \[\Per_1(1)\] is represented as
the \[-a^2$-plane for the family of maps \[f_a(z)=
z+a+1/z\], having a double fixed point of multiplier \[+1\] at infinity,
and having critical points at \[\pm 1\]. (Compare Appendix C.)
The
other fixed point has multiplier \[1-a^2\], and the period 2 orbit has
multiplier \[9-4a^2\]. Curves of
the form $${\rm Re}(f_a^{\circ n}(\pm 1)/a)- n={\rm constant}$$
for large \[n\] are also shown.\par}
\vskip -.1in
\endinsert

\pageinsert
\insertRaster 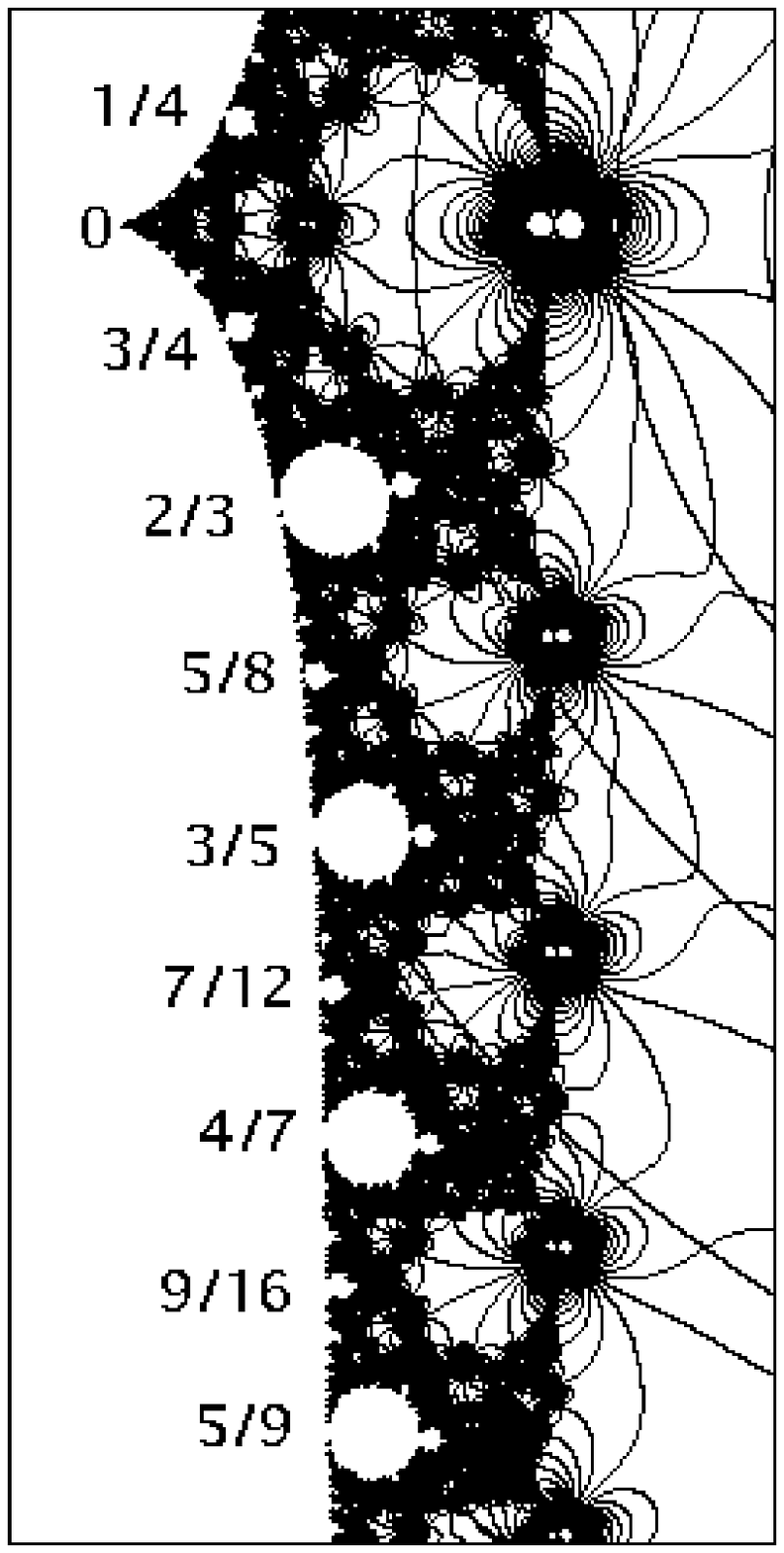 pixels 240 by 480 scaled 450
\smallskip
{\QP\bit Figure 6. The bifurcation locus in \[\Per_1(-1)=\Per_2(1)\]
is a distorted copy of the Mandelbrot set boundary with the \[1/2$-limb
missing (collapsed to the point at infinity). This picture shows
part of the \[-a^2$-plane for the family of maps \[z\mapsto
-(z+a+z^{-1})\], having a fixed point of multiplier \[-1\]. The notation
\[p/q\] indicates that the map also has a fixed point of multiplier \[e^{2\pi
 ip/q}\]. (The case \[p/q=1/2\] would correspond
to the point at infinity.)
\par}
\insertRaster 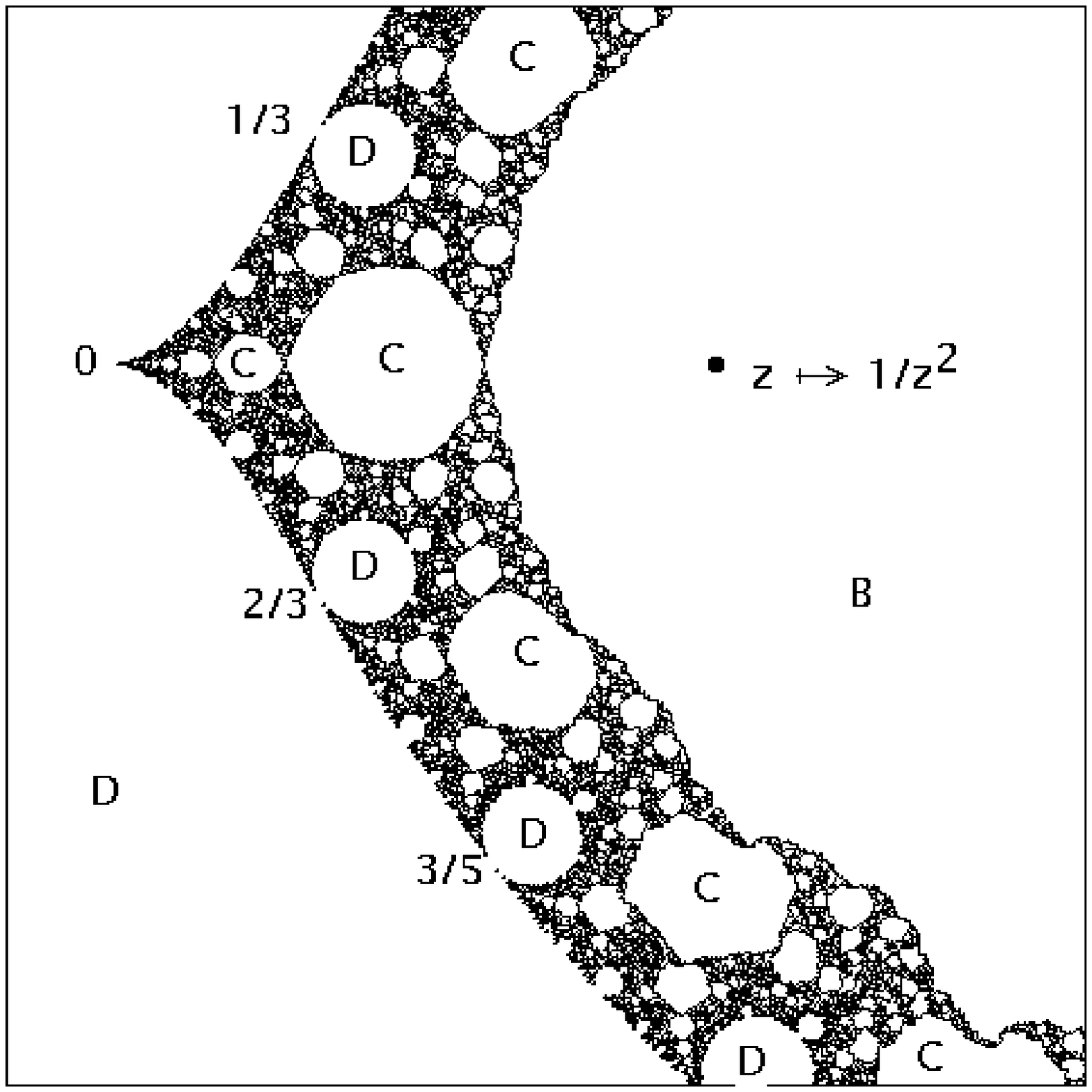 pixels 480 by 480 scaled 450
\smallskip
{\QP\bit Figure 7. Bifurcation locus in \[\Per_2(0)\]
(part of the \[-a^3$-plane for the family of maps \[z\mapsto
a/z+1/z^2\], with a critical point of period 2 at \[z=0\]).
This appears to be a homeomorphic copy of Figure 6. The types
of some of the more conspicuous hyperbolic components are indicated.
\par}
\endinsert

{\bf Period 1.} As noted in 3.4, the curves \[\Per_1(\mu)\]
are particularly easy to describe: Each locus
\[\Per_1(\mu)\] is a straight line of slope \[\mu+\mu^{-1}\]
in the \[(\sigma_1\,,\,
\sigma_2)$-coordinate plane, given for example by the equation
$$	\mu^3-\sigma_1\,\mu^2+\sigma_2\,\mu+2-\sigma_1~=~0~. \leqno
	\qquad \Per_1(\mu): $$
First consider the case \[\mu=0\]. Putting the fixed critical point
at infinity, we can use
the normal form \[z\mapsto z^2+c\]. (Here \[\sigma_1=2\] and
\[\sigma_2=4c\].) The corresponding bifurcation locus\footnote{$^1$}
{\rm By the {\bit bifurcation locus\/} for a parametrized
family of maps \[\{f_a\}\] we mean
the closed set consisting of parameter values \[a\]
for which the associated Julia
set \[J(f_a)\] does not vary continuously
under deformation of \[a\]. Compare [MSS], [Ly].}
in the \[c$-plane 
is the boundary of the familiar Mandelbrot set. The hyperbolic
components which intersect \[\Per_1(0)\] all have type D or E.
(Remark: If we lift to the
critically marked moduli space \[\M\cm\], then the locus \[\Per_1(0)\]
lifts to a union of two lines which intersect transversally. In fact,
using coordinates \[A,\,B,\,C\] as in 6.1,
the locus \[\Per_1\cm(\mu)\] splits as the union of a line
\[A=B=0\] for which the first critical point is fixed and a line \[A=C=0\]
for which the second critical point is fixed. These two lines intersect
transversally at the point \[A=B=C=0\], corresponding to the map \[z\mapsto
z^2\] with both critical points fixed.)

As noted in the proof of 8.5, the bifurcation locus
in the line \[\Per_1(\mu)\] varies by a continuous isotopy
as \[\mu\] varies within the open unit disk. As shown in Figure 5, it even
seems to retain the same topology as \[\mu\] tends non-tangentially to \[+1\].
On the other hand, as
\[\mu\to -1\] non-tangentially,
the central region of the Mandelbrot set opens out so as to
contain a full half-plane, and the entire \[1/2$-limb of the Mandelbrot set
disappears to infinity. 
Also, a number of the tentacles of the Mandelbrot set join together, so
as to enclose new regions. Compare Figure 6.

\midinsert\vskip -.1in
\insertRaster 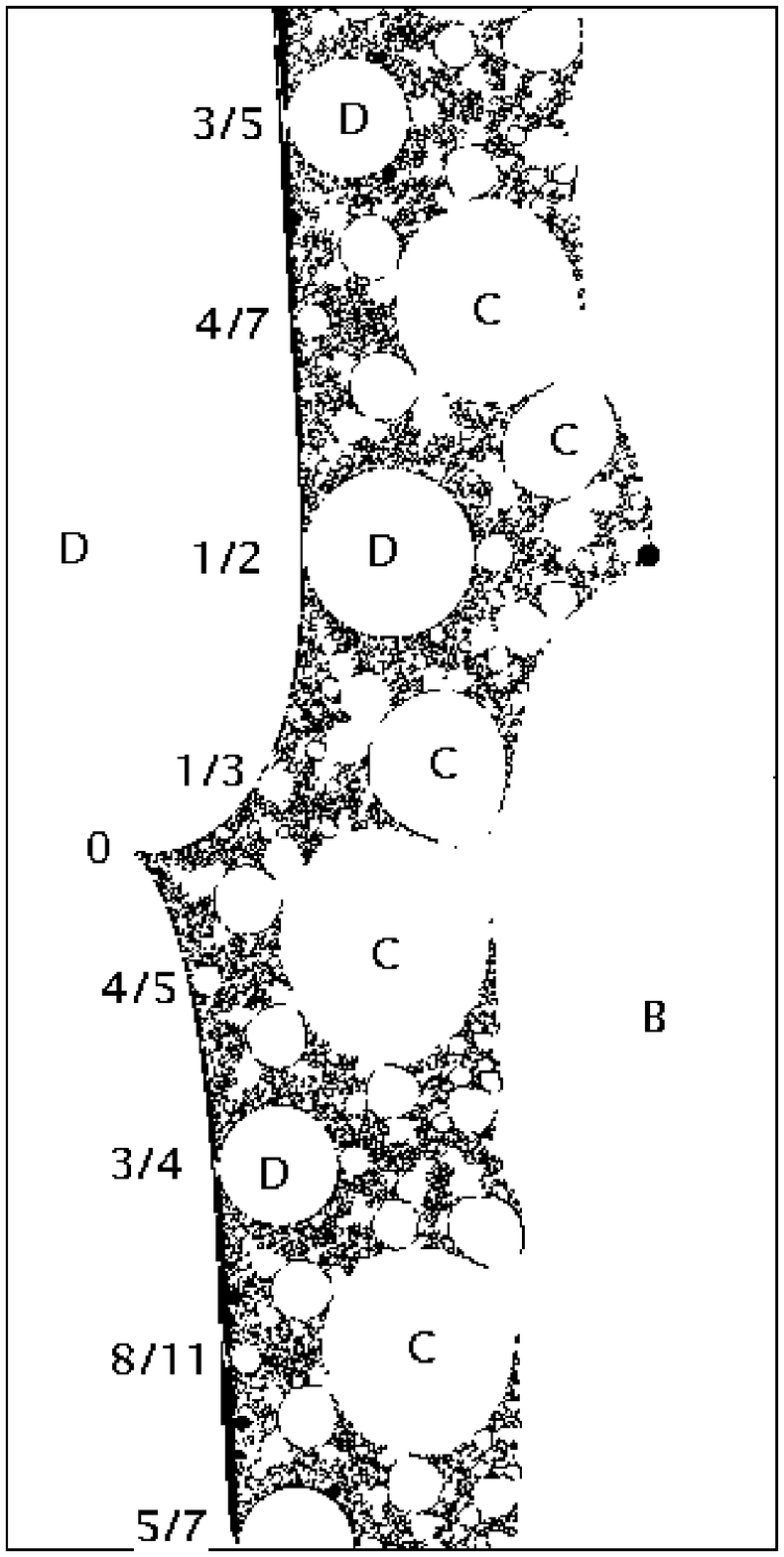 pixels 320 by 640 scaled 450
{\QP\bit Figure 8. Bifurcation locus in \[\Per_1(e^{2\pi i/3})\].
Here the \[2/3$-limb of the Mandelbrot set
has disappeared to infinity. Again various tentacles of
the Mandelbrot set have come together to enclose new regions. Under a slight
perturbation, the fixed point of multiplier \[e^{2\pi i/3}\] will
become repelling but split off an attracting period
three orbit, yielding hyperbolic components of type B, C, and D, as indicated.
(This picture looks rather different from Figure 6 since a different
algorithm was used.) The intersection point of this line with the line
\[\Per_2(-3)\] of Figure 10 is indicated by a heavy dot (upper right).
\par}\vskip -.1in
\endinsert

Similarly, as \[\mu\] tends to
\[e^{-2\pi ip/q}\] non-tangentially,
the entire \[p/q$-limb disappears to infinity. The case \[p/q=2/3\] is
illustrated in Figure 8. This
disappearance is one aspect of Tan Lei's observation that
mating between two quadratic polynomials is never possible when the two
belong to complex conjugate limbs of the Mandelbrot set.

A more precise statement has recently been given by Petersen [Pe].
If \[f\in\Rat_2\] has an attracting fixed point of multiplier \[\mu=\mu_1\]
and a repelling fixed point with combinatorial rotation number \[p/q\]
and multiplier \[\mu_2\], then Petersen shows that a suitable choice of
\[\log(\mu_2)\] is contained in the disk which is
tangent to the imaginary axis at \[2\pi ip/q\], and which contains
\[\log(1/\mu_1)\] as boundary point. Hence, as \[\mu_1\] tends to
\[e^{-2\pi ip/q}\] non-tangentially, it follows that \[\mu_2\] tends to
\[e^{2\pi ip/q}\]. If \[0<p/q<1\], then it follows from (2) that the
multiplier \[\mu_3\] at the third fixed point tends to infinity. {\it Thus, the
entire \[p/q$-limb of the Mandelbrot set must be missing in the locus
\[\Per_1(e^{-2\pi ip/q})\].}\smallskip


{\bf Period 2.}
According to Lemma 3.6, the curve \[\Per_2(\mu)\] is also a straight line
in the \[(\sigma_1\,,\,\sigma_2)$-coordinate plane, given by the equation
\[2\sigma_1+\sigma_2=\mu\]. 
For \[\mu=+1\], this locus \[\Per_2(1)\] coincides with the
locus \[\Per_1(-1)\] of Figure 6. As \[\mu\] varies from \[+1\] to zero,
the right hand region of Figure 6 opens up even further, as shown in Figure 7.
There is some difficulty in finding a good normal form
in the case \[\mu=0\], since the curve \[\Per_2(0)\] includes the
special point \[\langle z\mapsto 1/z^2\rangle\], which has an anomalously
large group of automorphisms. If we place the period 2 critical point at
the origin with its image at infinity, and normalize by an appropriate
scale change, then we get the formula
$$	f(z)~=~{a\over z}+{1\over z^2}~. $$
This normal form is not unique, since  for any cube root of unity \[\lambda\]
the map \[f(\lambda\,z)/\lambda\] will have the same form, but with the
parameter \[a\] replaced by \[\lambda\,a\]. Thus to obtain an invariant
we must pass to \[a^3\], which ranges over \[\C\].
In fact it turns out that \[a^3-6=\sigma_1\].\smallskip

\midinsert
\insertRaster 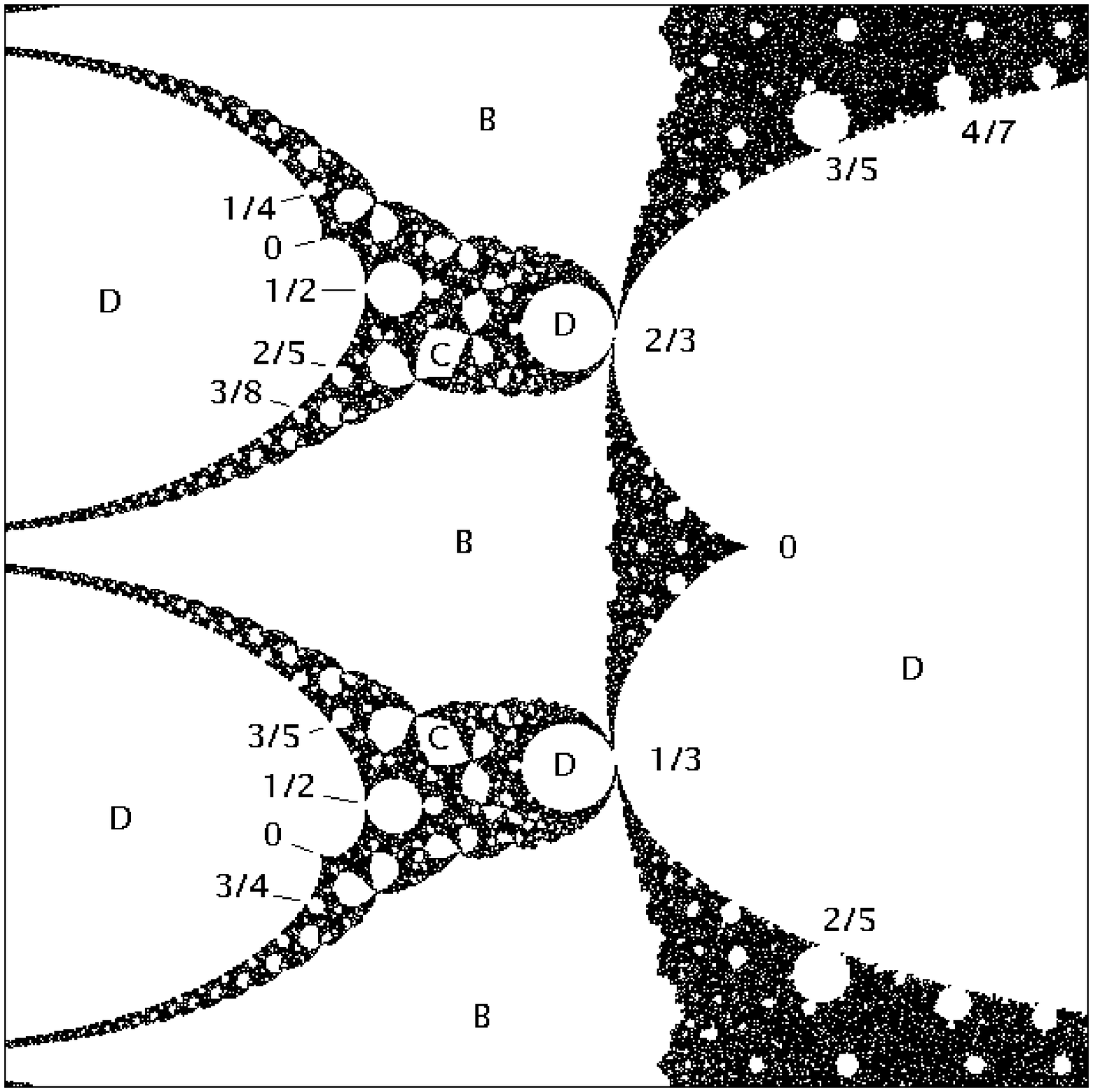 pixels 608 by 608 scaled 450
\smallskip
{\QP\bit Figure 9. Picture in \[\Per_3(0)\], showing
the \[\log(c)$-plane for the family of maps
$$	z~\mapsto~ 1-(1+c)/z+c/z^2 $$
with a critical point of period 3 at \[z=0\]. (The top and bottom are
to be identified, with a little bit of overlap.)
There are only two components of type B in the cylinder \[\C/2\pi i\Z\],
but many components of type C and D. For the three
most prominent components of type D, corresponding to the three period 3
components in the Mandelbrot set, internal angles have been indicated.
Note that the internal angles \[1/3\] respectively \[2/3\] respectively
\[1/2\] correspond to points at infinity.\par}
\endinsert

{\bf Period 3.} The curve \[\Per_3(0)\] also has genus zero, being conformally
isomorphic to a punctured plane or cylinder.
We may put the periodic critical point at the origin, with
orbit \[0\mapsto\infty
\mapsto 1\mapsto 0\]. This yields the unique normal form
$$	z~~\mapsto~~1~-~{1+c\over z}~+~{c\over z^2}~~, $$
where now the parameter \[c\] ranges over \[\C\ssm\{0\}\].
Compare Figure 9, which shows the \[\log(c)$-plane. The lower left hand
portion of this picture can be understood as a perturbation of Figure 8,
and similarly the upper left portion can be understood as a perturbation
of a complex conjugate figure.
In fact, as \[\mu\] varies from \[0\] to \[1-\epsilon\] the bifurcation
locus in \[\Per_3(\mu)\] varies continuously. However, in the limiting case as
\[\mu\to1\] the locus \[\Per_1(\mu)\] degenerates into a union of three
straight lines. (Compare 4.3.) Two of these lines correspond to the locus
\[\Per_1(e^{2\pi i/3})\] and its complex conjugate, while the third line
can be identified with the locus \[\Per_2(-3)\]. Thus, in order
to understand the right half of Figure 9, we must also study the bifurcation
locus in \[\Per_2(-3)\], as shown in Figure 10. Unfortunately, it seems to
require considerable imagination to see how to cut and paste Figure 8,
its complex conjugate, and Figure 10 so as to obtain Figure 9.

For more detailed studies of the curve \[\Per_3(0)\] and of the associated
dynamics, see [R4], [R5], [W].
\smallskip

\midinsert
\insertRaster 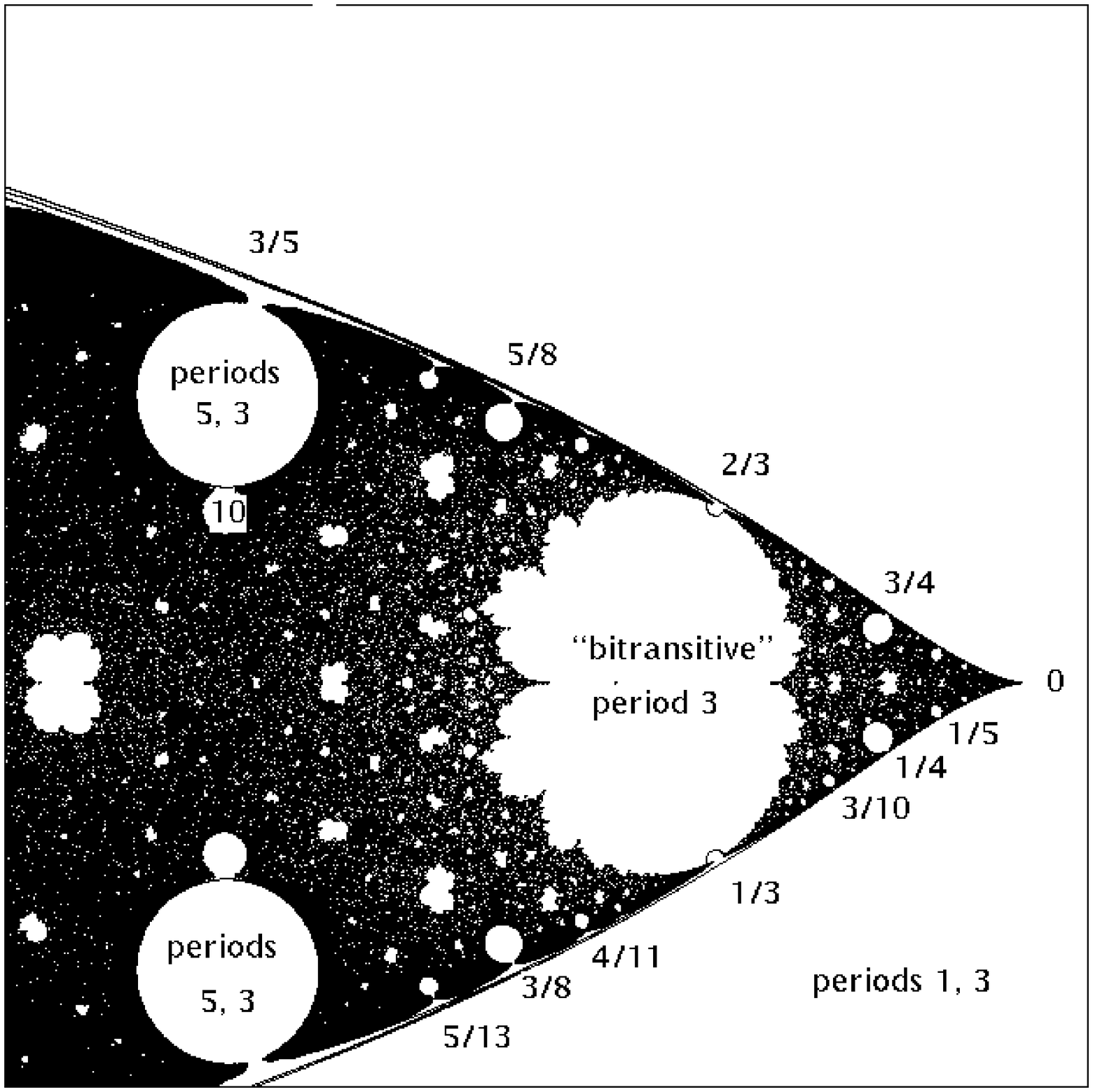 pixels 640 by 640 scaled 450
{\QP\bit Figure 10. Bifurcation locus in the line \[\Per_2(-3)\], which
coincides with the locus of \[\langle f\rangle\] for which \[f\] has a
non-degenerate orbit of period \[3\] and multiplier \[+1\]. (Caution:
The algorithm used to draw this Figure is not too satisfactory.\break In fact,
the boundary of the outer region looks rather ragged
because of very slow convergence.)
\par}
\endinsert

{\bf Period 4.} The curve \[\Per_4(0)\] is isomorphic to a three times
punctured plane. Again we place the periodic critical point at the origin,
and now suppose that it has orbit\break
\[0\mapsto\infty\mapsto 1\mapsto\rho\mapsto 0\]. Here
the parameter \[\rho\] can be described as the
cross-ratio of the points in this orbit. It is not difficult to check
that this invariant \[\rho\] determines the quadratic map uniquely, and
that \[\rho\] ranges over the set \[\C\ssm\{0\,,\;1/2\,,\;1\}\].
The value \[1/2\] is excluded since as \[\rho\to 1/2\] the associated
quadratic map degenerates to a fractional linear map \[z\mapsto 1-(2z)^{-1}\],
and the class \[\langle f\rangle\] tends to the ideal point \[\{i\,,\,-i\,,\,
\infty\}\].\smallskip

\midinsert
\insertRaster 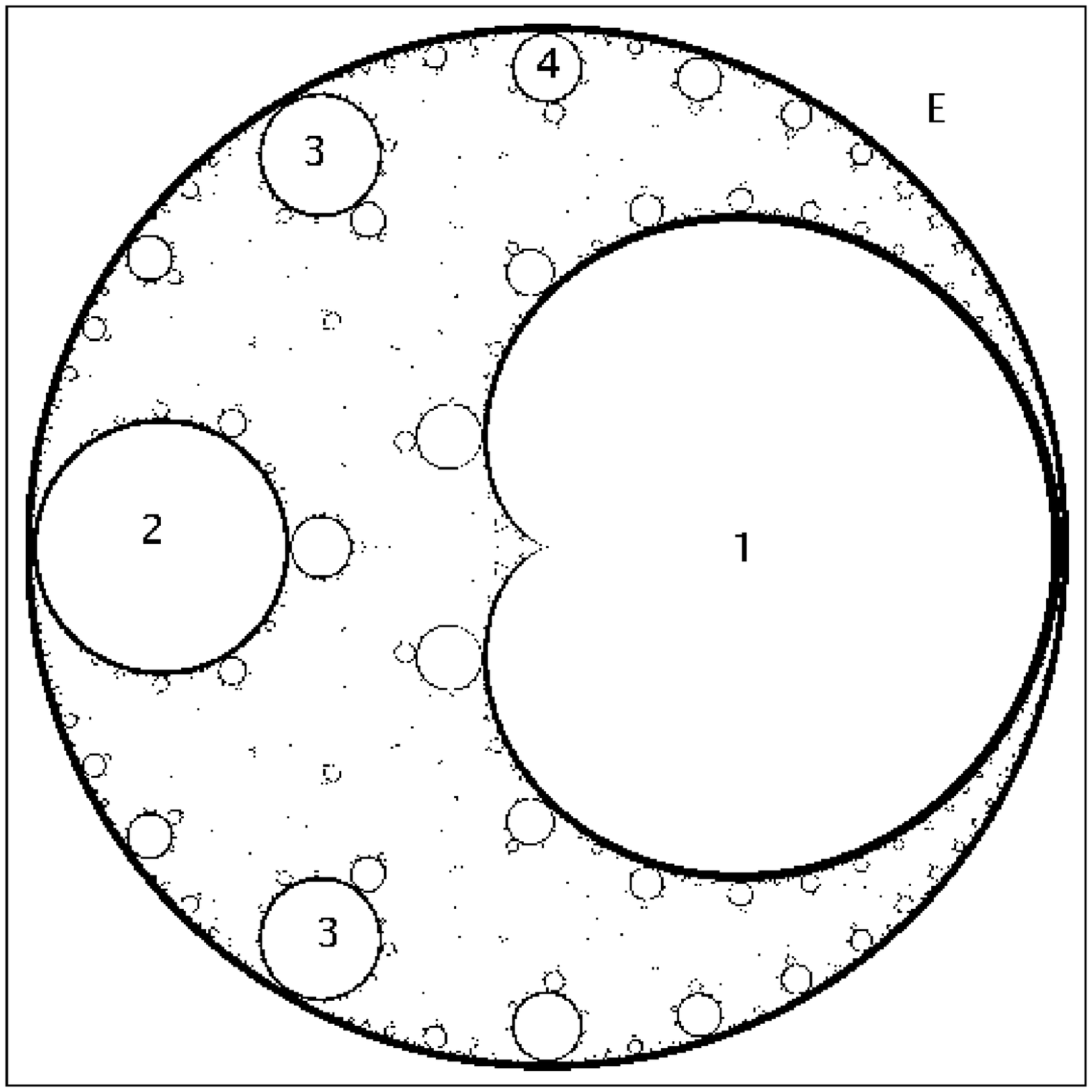 pixels 480 by 480 scaled 450
\smallskip
{\QP\bit Figure 11. Picture for the \[\lambda$-plane for the family of maps
\[z\mapsto\lambda(z+2+z^{-1})\] with pre-periodic critical orbit
\[-1\mapsto 0\mapsto\infty\mapsto\infty\]. In the more prominent regions
with an attracting orbit, the period is indicated.\par}
\endinsert

\vfil\eject
{\bf Pre-periodic critical orbits.}
Another important class of parameter slices are those for which one critical
point is pre-periodic. Since a non-periodic
critical point cannot map directly to a periodic
point when the degree is 2, the simplest possibility is that
the second forward image \[f^{\circ 2}(\omega_1)\]
is fixed. If we place the critical points at \[\pm 1\]
and this fixed point at infinity, this leads to the normal form
$$	w~\mapsto~\lambda\,(w+2+w^{-1})~, $$
with \[-1\mapsto 0\mapsto\infty\]. Note that \[\lambda^{-1}\] is the multiplier
of the fixed point at infinity. A corresponding picture in the \[\lambda$-plane
is shown in Figure 11. Here values of \[\lambda\] outside of the unit disk
correspond to maps in the escape locus \[E\]. However points within the
unit disk correspond to non-hyperbolic maps, since one critical orbit
lands on a repelling fixed point.\smallskip

{\bf Maps with \[J=\hat\C\].}
All of our previous curves have included only maps with at least one
attracting or parabolic orbit. Hence the associated Julia set has always
been a proper subset of the Riemann sphere. However, a map in the present
family\break \[\{w\mapsto \lambda(w+2+w^{-1}\}\] need not have any attracting
or parabolic fixed point. In fact the orbits of {\it both\/} critical
points \[\pm 1\] may eventually
land on repelling cycles. This happens for example for the value \[\lambda_0
=-1/4\], with \[1\mapsto-1\mapsto 0\mapsto\infty\]. The Julia set for
such a map is necessarily the entire sphere. In fact, as
an immediate consequence of a Theorem of Rees [R2] we have
the following much sharper statement: {\it Every neighborhood of
such a point \[\lambda_0\]
contains a set of parameter values \[\lambda\] of positive area
(2-dimensional Lebesgue measure) for which the associated
Julia set is the entire Riemann sphere \[\hat\C\].}

\midinsert\vskip -.1in
\insertRaster 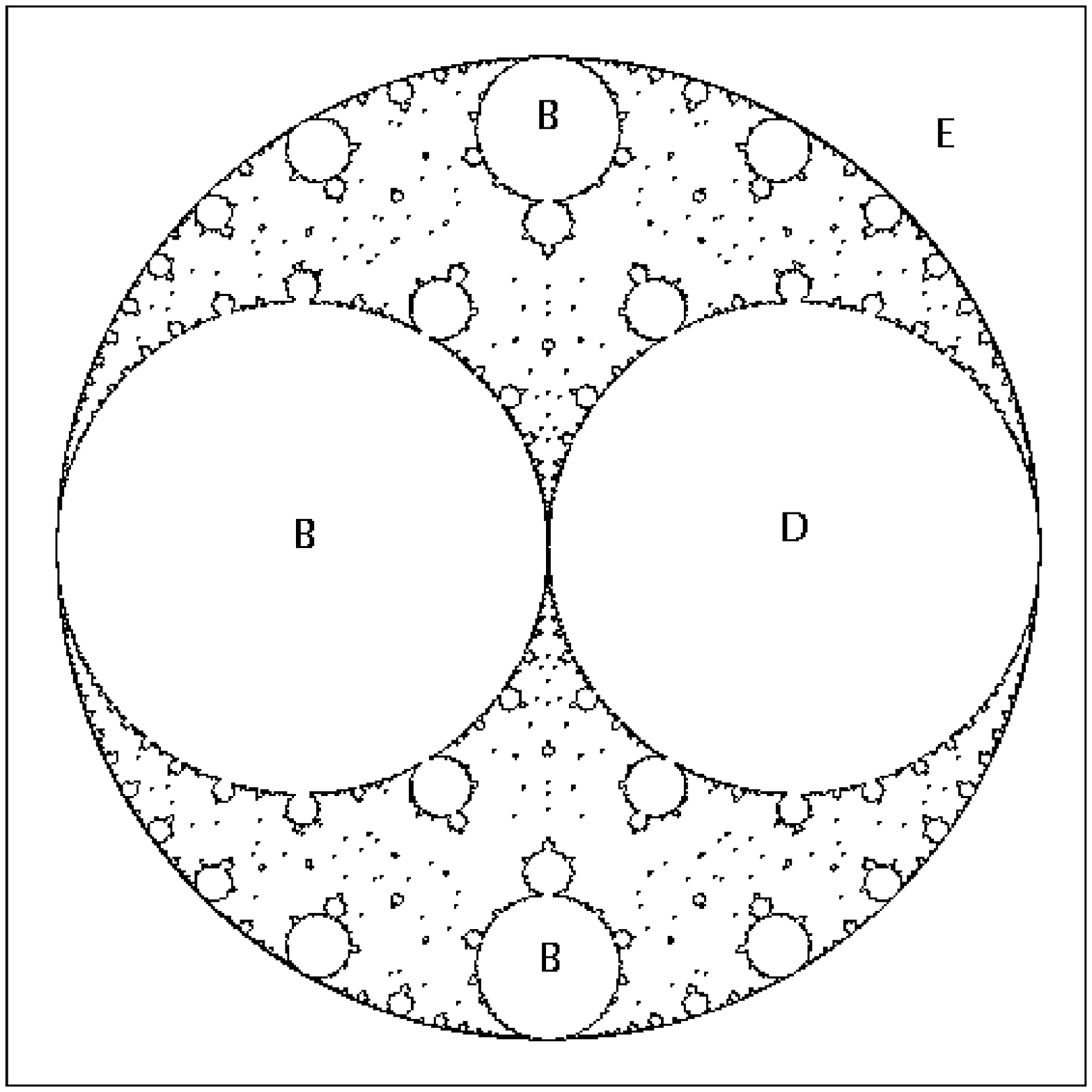 pixels 480 by 480 scaled 450
\smallskip
{\QP\bit Figure 12. Picture of the symmetry locus \[\sl\], represented as
the \[k$-plane for the family of maps \[z\mapsto
k\,(z+z^{-1})\]. The component of type D
on the right is centered at \[k=1/2\], with both critical points fixed,
while the component of type B on the left is centered at \[k=-1/2\], with
critical points \[-1\leftrightarrow 1\].
\par}\vskip -.1in
\endinsert

\vfil\eject
{\bf The symmetry locus.} An apparently quite different example is provided by
the {\bit branch locus\/} or {\bit
symmetry locus\/}, consisting of all conjugacy classes \[\langle f\rangle\]
for which \[f\] has a non-trivial automorphism. By Lemma 5.1 this locus
can be described parametrically as the set of \[\langle f_k\rangle\] with
$$	f_k(z)~=~k\,(z+z^{-1}) ~, $$
where the parameter \[k\] ranges over \[\C\ssm\{0\}\]. A picture of the
\[k$-plane is shown in\break Figure 12. (Compare [HP].)
Note that this figure is centrally
symmetric. In fact, since\break
\[f_k(-z)=-f_k(z)=f_{-k}(z)\] it follows that \[f_k\circ f_k=f_{-k}\circ
f_{-k}\]. Hence the Julia set \[J(f_k)\] is precisely equal to \[J(f_{-k})\],
and the bifurcation diagram in the \[k$-plane, describing different forms
of dynamical behavior, is essentially invariant under the involution
\[k\mapsto-k\]. This does not mean that \[f_k\] has precisely the
same dynamics as \[f_{-k}\]. For example \[f_1\] has a triple fixed point
while \[f_{-1}\] has three distinct fixed points.)
As in the previous example, this symmetry locus
contains maps like \[z\mapsto i(z+z^{-1})/2\] for which both critical
orbits eventually land on repelling cycles
\[(\pm 1\mapsto\pm i\mapsto 0\mapsto\infty)\].

C. Petersen has pointed out to me that this family of maps is
closely related to the locus of Figure 11. In fact the two-to-one
substitution \[w=z^2\] semi-conjugates the family \[z\mapsto k\,(z+z^{-1})\]
to the family
$$	w~\mapsto~\big(k\,(\sqrt w+\sqrt w^{-1})\big)^2~=~\lambda
	(w+2+w^{-1})~, $$
with \[\lambda=k^2\]. Note that symmetric pairs of hyperbolic components in 
the \[k$-plane of\break Figure 12 correspond
to ``sub-hyperbolic components'', where only one critical orbit\break converges
to a periodic attractor, in the \[\lambda$-plane of Figure 11.
\eject

\midinsert
\insertRaster 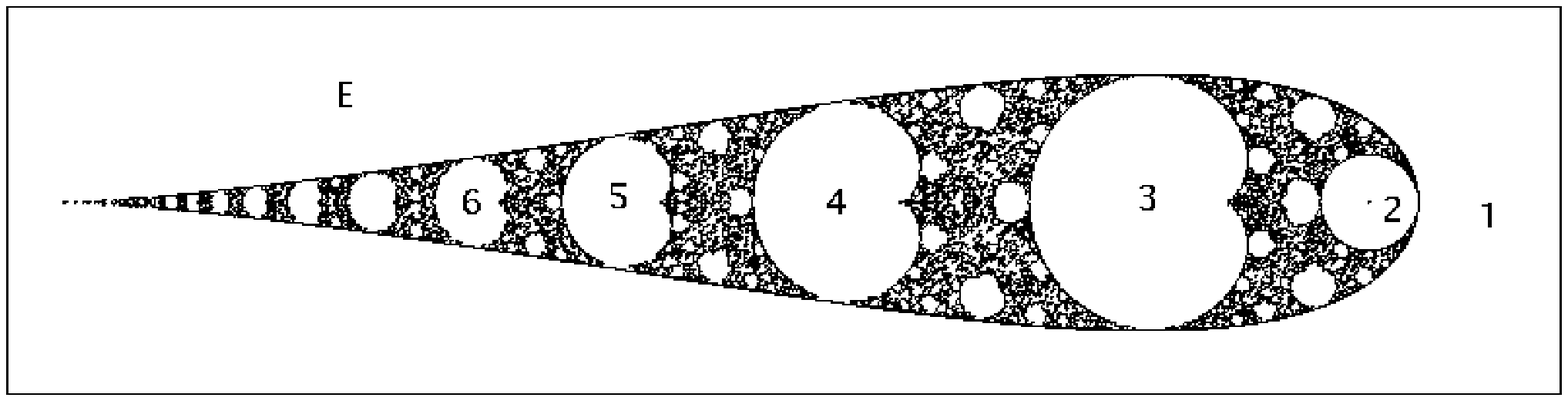 pixels 800 by 200 scaled 450
\smallskip

{\QP\bit Figure 13. \[a^3$-plane for the family of maps \[z\mapsto a+1/z^2\],
with one critical point mapping directly to the other. In the more
conspicuous hyperbolic components (type B or E), the period of the attractor
is indicated.\bigskip}
\endinsert

{\bf Maps for which one critical point maps eventually to
the other.} The simplest possibility here
is that one critical point maps directly to the other. Taking the critical
points to be \[0\to\infty\], we obtain the normal form \[z\mapsto a+1/z^2\]
of Figure 13.

As in the previous two examples, this family contains maps for which
the critical orbits eventually land on a repelling cycle.
This happens for example if \[f(z)= a+1/z^2\] with
\[a^3=-0.5\], so that \[0\mapsto\infty\mapsto a\mapsto-a\mapsto-a\].
Again, it follows from [R2] that there is a set of parameter
values of positive area for which the associated Julia set
is the entire Riemann sphere.
\vfil\eject

\centerline{\bf\S10. Real Quadratic Maps.}\medskip

Consider a quadratic rational map
which has real coefficients, and hence induces a map from the circle
\[\R\cup\infty\] into itself. We will distinguish seven different cases,
depending on the topology of this map from the circle to itself.\smallskip

First suppose that the two critical points of \[f\] are conjugate complex, then
it is easy to check that \[f\] induces a two-to-one covering map from
\[\R\cup\infty\] onto itself. In this case, we must must distinguish between
the positive derivative case, where \[f\] maps this circle with degree \[+2\],
and the negative derivative case with degree \[-2\].\smallskip

On the other hand, suppose that the two critical points are both real.
Then \[f\] maps the entire circle \[\R\cup\infty\]
onto a closed interval \[I=f(\R\cup\infty)\] which is bounded
by the two critical values. Evidently, in
order to study both critical orbits, we need only
study the dynamics of \[f\] restricted to this interval \[I\].

{\bf Definition.} We will say that a real quadratic map with real critical
points is either {\bit  
monotone\/} or {\bit unimodal\/} or {\bit bimodal\/} according as
the interior of this interval \[I=f(\R\cup\infty)\] contains
either no critical points, one critical point, or both critical points.
In the monotone case, we must distinguish
between monotone increasing and monotone decreasing. Similarly, in the
bimodal case we distinguish between \[(+\,-\,+)$-bimodal and
\[(-\,+\,-)$-bimodal according to the pattern of increasing and decreasing
laps.

The unimodal case
is illustrated in Figure 14, where just
one critical point \[\omega_1=0\] belongs
to the interior of the interval \[I=[-1\,,\,1]\]. (Since the other critical
point \[\omega_2=1\] belongs to the boundary of \[I\], we see that this
particular map lies on the boundary between the unimodal and the
\[(+\,-\,+)$-bimodal regions of the ``real moduli space''.)

\midinsert
\centerline{\psfig{figure=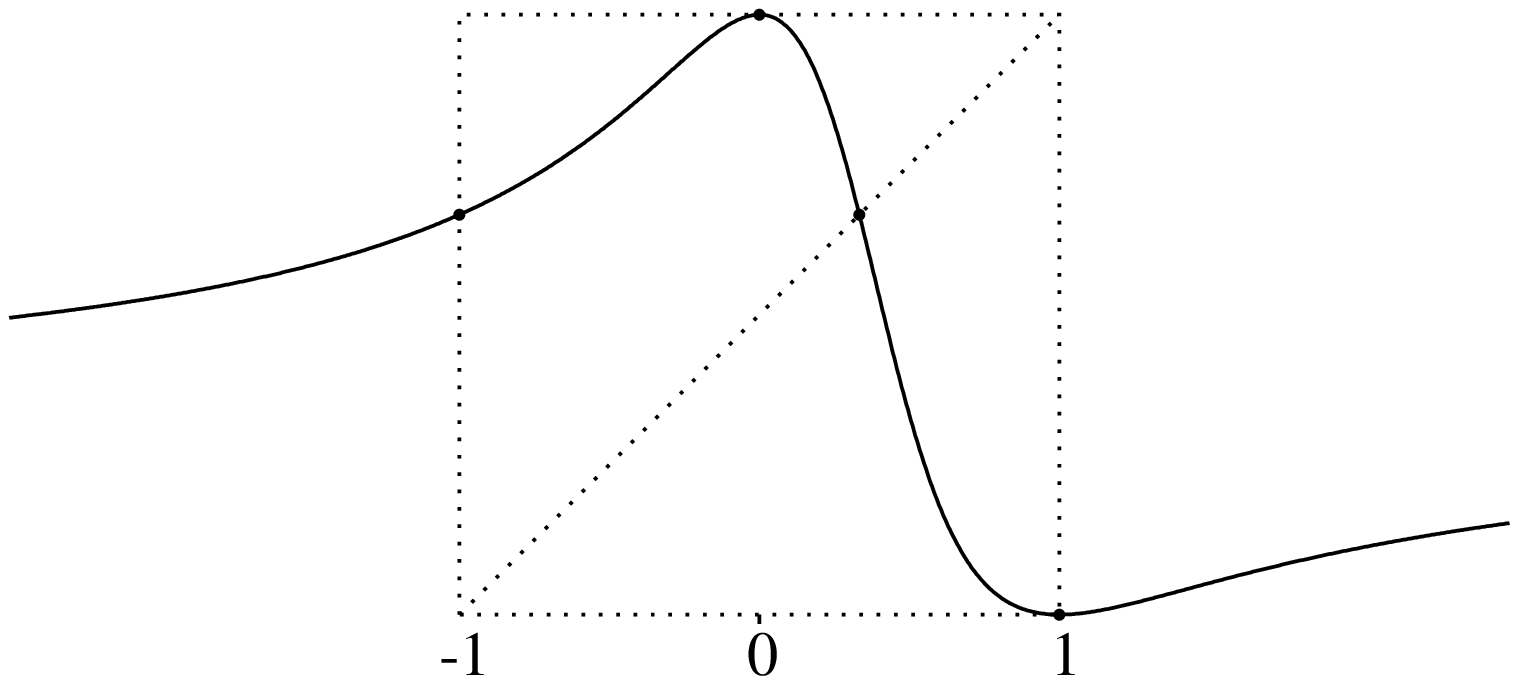,height=2in}}
{\QP\bit Figure 14. Graph of a borderline ``unimodal'' real quadratic map
$$	f_0(x)~=~(1-2x-x^2)/(1-2x+3x^2) $$
with \[f_0(\R\cup\infty)=[-1\,,\,1]\]. In this example,
both critical orbits land on the repelling fixed point \[1/3\], hence the
(complex) Julia set \[J(f_0)\] is the entire Riemann sphere.\par}
\vskip -.1in\endinsert

\midinsert
\centerline{\psfig{figure=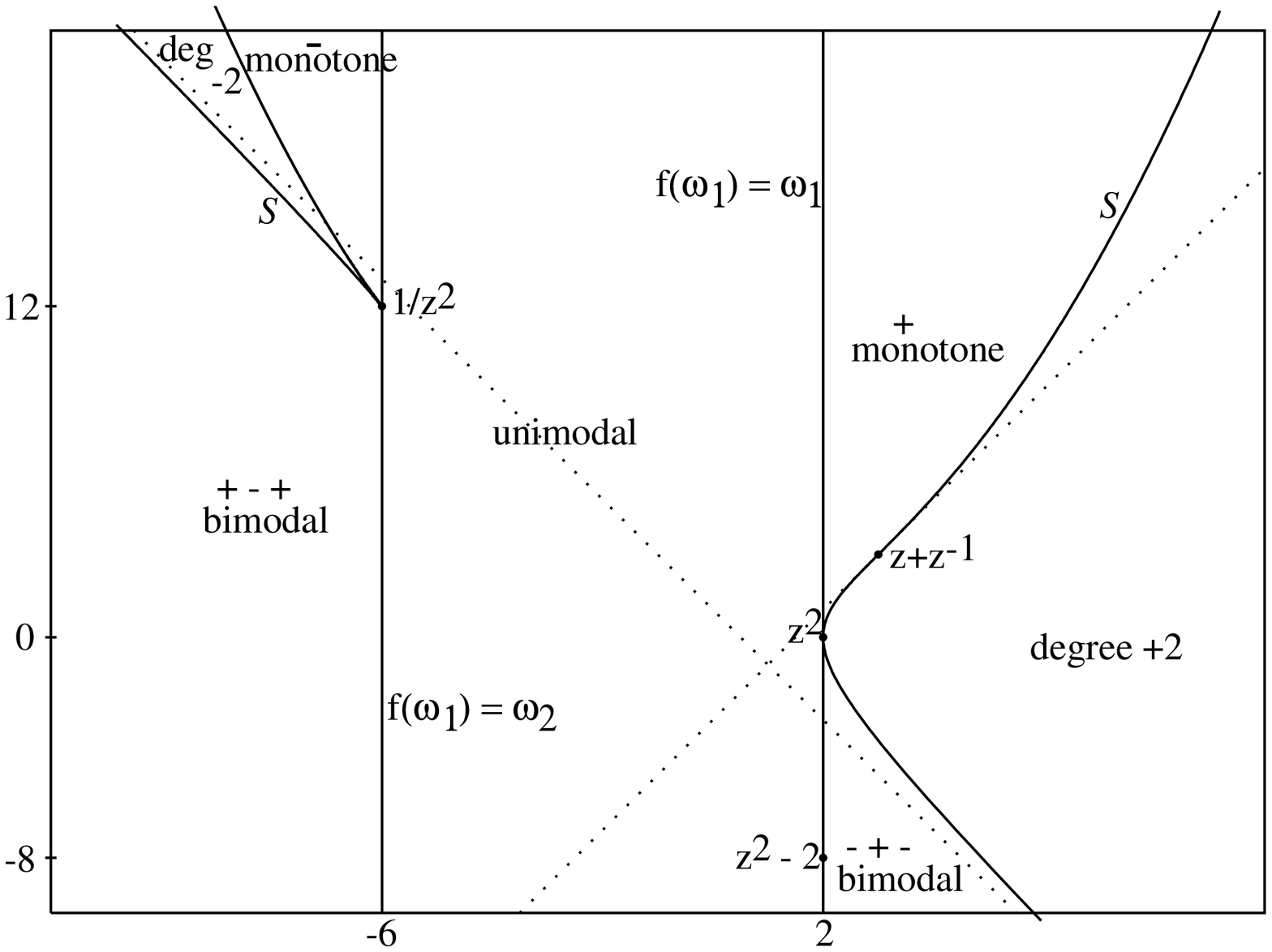,height=3.5in}}
\medskip
{\QP\bit Figure 15. Rectangle
\[[-12,10\times[-10,22]\] in the real \[(\sigma_1\,
,\,\sigma_2)$-plane, showing the seven regions corresponding to
different topological descriptions of \[f\] on \[\R\cup\infty\].
(Drawn to the same scale as Figures 1 and 16.) These seven regions are
bounded by the symmetry locus \[\sl\] together with the vertical lines
\[\sigma_1=-6\] and \[\sigma_1=2\]. The loci
\[\Per_1(\pm 1)\] have been drawn in as dotted lines.\par}
\vskip -.1in
\endinsert

It is not difficult to check that a conjugacy class \[\langle f\rangle\in\M\]
possesses a representative with real coefficients if and only if the two
invariants \[\sigma_1\] and \[\sigma_2\] of \S3 are both real. Hence we define
the {\bit real moduli space\/} to be simply the real
\[(\sigma_1\,,\,\sigma_2)$-plane. In general, the real representative
for \[\langle f\rangle\in\M\] is
unique up to real conjugacy. However, in the special case where \[\langle f
\rangle\] belongs to the real part of the symmetry locus \[\sl\]
of \S5 this is not
true. In fact, for the conjugacy class of \[f(x)=\lambda\,(x+x^{-1})\]
where \[\lambda\in\R\ssm\{0\}\], we can choose \[f\] itself as real
representative,
with critical points \[\pm 1\] and degree zero. But we can also take
\[f(ix)/i=\lambda(x-x^{-1})\] as representative,
with critical points \[\pm i\]
and with degree equal to \[2\,{\rm sgn}(\lambda)\].

The classification of real quadratic maps into seven different
topological types
corresponds to a partition of real moduli space into seven different
regions, as indicated in Figure 15. More precisely,
the real \[(\sigma_1\,,\,\sigma_2)$-plane is cut up into seven pieces
by the two real branches of the symmetry locus \[\sl\], and by the vertical
lines \[\sigma_1=2\] and \[\sigma_1=-6\] which correspond to the classes
\[\langle f\rangle\] for which \[f(\omega_1)=\omega_1\] respectively
\[f(\omega_1)=\omega_2\].

\midinsert
\insertRaster 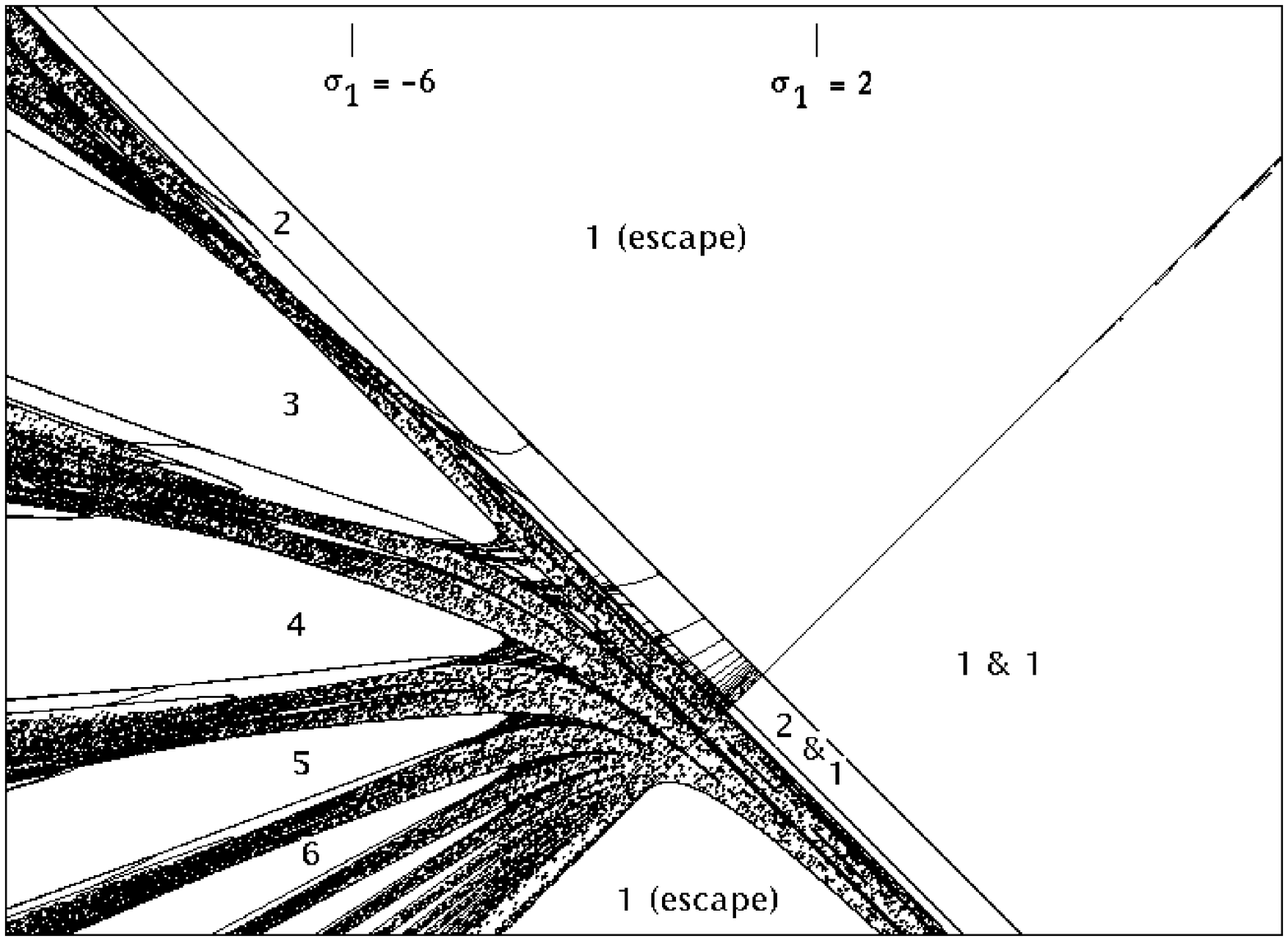 pixels 720 by 524 scaled 450
\smallskip
{\QP\bit Figure 16. Picture of the rectangle \[[-12,10]\times[-10,22]\] in the
real \[(\sigma_1\,,\,\sigma_2)$-plane, emphasizing dynamic behavior.
(Horizontal scale exaggerated.) 
The periods of (possibly complex) attracting orbits
for some of the more prominent hyperbolic components are indicated.
As we traverse any vertical line with
\[-6\le\sigma_1\le 2\], we follow a full family of unimodal maps.
For maps within the top component of the escape locus, there are no
real periodic orbits other than the attracting fixed point. However,
for maps in the bottom component, the entire Julia set is real.
%
\medskip}
\endinsert

\midinsert
\insertRaster 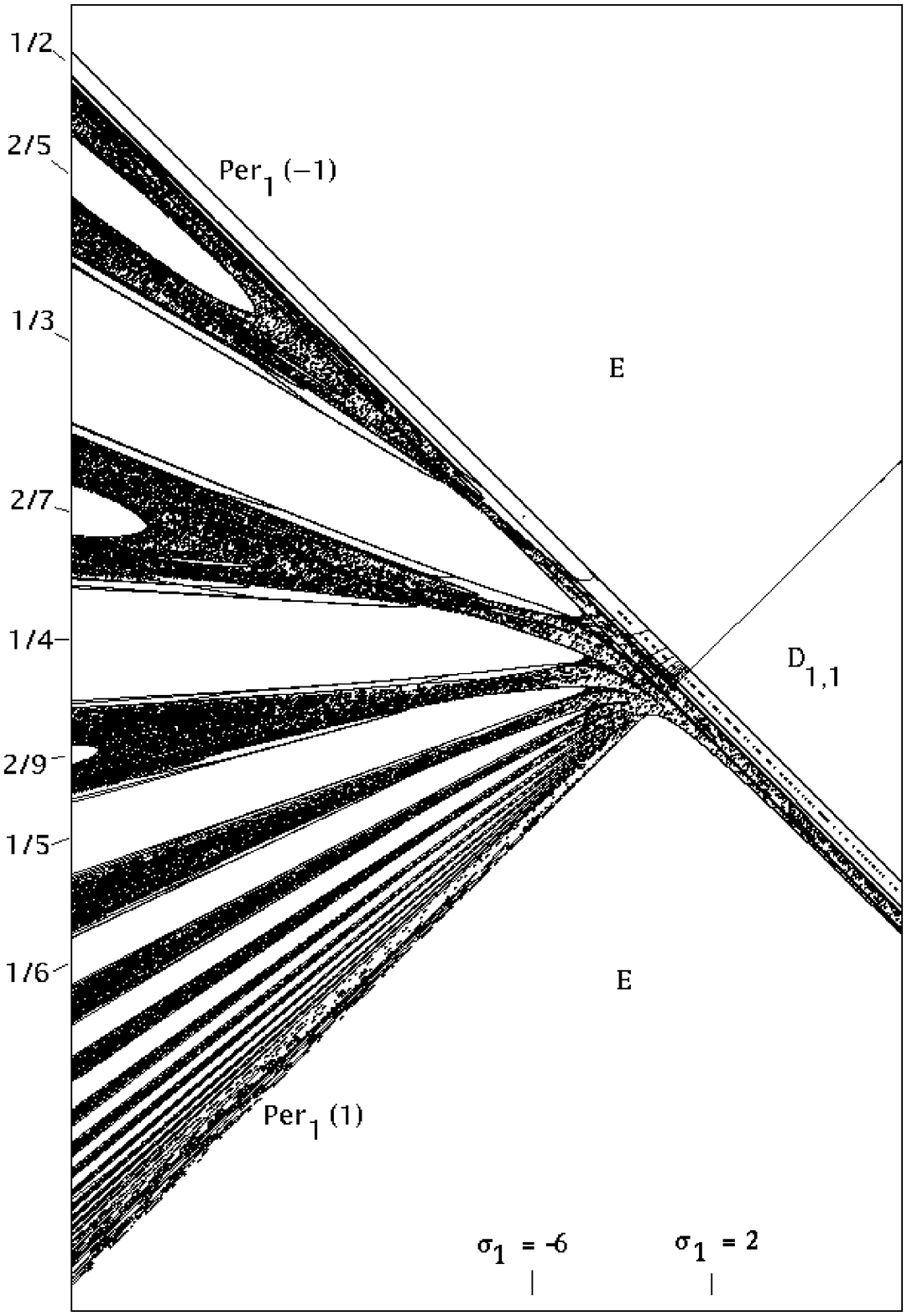 pixels 608 by 880 scaled 450
\smallskip
{\QP\bit Figure 17. Picture of the much larger
rectangle \[[-25,10]\times[-55,55]\] in the
real \[(\sigma_1\,,\,\sigma_2)$-plane. 
For
\[\sigma_1\] very negative, with \[\sigma_2/\sigma_1\approx 2\cos(2\pi ip/q)\]
we are in a component of type B with attracting orbit of period \[q\].
Some representative values of \[p/q\] are indicated along the left margin.
\smallskip}
\endinsert

The bifurcation locus in the real \[(\sigma_1\,,\,\sigma_2)$-plane is
plotted in Figures 16 and 17. Figure 16 shows the same region as Figures 1
and 15. Note that the vertical line
\[\sigma_1=2\] parametrizes the family of real quadratic
polynomials \[x\mapsto x^2+{\rm constant}\], while the line
\[\sigma_1=-6\] parametrizes the family of maps which carry one critical
point to the other --- that is, the real axis of Figure 13. As in Figure 13,
the numbers label hyperbolic components of type B, with attracting period
as indicated. Figure 17 shows a
much larger region, so that we can begin to see the limiting behavior as
\[(\sigma_1\,,\,\sigma_2)\] tends to the line at infinity \[\Per_1(\infty)\].

Significant features of both Figures are the following. The two
straight lines \[\Per_1(1)\] of slope \[+2\] and \[\Per_1(-1)\] of slope \[-2\]
cut the plane into four quadrants:

{\bf Top quadrant.} This is part of the
escape locus. If we restrict attention to maps in this quadrant
with two real critical points, then the real dynamics is completely trivial:
the successive images of \[\R\cup\infty\] converge uniformly to the real
fixed point.

{\bf Right hand quadrant.} These are the maps with two attracting fixed points.
They form the real part of the hyperbolic component \[{\rm D}_{1,1}\]
consisting of all maps with two attracting fixed points. It is centered
at the point \[\langle z\mapsto z^2\rangle\].

{\bf Bottom quadrant.} As we pass from the right hand quadrant to the bottom
quadrant, one fixed point bifurcates to a period two orbit, then to period
four, and
so on, in the usual Sharkovskii pattern, until we again reach the escape
locus.  (Note in fact that the standard quadratic polynomial
family is just the section of this picture with the vertical
line \[\sigma_1=2\].) However, this lower component of the real escape
locus has much more complicated real dynamics, since the entire Julia 
Cantor set is contained in \[\R\cup\infty\].

{\bf Left hand quadrant.} Here we evidently have a much more complicated
situation. There are conspicuous hyperbolic components of type B,
bordered by hyperbolic components of type C. There are also hyperbolic
components of type D, at least within the bimodal region \[\sigma_1<-6\]
(compare Appendix F),
but these are very tiny, and are not visible in the pictures. For
\[\sigma_1\] very negative, with \[\sigma_2/\sigma_1\approx 2\cos(2\pi ip/q)\]
we are in a component of type B with attracting orbit of period \[q\].
Some representative values of \[p/q\] are indicated along the left margin.
For \[p=1\] note that
these are exactly the same as the components which occur in Figure 13.
In fact, the real axis in Figure 13 corresponds exactly to the vertical line
\[\sigma_1=-6\] in Figures 16, 17.\medskip

The real topological classification of Figure 15 can be compared
with the complex dynamic
classification of Figures 16, 17 as follows.\bigskip

\noindent{\bf Monotone case:} If \[f\] is \[+$monotone, clearly
the two critical orbits must converge, either to a single fixed point or to
two distinct fixed points.
If neither fixed point is parabolic, this means that
\[\langle f\rangle\] must
belong either to the hyperbolic component of type \[{\rm D}_{1,1}\] with
two attracting fixed points, or to the escape component E~.
Similarly, in the \[-$monotone case the two critical orbits
must converge to a common fixed point or period two orbit. In the
non-parabolic case, this means that \[\langle f\rangle\] belongs either to
the escape component or to the period two hyperbolic component of type
B~.
\bigskip

\noindent{\bf Unimodal case:} Here there is a wide range of possible
dynamic behavior. In the\break
hyperbolic case, \[\langle f\rangle\] may belong to a component
of any one of the four types B, C, D, or E. One essential restriction is
that there cannot be two
attracting orbits of period greater than one.\bigskip

\noindent{\bf Bimodal case:} Not all bimodal
maps can arise as quadratic maps, since the condition that each point has at
most two pre-images imposes an essential restriction. In particular, the
topological entropy of the map \[I\to I\] can be at most \[\log 2\]. In the
\[(- + -)\] case, there is another quite strong restriction:

{\QP{\bf Lemma 10.1.} \it Every \[(- + -)$-bimodal quadratic map must have
an attracting fixed point.\medskip}

\noindent{\bf Proof.} Clearly such a map must have three fixed points on the
interval \[I=f(\R\cup\infty)\], and two of the three must have negative
multipliers, say \[\mu_1\,,\,\mu_3<0\]. If these two fixed points were
both repelling or parabolic, \[\mu_1\,,\,\mu_3\le -1\], then the formula
\[\mu_2=(2-\mu_1-\mu_3)/(1-\mu_1\mu_3)\] would imply that the third fixed
point also had negative (or infinite) multiplier, which is
impossible.\QED\medskip

\noindent
{\bf Degree $\bf\pm 2$ case:} Here the dynamics is again sharply restricted.

{\QP{\bf Lemma 10.2.} \it If the real quadratic map \[f\] has complex conjugate
critical points, so that it induces a covering
mapping of degree \[\pm 2\] from the circle \[\R\cup\infty\] onto itself, then
the (complex) Julia set \[J(f)\]
is either the circle \[\R\cup\infty\], or a Cantor
set contained in \[\R\cup\infty\].
Both critical orbits converge to attracting orbits of period
\[\le 2\], or to a parabolic fixed point.\medskip}

\noindent The proof is not difficult, and will be omitted.\QED

\end

\def\Rat{{\rm Rat}}
\def\Map{{\rm Map}}
\def\Per{{\rm Per}}
\def\result{{\rm result}}
\def\CP{{\bf CP}}
\def\RP{{\bf RP}}
\def\sl{{\cal S}}
\def\QP{\medskip\leftskip=.4in\rightskip=.4in\noindent}
\def\subarr{\hookrightarrow}
\def\darrow{\rightarrow\!\!\!\!\!\rightarrow}
\def\M{{\cal M}_2}
\def\fm{^{\rm fm}}
\def\cm{^{\rm cm}}
\def\tm{^{\rm tm}}
\ifnum\pageno<2 \pageno=44 \fi

\centerline{\bf Appendix A. Resultant and Discriminant.}\medskip

This will be a brief summary of standard material. Proofs may be
found, for example, in [vW]. Fix two integers \[m\,,\,n\ge 0\].
Let \[p\] range over the vector space \[V_m\] consisting of all
complex polynomial functions
$$ p(z)=a_0\,z^m+\cdots+a_m $$ of degree \[\le m\], and let \[q\] range over
the space \[V_n\] of polynomials $$ q(z)=b_0\,z^n+\cdots+b_n $$
of degree \[\le n\].
The {\bit resultant\/} 
\[\result_{m,n}(p,q)\], can be characterized as the unique
function from \[V_m\times V_n\] to \[\C\] which is bimultiplicative
$$\eqalign{	\result_{m_1+m_2,n}(p_1\cdot p_2\,,\,q)~&=~\result_{m_1,n}
	(p_1\,,\,q)\cdot \result_{m_2,n}(p_2\,,\,q)\cr 
	\result_{m,n_1+n_2}(p\,,\,q_1\cdot q_2)~&=~\result_{m,n_1}
	(p\,,\,q_1)\cdot \result_{m,n_2}(p\,,\,q_2)\cr} $$
and which coincides with the determinant in the degree \[1\,,\,1\] case
$$	\result_{1,1}(az+b\,,\,cz+d)~=~ad-bc~. $$
This function is \[(-1)^{mn}$-symmetric and
bihomogeneous of degree \[(n,m)\], that is
$$ \eqalign{&\result_{m,n}(p,q)~=~ (-1)^{mn}\,\result_{n,m}(q,p)\cr
& \result_{m,n}(\lambda p\,,\,\mu q)~=~\lambda^n\mu^m\,\result_{m,n}(p,q)~.}$$
Note that the subscripts \[m,n\] are essential, since we sometimes consider
polynomials with leading coefficient zero. However, in practice we
will suppress these subscripts whenever they are clear from the context.
The resultant can be expressed as
an \[(m+n)\times(m+n)\] determinant, for example
$$      {\result}(az^2+bz+c\,,\;pz^2+qz+r)\;=\;{\rm det}\left[\matrix{
a& b& c& 0\cr 0& a& b& c\cr p& q& r& 0\cr 0& p& q& r\cr}\right]\;.$$
In the special case \[m=0\] note that \[\result(a_0\,,\,q)=a_0^n\] is
completely independent of \[q\], while
for \[m=1\],
if \[p(z)=z-\xi\] then \[~\result(z-\xi\,,\,q)\;=\;q(\xi)\]. More
generally,
whenever the leading coefficient \[a_0\] is non-zero, we can factor \[p\] as
\[~	p(z)~=~a_0\,(z-\xi_1)\,\cdots\,(z-\xi_m)~, \]
and it follows that
\vskip -.2in
$${\result}(p,q)~=~a_0^n\,\prod_1^m q(\xi_i)~. $$
(On the other hand, if \[a_0=0\], then \[\result_{m,n}(p,q)=(-1)^n\,b_0\,
\result_{m-1,n}(p,q)\].)\break
Similarly, if \[b_0\ne 0\] then we can set
\[q(z)=b_0\,(z-\eta_1)\cdots(z-\eta_n)\] and write
$$      {\result}(p,q)~=~ (-1)^{mn}\,b_0^m\prod_1^n p(\eta_j)
~=~a_0^n\,b_0^m\,\prod_1^m\prod_1^n(\xi_i-\eta_j)~. $$
The resultant may well be non-zero even if one of the leading coefficients
is zero. However, if both \[a_0=0\] and \[b_0=0\] so that
\[{\rm deg}(p)<m\] and \[{\rm deg}(q)<n\],
then \[\result_{m,n}(p,q)=0\].\eject

\noindent { In fact the resultant
\[\result_{m,n}(p,q)\] is zero if and only if either:

$(1)$ the two polynomials \[p\] and \[q\] have a root in common, or

$(2)$ both \[{\rm deg}(p)<m\] and \[{\rm deg}(q)<n\].}
\bigskip

Closely related is the {\bit discriminant\/} of a polynomial of
degree \[m\] or less.
This is a polynomial function \[{\rm discr}_m:V_m\to \C\], which is homogeneous
of degree \[2m-2\],
$$	{\rm discr}_m(\lambda p)~=~\lambda^{2m-2}\,{\rm discr}(p)~. $$
If the leading coefficient \[a_0\] is non-zero,
the discriminant can be defined by the formula
$$      {\rm discr}_m(p)\;=\;(-1)^{m(m-1)/2}\,
	{\result}_{m,m-1}(p,\,p')/a_0\;, $$
where \[p'\] is the derivative. Alternatively, setting
\[~p(z)\,=\,a_0\,(z-\xi_1)\,\cdots\,(z-\xi_m)\], it can be
defined by the formula
$$      {\rm discr}(p)~=~a_0^{2m-2}\,\prod_{i<h}(\xi_i-\xi_h)^2\;. $$
(We omit the subscript \[m\] whenever it is clear from the context.)
Evidently, for such a polynomial with \[a_0\ne 0\], the discriminant
is non-zero if and
only if \[p\] has \[m\] distinct roots. On the other hand, if \[a_0=0\] then
\[	{\rm discr}_m(p)~=~a_1^2\;{\rm discr}_{m-1}(p)~. \]
In particular, for any \[m\ge 2\],
if both of the leading coefficients \[a_0\] and \[a_1\] are
zero, then the discriminant \[{\rm discr}_m(p)\] is zero.
Here are explicit formulas, for the cases \[m\le 3\] which will be of interest.
$$\eqalign{ {\rm discr}_1(az+b)~&=~1~,\cr
	{\rm discr}_2(az^2+bz+c)~&=~b^2-4ac~,\cr
  {\rm discr}_3(az^3+bz^2+cz+d)~&=~b^2c^2-4ac^3-4b^3d-27a^2d^2+18abcd~.\cr}
\eqno (17) $$

\bigskip\bigskip

\centerline{\bf Appendix B. The Space \[\Rat_d\] of Degree \[d\]
Rational Maps.}\medskip

This appendix will describe results from Segal [Se], as well as some
consequences of his work.
For each integer \[d\ge 0\] we can form the space \[\Rat_d\] consisting
of all rational maps \[f(z)=p(z)/q(z)\] where
$$	p(z)=a_0\,z^d+\cdots+a_d\quad{\rm and}\quad  q(z)=b_0\,z^d+\cdots+b_0
\eqno (18) $$
are complex
polynomials with no common root, and with leading coefficients \[a_0\]
and \[b_0\] which are
not both zero. Evidently \[\Rat_d\] is a Zariski open subset of the
complex projective space \[\CP^{2d+1}\] consisting of ratios
\[(a_0:\ldots:a_d:b_0:\ldots:b_d)\]. We will write such ratios
briefly as \[(p:q)\in\CP^{2d+1}\].
More explicitly, this open set can be described as the complement
$$	\Rat_d~\cong~\CP^{2d+1}\ssm W_d $$
where \[W_d\] is the algebraic variety consisting of all ratios \[(p:q)\]
whose resultant\break \[\result(p,q)=\result_{d,d}(p,q)\;\] is
equal to zero. (See Appendix A.)
{\it In particular, it follows easily that  \[\Rat_d\] is a smooth connected
complex manifold of complex dimension \[2d+1\].}

If we map each rational map \[f\in\Rat_d\] to the image \[f(\infty)=a_0/b_0\]
in the Riemann sphere \[S^2=\C\cup\infty\], then we obtain a fibration
$$	\Rat_d^0~\subarr~\Rat_d~\darrow~ S^2 $$
where \[\Rat_d^0\] is the fiber, consisting of rational maps which fix
the point at infinity in \[S^2\]. {\bf Definition.}
Let \[\Map_d\] be the infinite dimensional
space consisting of {\it all\/} continuous
maps from \[S^2\] to itself of degree \[d\], and let \[\Map_d^0\] be the
subspace consisting of maps which fix some base point. Then according to
Segal [Se]:

{\QP{\bf Assertion B.1.} \it The inclusions \[\Rat_d\subset\Map_d\] and
\[\Rat_d^0 \subset\Map_d^0\] induce homotopy equivalences through
dimension \[d\].\medskip}

As model for the fiber \[\Rat_d^0\], Segal uses the space of rational maps
\[f(z)=p(z)/q(z)\] where \[p\] and \[q\] are {\bit monic\/} polynomials of
degree \[d\], so that \[f(\infty)=+1\]. With this normalization, using an
argument due to J. D. S. Jones, he proves the following.

{\QP{\bf Assertion B.2.} \it For \[d\ge 1\] the fundamental group \[\pi_1(
\Rat_d^0)\] is free cyclic. In fact let \[(p,q)\] range over pairs of
monic degree
\[d\] polynomials with no common root. Then the correspondence
\[~p/q\mapsto\result(p,q)\in\C\ssm\{0\}~\] defines a fibration
$$	\Rat_d^0~\to~\C\ssm\{0\} ~,\eqno (19) $$
which induces an isomorphism of fundamental groups.
\medskip}

The fiber \[F_d\] of this fibration
consists of all pairs \[(p,q)\] of monic degree
\[d\] polynomials which have resultant \[\result(p,q)\] equal to \[1\].
Thus \[F_d\] can be considered as an algebraic hypersurface in the
coordinate space \[\C^{2d}\]. As an immediate consequence:

{\QP{\bf Corollary.} \it This fiber \[F_d\] is
simply connected. Furthermore, the universal covering space
\[\widetilde\Rat_d^0\] is homeomorphic to \[F_d\times\C\].\medskip}

Developing these ideas further, we have the following statement. Let
$$	E_d~=~\{(p,q)~:~\result(p,q)=1\}  $$
be the affine hypersurface in \[\C^{2d+2}\]
consisting of {\it all\/} pairs of
polynomials \[(p,q)\] of the form (18) which have resultant \[\result(p,q)\]
equal to \[1\]. Note that there is a natural fibration
$$	F_d~\subarr~E_d~\darrow~\C^2\ssm\{(0,0)\}~, \eqno (20) $$
where the projection map from \[E_d\] to \[\C^2\ssm\{(0,0)\}\] carries each
pair \[(p,q)\] of polynomials with resultant \[1\] to the pair \[(a_0\,
,\,b_0)\] of leading coefficients.
In fact each matrix \[\left[\matrix{a&b\cr c&d\cr}\right]\]
in the group \[{\rm SL}(2,\C)\]  acts on the total space \[E_d\],
carrying each \[(p,q)\] with resultant \[1\] to \[(ap+bq\,,\,cp+dq)\].
Using the fact that \[{\rm SL}(2,\C)\] is generated by elementary matrices,
we see easily that this action preserves the resultant. Clearly this same group
acts on the base space, and the projection map is equivariant. It follows
that we do indeed have a fibration.

{\QP{\bf Assertion B.3.} \it The fundamental group \[\pi_1(\Rat_d)\] is cyclic
of order \[2d\].\break
Furthermore, the universal covering \[\widetilde\Rat_d\]
can be identified with this affine hypersurface \[E_d\].
\smallskip}

{\bf Proof.} Since the resultant of \[p\] and \[q\] is bi-homogeneous of
degree \[(d,d)\] in the coefficients of \[p\] and \[q\], we have the identity
$$	\result(\lambda p\,,\,\lambda q)~=~\lambda^{2d}\,\result(p,q)\;. $$
Thus whenever the resultant of \[p\,,\,q\] is non-zero, so that
\[p/q\] is rational of degree \[d\], there are exactly \[2d\]
choices of \[\lambda\] for which this expression \[\lambda^{2d}\,\result(p,q)\]
is equal to \[+1\]. Evidently the
group of \[2d$-th roots of unity operates freely on the locus \[E_d\],
with quotient space isomorphic to \[\Rat_d\]. Thus \[E_d\] is a \[2d$-fold
covering manifold of \[\Rat_d\]. But \[E_d\] is simply connected, since
it is the total space of a fibration (20) which has simply connected fiber
and simply connected base.\QED
\smallskip

{\bf Proof of 2.1.}
Now let us specialize to the case \[d=2\], and prove Theorem 2.1, as stated
in \S2. Every quadratic rational map has two distinct critical points
\[\omega_1\ne \omega_2\] in the Riemann sphere \[S^2\]. We will use the
notation \[{\cal D\!S}_2(S^2)\] for the {\bit deleted symmetric product\/},
consisting of all unordered pairs of distinct elements in \[S^2\]. Evidently
the correspondence\break \[f\mapsto\{\omega_1\,,\,\omega_2\}\], which assigns
to each \[f\in\Rat_2\] its set of critical points, defines a\break continuous
map \[\Rat_2\to{\cal D\!S}_2(S^2)\].

On the other hand, the group \[\Rat_1\cong {\rm PSL}(2,\C)\] of all
M\"obius transformations of \[S^2\] acts freely on \[\Rat_2\] by left
composition. That is, we have
$$	\Rat_1\times\Rat_2~\to~\Rat_2\qquad{\rm defined~by}\qquad
	(g,f)~\mapsto~ g\circ f~. $$
It is easy to check that two maps in \[\Rat_2\] belong to the same orbit
under this action if and only if they have the same critical points. {\it Thus
our correspondence \[f\mapsto\{\omega_1\,,\,\omega_2\}\] is the projection
map of a principal fiber bundle, having the group \[\Rat_1\cong {\rm PSL}(2,\C)\]
as fiber.}

If we are interested only in the homotopy theory of this bundle, then we may as
well pass to a compact principal sub-bundle
$$	{\rm SO}(3)~\subarr~M^5~\darrow~\RP^2 \eqno (21) $$
which is embedded as a deformation retract. To see this, note that the
projective unitary group \[{\rm PSU}(2)\cong{\rm SO}(3)\] is embedded in
the M\"obius group \[\Rat_1\] as a deformation retract, and that the projective
plane \[\RP^2\] consisting of pairs of antipodal points in \[S^2\] is
embedded in \[{\cal D\!S}_2(S^2)\] as a deformation retract. It is not
difficult to check that the corresponding
total space \[M^5\] is indeed a smooth manifold,
embedded in \[\Rat_2\] as a\break deformation retract.\smallskip

Such \[{\rm SO}(3)$-bundles over \[\RP^2\] are classified by homotopy
classes of mappings\break
\[\RP^2\to {\rm BSO}(3)\], or (using obstruction theory or Hopf's Theorem)
by cohomology classes in the group
$$ H^2\big(\RP^2\,;\;\pi_1({\rm SO}(3))\big)~\cong~\Z/2\;. $$
In fact the bundle
\[(21)\] must be the non-trivial \[{\rm SO}(3)$-bundles over \[\RP^2\].
For if \[M^5\] split as
a product, then the fundamental group \[\pi_1(M^5)\cong\pi_1(\Rat_2)\] would
be the direct sum of two cyclic groups of order 2. But according to B.3
this group is actually cyclic of order 4.\QED

As one immediate consequence of 2.1, we see that  \[\Rat_2\] has the
rational homology of a 3-sphere. As another easy consequence:

{\QP{\bf Corollary B.4.} \it The 2-fold orientable covering manifold of \[M^5\]
is homeomorphic to the product \[{\rm SO}(3)\times S^2\]. Hence the universal
covering of \[M^5\] is homeomorphic to \[S^3\times S^2\].\medskip}

{\bf Proof.} Geometrically, this 2-fold covering can be described as the
space of triples\break \[(f\,,\,\omega_1\,,\,\omega_2)\] consisting 
of a quadratic
map in \[M^5\] with {\bit marked critical points\/}. (Compare \S6.)
The statement
that the lifted \[{\rm SO}(3)$-bundle over \[S^2\] splits as a product
means that there exists some cross-section \[\omega\mapsto
f_\omega\] which assigns to each point \[\omega\in S^2\] a quadratic rational
map having \[\omega\] as critical point. This follows from an easy cohomology
argument, or can be constructed
explicitly as follows. For each \[\omega\] in the
finite plane \[\C\], consider the M\"obius transformation
\[\;	g_\omega(z)~=~(\bar\omega z+1)/(-z+\omega) \;.\]
This belongs to the compact subgroup \[{\rm SO}(3)\subset\Rat_1\]. It
carries \[\omega\] to \[\infty\] and carries
the antipodal point \[-1/\bar\omega\] to zero. The composition
\[f_\omega(z)=\phi_{\omega^2}^{-1}(\phi_\omega(z)^2)\] is then a quadratic
rational map in \[M^5\]
having \[\omega\] and \[-1/\bar\omega\] as critical points.
This construction has been chosen so that \[f_\omega(z)\] tends to a well
behaved limiting value \[f_\infty(z)=z^2\] as \[\omega\to\infty\].\QED

{\bf Remark.} Of course such a correspondence \[\omega\mapsto f_\omega\],
which assigns to each point \[\omega\in S^2\] a rational map \[f_\omega\in
\Rat_2\] having \[\omega\] as critical
point, cannot possibly be holomorphic. For the associated
correspondence \[\omega\mapsto\omega'\] which assigns to \[\omega\] the
{\it other\/} critical point of \[f_\omega\] has degree \[-1\],
and hence is non-holomorphic.
\bigskip\bigskip

\centerline{\bf Appendix C. Normal Forms and Relations Between Conjugacy
Invariants.}\medskip
This section will study three convenient normal forms for quadratic
rational maps. First, as in the proof of 3.1, we consider the following.
\smallskip

{\bf Fixed Point Normal Form.} Let
$$      f(z) = z{z+\mu_1\over \mu_2 z+ 1}~, \eqno (22)$$
with \[\mu_1\mu_2\ne 1\]. Here the origin is
a fixed point of multiplicity \[\mu_1\], and infinity is a fixed point
of multiplicity \[\mu_2\]. The third fixed point lies at
\[(1-\mu_1)/(1-\mu_2)\], and has multiplier \[\mu_3=(2-\mu_1-\mu_2)/
(1-\mu_1\mu_2)\]. The critical points for this map are the two roots
\[\omega=(-1\pm\sqrt{1-\mu_1\mu_2})/\mu_2\] of the equation
$$      \mu_2\,\omega^2+2\,\omega+\mu_1~=~0~, $$
(with one critical point at infinity when \[\mu_2=0\]). A
brief computation shows that the corresponding critical values
are given by the equation
$$      f(\omega)~=~-\omega^2~. $$
\smallskip

{\bf Mixed Normal Form.} Now consider the normal form
$$      f(z) ~=~ {z + z^{-1} \over\mu}+a \eqno (23)$$
with critical points \[\pm 1\] and with a fixed point of multiplier
\[\mu\ne 0\] at infinity. (Compare 5.1 and 8.3, as well as
Figures 5, 6, 12, 13.) Here the critical values are \[f(\pm 1)=
a\pm 2/\mu\]. The two finite fixed points \[z_1\] and \[z_2\]
can be described as the roots of the equation
$$      (1-\mu)\,z_i^2+a\,\mu\, z_i+ 1~=~0~. $$
Either by direct computation, or using C.1 below, one finds that the
corresponding\break multipliers \[\mu_i=f'(z_i)=(1-z_i^{-2})/\mu\] satisfy
the relation
$$ \sigma_1~=~\mu_1+\mu_2+\mu~=~\mu(1-a^2)-2+4/\mu~.$$
As in 3.2,
the symmetric function \[\sigma_2=\mu_1\mu_2
+\mu_1\mu_3+\mu_2\mu_3\] can then be computed from the equation
\[0=\mu^3-\sigma_1\mu^2+\sigma_2\mu-\sigma_3=\mu^3-\sigma_1\mu^2+\sigma_2\mu
-\sigma_1+2\], hence
$$      \sigma_2=(\mu+\mu^{-1})\sigma_1-(\mu^2+2/\mu)~. $$
The parameter \[\tau=\mu_1\mu_2\] along the line \[\Per_1(\mu)\]
is equal to \[\sigma_3/\mu=(\sigma_1-2)/\mu\],
hence
$$      \tau~=~\Big({2-\mu\over\mu}\Big)^2-a^2~. \eqno (24) $$
(Compare 6.8.) {\it Thus, for \[\mu\ne 0\],
the parameter \[\tau\] along the line \[\Per_1(\mu)\]
coincides with \[-a^2\], up to translation by a constant depending only
on \[\mu\].}
\smallskip

{\bf Critical Point Normal Form.}
If we put the critical points at zero and infinity, so that
\[f(z)=f(-z)\], then
$$      f(z)\;=\; {\alpha z^2+\beta\over \gamma z^2+\delta}~.\eqno(25) $$
This normal form is particularly convenient for computations and for Julia
set pictures. (See Figures 2, 3, 4. In fact Figures 8, 10, 16, 17 were made
by first translating into these coordinates, using C.4.)
As in 6.1, the quantities
$$	A={\alpha\,\delta\over\alpha\,\delta-\beta\,\gamma}={1+\beta\,\gamma
\over\alpha\,\delta-\beta\,\gamma}~,\qquad B={\alpha^3\beta\over
 (\alpha\,\delta-\beta\,\gamma)^2}~,\qquad C={\gamma\,\delta^3\over
(\alpha\,\delta-\beta\,\gamma)^2} $$
are invariant under holomorphic conjugacies which fix the two critical
points. Hence \[A\] and \[\Sigma=B+C\] are holomorphic conjugacy invariants.
As noted in 6.3, we can use either the invariants \[\sigma_1\,,\,
\sigma_2\] of \S3 or these invariants \[A\] and \[\Sigma\] as coordinates
on the moduli space \[{\cal M}_2\], in order to show
that \[{\cal M}_2\] is isomorphic to \[\C^2\]. The coordinates \[A,\,\Sigma\]
are apparently easier to compute, since we need only
solve a quadratic equation to find the critical points and reduce to
the normal form (25), while we must solve a cubic equation to find the
fixed points. In fact we can always compute the invariants without solving
any polynomials equations, and
it is not too difficult to compute one pair of invariants
in terms of the other.\smallskip

We will first prove the following. Let \[\omega_1\,,\, \omega_2\] be the critical
points of \[f\], and let \[v_i=f(\omega_i)\] be the corresponding critical values.
The relative position of these four points is determined by the
cross-ratio
$$	\rho(f)\;=\;{(v_1-\omega_1)\,(v_2-\omega_2)\over(\omega_1-\omega_2)
\,(v_1-v_2)}\;, \eqno (26) $$
which is always a finite complex number, since the denominator cannot vanish.
Note that \[\rho=0\] if and only if one of the critical points is fixed by
\[f\], so that \[f\] is conjugate to a polynomial map. Similarly, it is
not hard to check that \[\rho=1\] if and only if one critical point maps
to the other.

{\QP{\bf Lemma C.1.} \it This cross-ratio \[\rho\] is related to the invariants
\[\sigma_i\] of Lemma 3.1 and to the invariant \[A\] of Lemma 6.1 by the linear
equation
$$	\rho(f)~~=~~-\sigma_3/8~~=~~(2-\sigma_1)/8~~=~~1-A~. $$}

{\bf Proof outline.} If we use the normal form (25), then
$$ \omega_1=0\;,\quad \omega_2=\infty\qquad{\rm and}\qquad v_1=\beta/\delta\;,
\quad v_2=\alpha/\gamma\;,$$
so the cross-ratio (26) reduces to the form \[v_1/(v_1-v_2)=-\beta\gamma
=1-A\].
On the other hand, if we use the normal form (4), then computation shows
that the critical points and critical values are given by
$$ \omega_\pm~=~\big(-1\pm\sqrt{1-bc}\;\big)/c\qquad{\rm and} \qquad f
	(\omega_\pm)	=-\omega_\pm^2~, \eqno (27)$$
with \[b=\mu_1\,,~c=\mu_2\]. The cross-ratio (26)
then works out to be \[bc\,(b+c-2)/(8-8bc)\].
Using (3) and $(2)$ this reduces to \[-\sigma_3/8=(2-\sigma_1)/8\],
as required.\QED

In order to relate the remaining invariants, we will make use of the
discriminant of a polynomial, and the resultant of two polynomials.
(Compare Appendix A.) Note that the three fixed points of a quadratic
rational function \[f(z)=p(z)/q(z)\] are
the roots of the polynomial equation \[z\,q(z)-p(z)=0\]. (This equation
is understood to have a ``root at infinity'' if and only if its degree is
strictly less than 3.) The {\bit discriminant
\/} of this equation is a homogeneous function of degree 4 in the coefficients
of \[p\] and \[q\], which vanishes if and only if there is a multiple root.
In order to make this discriminant into an invariant of the rational function
\[p(z)/q(z)\], we must divide by the resultant of \[p\] and \[q\], which is
also a homogeneous function of degree 4 in the coefficients of \[p\] and \[q\].
This resultant is never zero, since these two polynomials have no common root.

{\QP{\bf Lemma C.2.} \it For any quadratic rational function \[f=p/q\], the
{\bit discriminant ratio}
$$	{\rm drat}(f)\;=\;{{\rm discr}(zq-p)\over{\result}(p\,,\,q)}
\eqno (28) $$
is a conjugacy invariant
which can be expressed in terms of the invariants\break \[\sigma_i\] of \S3
by the formula
$${\rm drat}(f)\;=\;
\prod(\mu_i-1)\;=\;\sigma_3-\sigma_2+\sigma_1-1\;=\;
2\sigma_1-\sigma_2-3\;, $$
and in terms of the invariants \[A\] and \[\Sigma=B+C\] of \S4 by the formula
$$	{\rm drat}(f)\;=\;-8A^2+36A-4\Sigma-27\;. $$
\smallskip}

{\bf Proof.} Clearly this ratio is a well defined complex number which
depends only of the function \[f(z)=p(z)/q(z)\]. To show that it is a conjugacy
invariant, first let us try replacing \[f(z)\] by \[f(z+c)-c\]. Then the
polynomials \[p\] and \[q\] will be replaced by \[p(z+c)-cq(z+c)\] and
\[q(z+c)\] respectively. Similarly the polynomial \[F(z)=zq(z)-p(z)\]
will be replaced by \[F(z+c)\]. It is then easy to check that both the
numerator \[{\rm discr}(F)\] and the denominator \[\result(p,q)\] remain
invariant. Similarly, if we replace \[f(z)\] by \[1/f(1/z)\], then we must
replace \[p\,,\;q\] and \[F\] by \[z^2q(1/z)\,,\;z^2p(1/z)\] and \[-z^3F(1/z)\]
respectively. Again it is easy to check that both \[{\rm discr}(F)\] and
\[\result(p,q)\] remain invariant. Finally, if we change scale, replacing
\[p(z)\,,\;q(z)\] and \[F(z)\] by \[p(\lambda z)\,,\;\lambda q(\lambda z)\]
and \[F(\lambda z)\], then it is not hard to check that both numerator
and denominator will be multiplied by \[\lambda^6\]. This proves that
the ratio (28) is a holomorphic conjugacy invariant.
\smallskip

In both cases, we can now simply make a direct computation.
Thus if \[p(z)=\alpha z^2+\beta\] and \[q(z)=\gamma z^2+\delta\]
as in (25), then it is easy to check that
$$	{\rm resul}(p,q)\;=\;(\alpha\delta-\beta\gamma)^2\;=\;1\;. $$
The classical formula (17) for the discriminant of
\[zq(z)-p(z)=\gamma z^3-\alpha z^2+\delta z-\beta\] is
$$	{\rm discr}(\gamma z^3-\alpha z^2+\delta z-\beta)\; =\; -27\beta^2\gamma^2
-4\delta^3\gamma+18\alpha\beta\gamma\delta+\alpha^2\delta^2-4\alpha^3\beta\;.
 $$
(See [vW].) In terms of the invariants (13) we can easily reduce this
to the expression
$$	{\rm drat}(f)\;=\;{\rm discr}(zq-p)\;=\;-8A^2+36A-4(B+C)-27\;, $$
as required.\smallskip

Similarly, direct computation using the normal form (4) together with (3)
shows that the invariant \[{\rm drat}(f)\] is equal to \[\prod(\mu_i-1)\].\QED

{\bf Remark C.3.} Lemma C.1 could also be expressed by
an analogous explicit formula
$$	\sigma_3\;=\;{\result
(zq-p\,,\,p'q-q'p)\over\result(p,q)^2}\;.\eqno (29) $$
\smallskip

Combining Lemmas C.1 and C.2, we obtain the following.

{\QP{\bf Corollary C.4.} \it The invariants \[\sigma_i\,,\;A\] and \[\Sigma=B+C\]
are related by the identities
$$	\sigma_1\;=\;8A-6\;,\qquad \sigma_2\;=\; 8A^2-20A+
4\Sigma+12\;\, . \eqno (30) $$\smallskip}

\bigskip\bigskip\eject

\centerline{\bf Appendix D. The Geometry of Periodic Orbits.}\medskip

First some elementary number theory. Let \[\nu_d(n)\] be the number of periodic
points of period \[n\] for a generic polynomial map \[f:\C\to\C\]of
degree \[d\ge 2\]. Since the fixed points of \[f^{\circ n}\], that is the roots
of the equation \[f^{\circ n}(z)-z=0\] of degree \[d^n\], are precisely
the periodic points with period \[m\] dividing \[n\], it follows that
$$	d^n~=~\sum_{m|n}\,\nu_d(m)~, $$
to be summed over all numbers \[m\ge 1\] which divide \[n\]. We can easily
solve inductively for \[\nu_d(n)\]. For example
$$	\nu_d(1)=d\;,\quad\nu_d(2)=d^2-d\;,\quad\nu_d(3)=d^3-d\;,\quad\nu_4
(d)=d^4-d^2~.$$
(Using the ``M\"obius Inversion Formula'', the solution can be written as
$$	\nu_d(n)~=~\sum_{m|n}\,\mu(n/m)\,d^m~, $$
where the M\"obius function \[\mu(\;)\] is defined by
\[\mu(p_1\cdots p_k)=(-1)^k\] for a product of \[k\] distinct primes with
\[k\ge 0\], and \[\mu(m)=0\] if \[m\] is not a product of distinct primes.
See for example [HW].)
It is easy to check that this number \[\nu_d(n)\] is always divisible
by both \[d\] and \[n\]. Henceforth, let us specialize to the case \[d=2\].
Thus the number of period \[n\] orbits for a generic
quadratic polynomial map is
\[\nu_2(n)/n\]. According to 4.2, the number of period
\[n\] hyperbolic components in the Mandelbrot set is \[\nu_2(n)/2\].
Here is a table.\medskip

\centerline{\vbox{\hrule\hbox{\vrule\vbox{
{\halign{\quad # & \quad# & # & # & # & # & # & # & # & # & # & # & #\cr
\noalign{\smallskip}
$n$ & 1 & 2 & 3 & 4 & 5 & 6 & 7 & 8 & 9 & 10 &11&12\cr
$\nu_2(n)/n$ & 2 & 1 & 2 & 3 &  6 &  9 & 18 &  30 &  56 &  99 &186&335\cr
$\nu_2(n)/2$ & 1 & 1 & 3 & 6 & 15 & 27 & 63 & 120 & 252 & 495&1023&2010 \cr
\noalign{\smallskip}}
}}\vrule}\hrule}}
\medskip

Now let us return to the study of quadratic rational maps.
To simplify the discussion, let us choose some normal form
for which the point at infinity cannot be periodic. For example, any
conjugacy class contains a representative of the form
$$	f(z)~ = ~{1\over a(z+z^{-1})+b} ~, $$
with \[\infty\mapsto 0\mapsto 0\], so that \[\infty\] is not periodic.
Writing \[f^{\circ n}(z)\] as a quotient \[p_n(z)/q_n(z)\] of two polynomials,
where \[~{\rm degree}(q)=2^n\ge{\rm degree}(p)\], 
it follows that the equation
$$	z\,q_n(z)-p_n(z)~=~0 $$
for fixed points of \[f^{\circ n}\] has degree exactly \[2^n+1\]. Using
an argument similar to the proof of 4.2, we see
that there are unique polynomials \[\Phi_n(z)\] so that
$$	z\,q_n(z)-p_n(z)~=~\prod_{m|n}\Phi_m(z)~~~\Longleftrightarrow~~~
\Phi_n(z)~=~\prod_{m|n}\,(z\,q_n(z)-p_n(z))^{\mu(n/m)}~. $$
Here \[\Phi_n(z)\] has degree \[\nu_2(n)\] for \[n>1\], and \[\Phi_1(z)\]
has degree \[3\]. To simplify the notation, let us write \[N=\nu_2(n)/n\]
provided that \[n>1\], and \[N=3\] if \[n=1\].
By definition, the \[Nn\] roots of \[\Phi_n(z)\],
counted with multiplicity, will be called the (formal) collection of
periodic points of period \[n\]. In the generic case, these period \[n\]
points will all be distinct, but in special cases it may happen either that
two period \[n\] orbits come together, or that a period \[n\] orbit degenerates
to an orbit of lower period.

In any case, since the map \[f\] carries this collection of \[Nn\] points
(counted with multiplicity) into itself,
we can group these \[Nn\] points
into \[N\] orbits of formal period \[n\] (and actual period dividing
\[n\]), and form the collection of
multipliers \[\{\,\eta_1^{(n)}\,,~\cdots\,,~\eta_{N}\,
^{(n)}\}\] of these period \[n\] orbits. Let
$$	\sigma_1^{(n)}\,,~\cdots\,,~\sigma_{N}^{(n)} $$
be the elementary symmetric functions of these multipliers.

{\QP{\bf Lemma D.1.} \it Each \[\sigma_i^{(n)}\] can be expressed as
a polynomial function of the two invariants \[\sigma_1\] and \[\sigma_2\]
of \S3.\smallskip}

{\bf Proof outline.} Since \[\sigma_1\] and \[\sigma_2\] form a complete
set of conjugacy invariants, we can express \[\sigma_i^{(n)}\] as
a single valued function of \[(\sigma_1\,,\,\sigma_2)\]. This function is
continuous, since the collection of roots of a polynomial depends continuously
on the polynomial.\break
It is algebraic (that is its graph is an affine variety
in \[\C^3\]) since the construction is purely algebraic. But a continuous
algebraic function from \[\C^2\] to \[\C\] is necessarily a polynomial.\QED
\vskip - .1in

{\QP{\bf Corollary D.2.} \it Each locus \[\Per_n(\eta)\subset \M\]
is an affine algebraic variety, defined by the polynomial equation
$$ \eta^N-\sigma_1^{(n)}\eta^{N-1}+-\cdots\pm\sigma_N^{(n)}~=~0~.$$\smallskip}

The proof is immediate.\QED\smallskip

{\bf Examples.} For \[n=2\] there is just one period \[2\] orbit. According
to 3.6, its multiplier \[\eta_1^{(2)}\] is equal to
$$      \sigma_1^{(2)}~=~2\sigma_1+\sigma_2~. $$
For \[n=3\] there are two period \[3\] orbits. The sum and product of their
multipliers are given by the formulas
$$\eqalign{     \sigma_1^{(3)}~&=~\sigma_1\,(2\sigma_1+\sigma_2)+3\sigma_1+3~,
\cr        \sigma_2^{(3)}~&=~(\sigma_1+\sigma_2)^2\,(2\sigma_1+\sigma_2)
        -\sigma_1\,(\sigma_1+2\sigma_2)+12\sigma_1+28}~. $$
In particular,
\[\sigma_1^{(3)}\] is a quadratic polynomial in
\[(\sigma_1\,,\,\sigma_2)\], and \[\sigma_2^{(3)}\] is a cubic polynomial.
I don't know any neat proof of these formulas. However, they can be verified
by first studying the two multipliers \[\eta_i^{(3)}\]
asymptotically as \[\sigma_1\,,
\,\sigma_2\to\infty\] in order to determine the degree and the highest degree
terms of these polynomials, and then computing enough explicit examples to
uniquely determine the lower degree terms. Here are six
easily computed examples, which can be used to complete this
computation. Let \[\omega=e^{2\pi i/3}=(-1+i\,\sqrt 3)/2\].\bigskip

\centerline{\vbox{\hrule\hbox{\vrule\vbox{
\halign{\quad\hfil #\hfil &\qquad\hfil #\hfil & # & # &\qquad \hfil #\hfil
 & # &\hfil #\hfil\quad\cr
 $ f(z)$& $\{\mu_i\}$ & $\sigma_1$
 & $\sigma_2$ & $\{\eta_i^{(3)}\}$ &
 $\sigma_1^{(3)}$ & $\sigma_2^{(3)}$\cr
\noalign{\smallskip\hrule\smallskip}
$ z^2$ & $\{0,~0,~2\}$ & ~2 & ~0 & $\{8,~8\}$ & ~16 & 64\cr
$ z^2-1.75$ & $\{0,\,1\pm\sqrt 8\}$ & ~2 & $-7$ & $\{1,~1\}$ &
 ~2 & 1\cr
$ z^2-2$ & $\{0,-2,~4\}$ & ~2 & $-8$ & $\{8,-8\}$ & ~0 & $-64$\cr
$ 1/z^2$ & $\{-2,-2,-2\}$ & $-6$ & ~12 & $\{-8,-8\}$ & $-16$ & 64\cr
$z(z+\omega)/(\omega z+1)$
& $\{\omega,\,\omega,-2\omega\}$ & ~0 & $-3\bar\omega$ & $\{1,~1\}$ & ~2 & 1\cr
$z(z+\bar\omega)/(\bar\omega z+1)$
& $\{\bar\omega,\,\bar\omega,-2\bar\omega\}$ & ~0 & $-3\omega$ & $\{1,~1\}$ &
 ~2 & 1\cr
\noalign{\medskip\hrule}}}\vrule}}
}\medskip

\bigskip\medskip

\centerline{\bf Appendix E. Totally disconnected Julia sets in Degree \[d\].}
\smallskip

Let  \[f\] be a rational map of degree  \[d \ge 2\]. This Appendix will
prove the following two statements.

{\QP\bf Lemma E.1. {\it If all critical values of \[f\]
 belong to a single Fatou component, then the Julia set
 \[J\] of \[f\]  is totally disconnected.}\smallskip}

{\bf Remark.} Both
Przytycki [Pr2] and Makienko [Mak] have announced the sharper result that
such a totally disconnected Julia set containing
no critical points is homeomorphic, with its dynamics, to the one-sided
shift on \[d\] symbols.
(For the polynomial case compare [BDK], and for the quadratic case
compare [GK].)
For a study of totally disconnected Julia sets which do contain critical
points, see [BH].

{\QP{\bf Lemma E.2.} \it Whenever the Julia set \[J\] is
totally disconnected, all orbits in the Fatou set \[\C\ssm J\]
converge, either:

\noindent$~~(1)$ to an attracting fixed point, or

\noindent$~~(2)$ to a parabolic fixed point of multiplicity two.\smallskip}

{\bf Remark.} A fixed point \[z_0\ne\infty\] has multiplicity equal to
two if and only if the first two derivatives satisfy \[f'(z_0)=1\] and
\[f''(z_0)\ne 0\], so that there
is exactly one attracting petal in the associated Leau-Fatou flower.

Combining these two statements, we evidently obtain Lemma 8.1, as stated
in \S8.\medskip

{\bf Proof of E.1} (modeled on [R3]).
Suppose that all critical values belong to the Fatou component
\[U_1\] . Then \[U_1\] must be an immediate basin for an attracting or
parabolic periodic point, and
must belong to a cycle $$ U_1 \to U_2 \to \cdots \to U_p \to U_1 $$
of Fatou components with period \[p\ge 1\].
Evidently every Fatou component must eventually fall onto this cycle.
(Compare [Su].)\smallskip

{\bf Hyperbolic case:} If \[f\]  is hyperbolic, then the closure of the
post-critical set \[X\] splits as a disjoint union \[ \overline X=\overline X_1\cup\cdots\cup \overline X_p\]
where \[\overline X_i\] is a compact subset of \[U_i\]. Choose a closed topological
disk \[D\subset U_1\]  with  \[\overline X_1\]  contained in the interior of \[D\],
and let \[\Gamma\] be the boundary of \[D\]. Note that every critical value
of \[f\] is contained in \[X_1\subset D\], and note that \[\Gamma\] separates
\[X_1\] from the Julia set.\smallskip


{\bf Step 1.} We prove that no component
\[\Gamma_n\]  of \[f^{ -n}(\Gamma)\]  can separate points of any  \[X_i\].
If \[\Gamma_n\] is not contained in \[U_i\], then \[\Gamma_n\] certainly
cannot separate \[X_i\]. For
\[\Gamma_n\subset U_i\], we will prove this statement
by induction on \[n\]. If \[1\le i< p\], then
since \[f\] induces a homeomorphism from the pair \[(U_i\,,\,X_i)\] onto
\[(U_{i+1}\,,\,X_{i+1})\], the statement follows by induction.
Thus we may suppose that \[i=p\], so that
\[\Gamma_n\subset U_p\] with \[f(\Gamma_n)
=\Gamma_{n-1}\subset U_1\].
Suppose inductively that \[\Gamma_{n-1}\] does not separate \[X_1\].
Choose a small disk  \[\Delta\subset U_1\]  containing  \[X_1\]
and disjoint from \[\Gamma_{n-1}\], and let \[\Delta'\] be the closure
of the complementary disk \[\hat C\ssm\Delta\]. Since
\[\Delta'\] contains no critical values, it follows that \[f^{-1}(\Delta')\]
is a union of \[d\] disjoint disks, and hence that its complement
\[f^{-1}(\Delta)\] is connected. Therefore the loop
\[\Gamma_n\subset f^{-1}(\Delta')\]
cannot separate points of \[X_p\subset f^{-1}(\Delta)\]. This proves Step 1.

In particular, since \[\Gamma_{n-1}\] cannot separate two critical values,
it follows that each \[\Gamma_n\] maps homeomorphically onto \[f(\Gamma_n)=
\Gamma_{n-1}\].

\smallskip {\bf Step 2.} We use the Poincar\'e metric on the complement
of \[\overline X\]. Since the iterated pre-images of  \[\Gamma\]  remain
bounded away from \[\overline X\], and since each component of
\[f^{-n}(\Gamma)\] maps homeomorphically
onto \[\Gamma\], it follows that the Poincar\'e
diameters of the various components of \[f^{ -n}(\Gamma)\]
shrink uniformly to zero as \[n\to \infty\]. Therefore, the diameters of these
components in the spherical metric shrink to zero also.
\smallskip

{\bf Step 3.} If \[n\equiv 0~\mod p\], so that \[f^{\circ n}(X_1)\subset X_1\],
then \[f^{-n}(\Gamma)\] separates \[X_1\] from the Julia set. For otherwise,
if there were a path \[\pi\] from a point
of \[X_1\] to a point of \[J\] in the complement of \[f^{-n}(\Gamma)\],
then \[f^{\circ n}(\pi)\] would be a path from \[X_1\] to \[J\] in the
complement of \[\Gamma\], which is impossible.

\smallskip

{\bf Step 4.} Proof that \[J\] is totally disconnected, and that
\[U_1=\hat C\ssm J\]. Using the spherical metric,
choose  any \[\epsilon\] smaller than the diameter of  \[X_1\],
and choose  \[n\equiv 0~\mod p\]  so large that each component
of \[f^{-n}(\Gamma)\]  has diameter
less than  \[\epsilon\] . Since  \[X_1 \]  cannot fit inside
any of these loops, it follows that every component of \[J\] must fit inside
one of them. Similarly, any other Fatou component \[U'\] must
fit inside one of these loops. Since \[\epsilon\] can be arbitrarily small,
this proves that \[J\] is totally disconnected, and that \[U_1\]
is the unique Fatou component. In particular, \[p=1\].\smallskip

{\bf Parabolic case:} Here the argument is more delicate, since the
post-critical closure \[\overline X\] intersects the Julia set in a parabolic cycle.
Again we can choose a disk \[D\] which contains \[X_1\], the set of
post-critical points
of \[U_1\], in its interior, but now this disk \[D\] must have the
parabolic point of \[\overline X_1\]
on its boundary. Steps 1, 3, 4 go through very much
as in the hyperbolic case, but Step 2 requires more work. Again we use
the Poincare metric on the complement of \[\overline X\] .
This decreases under every branch of  \[f^{-1}\] , and decreases by a
definite factor less than 1 throughout any region bounded away from  \[\overline X\].
Let \[P\] be the union of suitably chosen repelling petals at the points of our
parabolic cycle. Then there is a unique branch of \[f^{-1}\] which maps
\[P\] into itself. For any compact curve segment in \[P\], the spherical
lengths of the successive images under this branch of \[f^{-1}\]
shrink to zero. Now consider an arbitrary sequence of  pre-images
$$ f:~	\cdots~\buildrel\cong\over\to~ \Gamma_2~\buildrel\cong\over\to~\Gamma_1~\buildrel\cong\over\to~\Gamma_0=\Gamma~, $$
where each \[\Gamma_n\] is a component of \[f^{-1}(\Gamma)\]. We must prove
that the diameter of \[\Gamma_n\]\break (in the spherical metric) shrinks to
zero as \[n\to\infty\]. Let \[g_n:\Gamma\buildrel\cong\over\to\Gamma_n\]
be the branch of \[f^{-n}\] which maps \[\Gamma\] to \[\Gamma_n\].

For any point \[z\in\Gamma\] there are two possibilities. If \[g_n(z)\in P\]
for all \[n\ge n_0\], then the same is true for all points in some neighborhood
\[N\] of \[z\]. (Here \[n_0\] is constant throughout this neighborhood.) It
follows that the
spherical length of \[g_n(N)\] shrinks to zero as \[n\to\infty\]. In this
way, we control the points belonging to some open subset of \[\Gamma\].

For points \[z\] in the complementary closed subset \[Y\subset\Gamma\],
there are
infinitely many images \[g_n(z)\] which do not lie in \[P\], and hence
are bounded away from the post-critical
closure \[\overline X\]. It is easy to check that the Poincar\'e arclength
of \[g_n(Y)\] is finite for large \[n\]. Thus the Poincar\'e length of
\[g_n(Y)\] can be expressed as the integral over \[Y\] of a function which
decreases monotonically to zero as \[n\to\infty\]. Hence
this length shrinks to zero as \[n\to\infty\], and
the spherical arc length of \[g_n(Y)\] must
also shrink to zero. The rest of the proof
proceeds as before.\QED\smallskip

{\bf Proof of E.2.} If \[J\] is totally disconnected and contains a parabolic
periodic point \[z_0\], then \[z_0\] must be
a fixed point must have multiplicity two. For
otherwise there would be more than one parabolic basin; but if \[J\] is totally
disconnected, then the Fatou set must be connected.\QED
\vfil\eject

\centerline{\bf Appendix F. A ``Sierpinski carpet'' as Julia set.}\smallskip

\centerline{written in collaboration with Tan Lei\footnote{$^1$}
{Address: Laboratoire Math., \'Ecole Norm. Sup.,
69364 Lyon, France.}}
\medskip

By definition, a compact subset of the plane is called a {\bit Sierpinski
carpet}\footnote{$^2$}{The term ``Sierpinski curve'' is commonly used
in the literature. However Mandelbrot's term ``Sierpinski carpet'' seems
more descriptive, and serves to distinguish this object from
other examples of fractals due to Sierpinski.}
if it is connected, locally
connected, nowhere dense, and has the property that the boundaries of the
various
complementary components are pairwise disjoint Jordan curves.
Any two Sierpinski carpets are homeomorphic to each other. (See [Wh]).
The standard example is obtained from the unit square by
subdividing into nine subsquares, removing the interior of the middle
subsquare, then removing the middle ninth of each of the remaining eight
subsquares, and continuing inductively.

It seems to be widely known among experts that such a Siepinski carpet
can occur as the Julia set for a rational function, but we are not aware
of any explicit example in the literature.

Wittner [W] noted that there is one and only one real quadratic rational map,
up to conjugacy, such that the two critical orbits have periods three and four
respectively. He showed that this example cannot be obtained by mating two
quadratic polynomials. {\bit We will show that the Julia set for this
rational map is homeomorphic to a Sierpinski carpet.}

\midinsert
\insertRaster 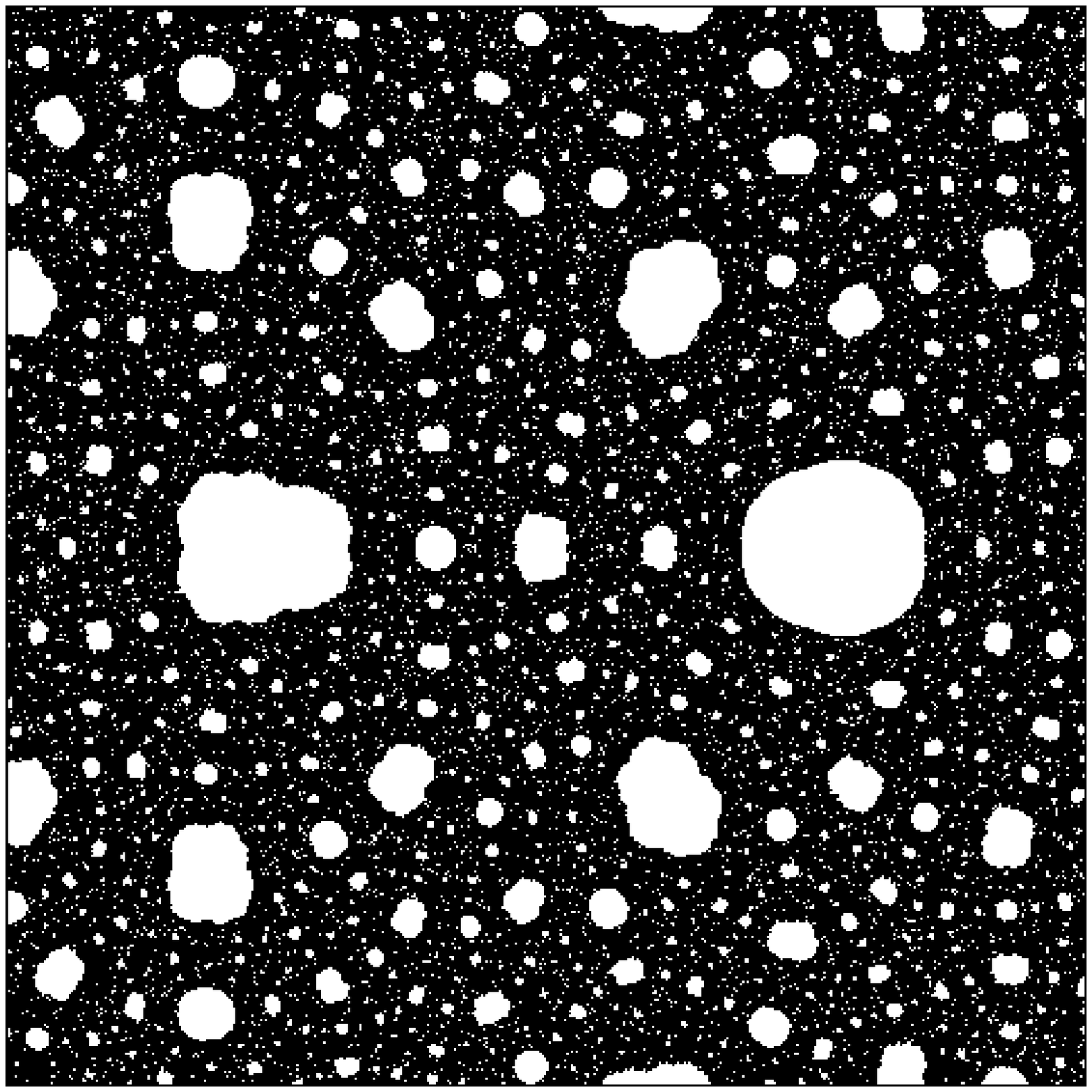 pixels 480 by 480 scaled 450
\centerline{Figure 18. The Julia set \[J(g)\].}
\endinsert

{\bf The example.} Consider a real quadratic map \[f(z)= a(z+z^{-1})+b\],
normalized so that the critical points are at \[\pm 1\], and so that there
is a fixed point of multiplier \[1/a\] at infinity. Computation shows that
there are unique values \[a= -.138115091\],\break \[b= -.303108805\] so that
the critical point \[x_0=1\] has period 4 and so that the critical point
\[y_0=-1\] has period 3. (The critical orbits are arranged along the real
line as \[x_3<y_0<x_1<y_1<x_2<x_0<y_2\], where for example
\[x_0\mapsto x_1\mapsto x_2
\mapsto x_3\mapsto x_0\].) It will be convenient to conjugate
\[f\] by the Moebius transformation
\[z\leftrightarrow (z+i)/(iz+1)\] which interchanges
the real axis \[\R\cup\infty\] and the unit circle \[S^1\]. The critical orbits
for the resulting map \[g\] are shown in Figure 19, and the corresponding
Julia set, drawn to the same scale, is shown in Figure 18. Note that the
critical values \[x_1\] and \[y_1\] divide the unit circle \[S^1\] into a
shorter arc \[I\] and a longer arc \[J=S^1\ssm {\rm interior}\,I\]. The map \[g\] carries
both the interval \[[-1,1]\] and the complementary interval \[\R\cup\infty\ssm
(-1,1)\] homeomorphically onto \[I\], and carries both the upper and lower unit
semi-circles homeomorphically onto \[J\]. It has a repelling fixed point at
the bottom point \[-i\] of the lower semi-circle.

\midinsert
\insertRaster 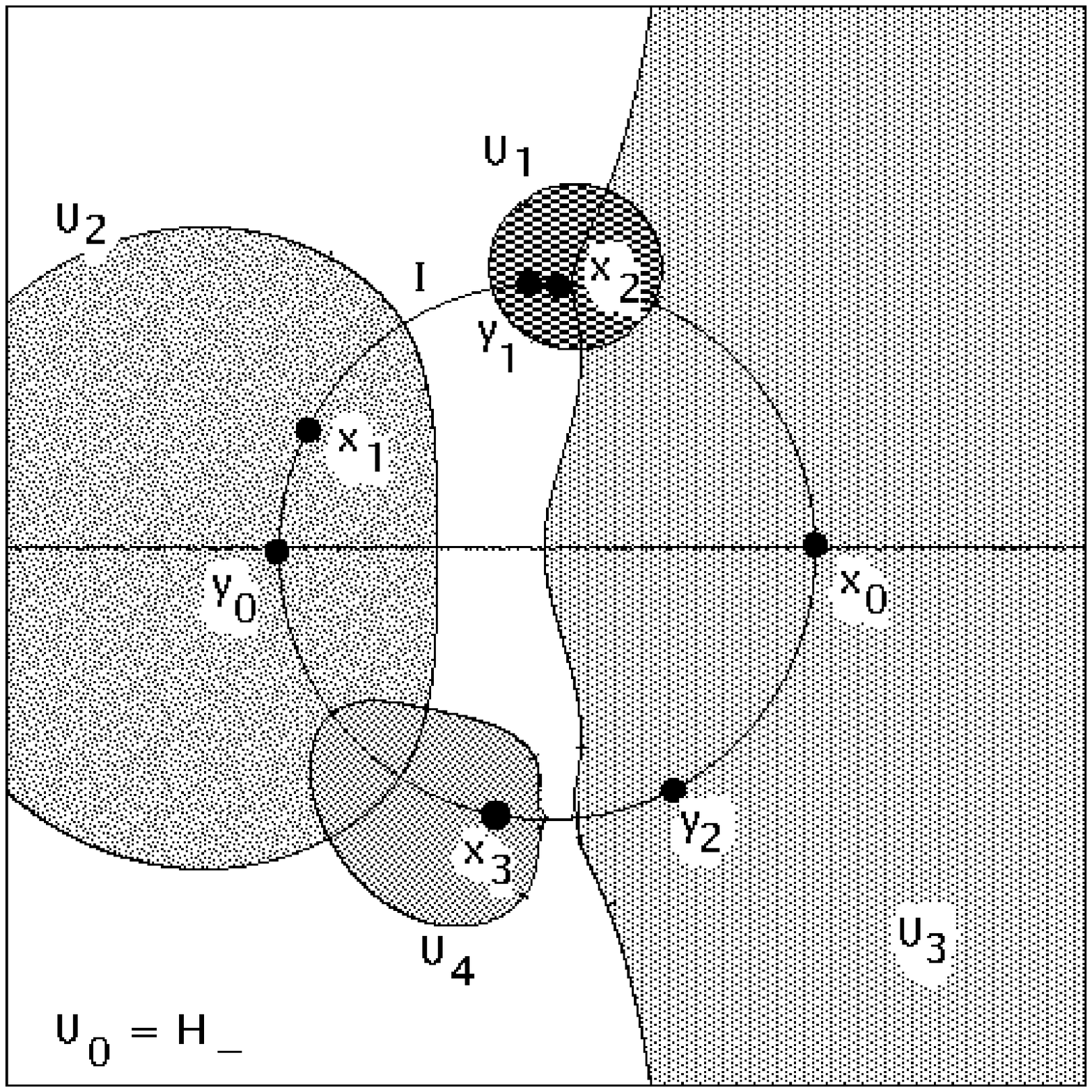 pixels 480 by 480 scaled 450
\centerline{Figure 19. Critical orbits for \[g\].}
\endinsert

It follows from Lemma 8.2 that the Julia set \[J=J(g)\] it connected,
so that each Fatou component is simply-connected. Since \[g\] is hyperbolic
on \[J\], it follows from 7.1 that \[J\] is locally connected.
In order to prove that each Fatou component is bounded by a Jordan
curve\footnote{$^1$}
{The property that every Fatou component is bounded by a Jordan curve
seems very common among hyperbolic rational maps with \[J\] connected.
It would be interesting to decide just when it is satisfied. Note however
that the boundary of a fixed fully invariant Fatou component (for example
the component containing infinity for a polynomial map) is usually
not a Jordan curve.}
we proceed as follows.

We will first
construct a simply connected neighborhood \[U_4\] of the point \[x_3\]
so that the 4-fold composition \[g^{\circ 4}\] maps \[U_4\] by a
polynomial-like map of degree 4 onto the lower half-space \[H_-\].
(For definitions, see [DH3].) In fact, for \[0\le j\le 4\], 
let \[U_j\] be the component of \[g^{-j}(H_-)\] which contains
the point \[x_{3-j}\] (where \[x_{-1}=x_3\]), so that
$$	U_4~\to~ U_3~\to~ U_2~\to~ U_1~\to~ U_0=H_-~~. $$
Note that the set \[U_0=H_-\]
is invariant under the inversion \[z\mapsto 1/\bar z\] in the circle \[S^1\].
Since the map \[g\] commutes with this inversion, and since each \[U_j\]
contains a point \[x_{3-j}\] of the circle, it follows inductively that
each \[U_j\] is invariant under inversion. Furthermore, since \[g(z)=g(1/z)\],
it follows that each \[g^{-1}(U_j)\] is invariant under \[z\mapsto 1/z\],
or equivalently under complex conjugation. We will also show inductively
that each \[U_j\] is simply-connected. We can distinguish two cases.
If \[U_j\] contains one critical value \[x_1\] or \[y_1\], then \[U_{j+1}=
g^{-1}(U_j)\] is a ramified 2-fold covering of \[U_j\]. In this case,
\[U_{j+1}\] must contain the real point \[x_0\] or \[y_0\], and must be
invariant under complex conjugation. On the other
hand, if \[U_j\] contains no critical value, then \[g^{-1}(U_j)\] splits into
two simply-connected components. Neither can intersect the real axis,
so one must lie in the upper half-plane and the other must be its complex
conjugate in the
lower half-plane. Just one of these two is the component \[U_{j+1}\] containing
\[x_{3-(j+1)}\]. (If \[U_j\] contained both critical values, then
\[g^{-1}(U_j)\] would not be simply connected; however, this case does not
occur.) It is now easy to check that the sets \[U_j\] are situated as
illustrated in Figure 19, containing post-critical points as shown,
and that the proper maps
\[	U_4\to U_3\to U_2\to U_1\to U_0 \]
have degrees \[1,\;2,\;2,\;1\] respectively. Thus the composition is proper
of degree 4.

Proof that the
region \[U_4\] is compactly contained in \[U_0=H_-\]. Note that the
successive images of \[\R\cup\infty\] under \[g\] are intervals around
the unit circle as follows:
$$	\R\cup\infty~~\to~~ I=[y_1\,,\,x_1]~~\buildrel\cong\over\to~~ [y_2\,
	,\,x_2]~~\to~~	[x_1\,,\, x_3]~~\to~~ [x_0\,,\,y_1]~, $$
where \[[\alpha\,,\,\beta]\] stands for the interval from \[\alpha\] to
\[\beta\] in the counter-clockwise direction around \[S^1\].
Suppose that the closure \[\bar U_4\] contains some point
\[u\] belonging to\break \[\partial H_-=\R\cup\infty\]. Then the third
forward image \[g^{\circ 3}(u)\] must belong
to the interval \[[x_1\,,\,x_3]\subset S^1\]. But it must also belong
to \[g^{\circ 3}(\bar U_4)=\bar U_1\], which is contained in the upper
half-space, hence it must belong
to the smaller interval \[[x_1\,,\,y_0]\]. Hence \[g^{\circ 4}(u)\in\bar U_0\]
must belong to \[g[x_1\,,\,y_0]=[x_2\,,\,y_1]\subset H_+\], which is
impossible. Thus \[g^{\circ 4}:
U_4\to U_0\] is polynomial-like.\smallskip

Let \[X_j\] be the immediate attractive basin for the super-attracting
point \[x_j\] (or equivalently let \[X_j\] be the Fatou component containing
the point \[x_j\]). Note that each \[X_j\] is invariant under the inversion
\[z\mapsto 1/\bar z\].
Note also that \[X_1\] and \[X_2\] are contained in the
upper half-plane \[H_+\], while \[X_3\] is contained in the lower half-plane
\[H_-\]. For if one of these sets intersected the real axis, then its image
\[X_2\] or \[X_3\] or \[X_0\] would intersect the interval \[g(\R)=I\]. Hence
this image, which is simply-connected and invariant under inversion, would
have a non-connected intersection with the unit circle; which is impossible.
Since \[X_3\subset H_-=U_0\], it follows by induction on \[i\] that
\[X_{3-i}\subset U_i\] for \[i=1,2,3,4\]
(taking the subscript \[3-i\] modulo 4). In particular, \[X_3\subset U_4\].

Let \[K=K(g^{\circ 4}| U_4)\] be the filled
Julia set for the polynomial-like mapping\break \[g^{\circ 4}:U_4\to U_0\].
By definition, \[K\] consists of those
points whose successive images under \[g^{\circ 4}\]
remain within the set \[U_4\]. 
Since \[X_3\subset U_4\], and since \[g^{\circ 4}\] maps the basin \[X_3\]
onto itself, it follows that \[~X_3\subset K\subset U_4~\].

(In fact \[X_3\] is precisely equal to that component of the interior of \[K\]
which contains the super-attractive point \[x_3\].
Interior points of \[K\] certainly always
belong to the Fatou set \[\hat\C\ssm J(g)\]. On the other hand,
it is not difficult to check that boundary points of \[K\] belong to
the Julia set \[J=J(g)\]. It follows that
each component of the interior of \[K\] is a Fatou component for \[g\],
that is a connected component of \[\C\ssm J(g)\].)

Proof that \[X_3\] is bounded by a Jordan curve. Since \[J\] is locally
connected, it follows from a Theorem of Caratheodory that
the various internal rays from \[x_3\] land at well defined
points of \[\partial X_3\], which depend continuously on the internal
angle. (Compare the discussion in [DH2, p. 20].) Thus,
to show that \[\partial X_3\] is a Jordan curve, we need only
show that no two internal rays land at a common point.
But two rays from \[x_3\] landing at a common point \[z\in\partial X_3\]
would form a simple closed curve \[\Gamma\]
within the closure \[\bar X_3\subset K\subset U_4\].
Since \[K\] is a full set (that is, since \[\hat\C\ssm K\] is
connected), it would follow that \[\Gamma\] must bound a region \[W\subset K\],
which is necessarily contained in the Fatou set
\[\hat\C\ssm J\]. Now some
intermediate ray from \[x_3\] must land at a distinct point \[z'\ne z\]
of \[\partial
X_3\]. (It follows from the reflection principle that there cannot be
an entire interval of internal angles with a common landing point.)
Evidently this landing point \[z'\] must belong both to the region
\[W\subset\hat\C\ssm J\], and to the boundary \[\partial X_3\subset J\],
which is impossible.

It follows easily that every iterated pre-image of \[X_3\] is also bounded
by a Jordan curve.
Now let \[Y_i\] be the Fatou component containing \[y_i\]. The proof that
\[Y_2\] is bounded by a Jordan curve is completely analogous. For this
proof, we construct a neighborhood \[V_6\] of \[y_2\] so that the map
\[g^{\circ 6}:V_6\to V_0=H_-\] is polynomial-like of degree 8. In fact
let \[V_i\] be the component of \[g^{-i}(H_-)\] which contains \[y_{5-i}\]
(taking the subscript \[5-i\] modulo 3). Then \[V_i=U_{i}\] for \[i\le 3\]
However \[V_4\] is the complex conjugate (and the reciprocal) of
\[U_4\]. The pre-image \[g^{-1}(V_4)=V_5\] is a neighborhood of \[y_0\]
containing no critical values, so its pre-image is the union of a component
\[V_6\] containing \[y_2\] and a complex conjugate component. A similar
argument shows that \[V_6\] is compactly contained in \[V_0=H_-\], and it
follows that \[Y_2\subset V_6\] is bounded by a Jordan curve. Since every Fatou
component is an iterated pre-image of either \[X_3\] or \[Y_2\], this
proves that every Fatou component is bounded by a Jordan curve.

Proof that the various regions \[X_i\] and \[Y_j\] have no boundary
points in common. We noted above that the regions \[X_1\] and \[X_2\]
are contained in the upper half-plane \[H_+\], and a similar argument shows
that \[Y_1\subset H_+\]. On the other hand, \[X_3\] and \[Y_2\]
are compactly contained in \[H_-\]. This proves that the closures \[\bar X_3\]
and \[\bar Y_2\] are disjoint from \[\bar X_1\] and \[\bar X_2\] and \[\bar
Y_1\]. It then follows easily that all of the sets \[\bar X_i\] and \[\bar
Y_j\] are pair-wise disjoint. For example, if \[\bar X_2\] intersected
\[\bar Y_1\] then \[g^{\circ 4}(\bar X_2)=\bar X_2\] would have to intersect
\[g^{\circ 4}(\bar Y_1)=\bar Y_2\].

Now suppose that any pair of Fatou components \[F\] and \[F'\] had a boundary
point in common. We claim first that \[g(F)\ne g(F')\]. For otherwise,
if \[g(F)=g(F')\], then two different internal rays from the center point
of \[g(F)\] would land on a common boundary point. Evidently
it would follow that the boundary of \[g(F)\] is not
a Jordan curve, contradicting our previous conclusion. Thus \[g(F)\ne g(F')\]
must be distinct Fatou components, with a boundary point in common. Taking
successive forward images under \[g\], we must eventually find that two
of the \[X_i\] and \[Y_j\] have a boundary point in common, which is
impossible. This completes the proof that \[J(g)\] is a Sierpinski carpet.\QED
\eject

\noindent{\bf References.}\smallskip

\ref [Bie],~~ B. Bielefeld, Questions in quasiconformal surgery,
in ``Problems in holomorphic dynamics'', ed. Bielefeld and Lyubich,
Stony Brook IMS preprint 1992/7 (see also 1990/1); to appear in
``Linear and Complex Analysis Problem Book'', ed. Havin and Nikolskii,
Springer Verlag.

\ref [BDK], P. Blanchard, R. Devaney and L. Keen, The dynamics of complex
polynomials and automorphisms of the shift, Invent. Math {\bf 104} (1991)
545-580.

\ref [BH] B. Branner and J. H. Hubbard, The iteration of cubic polynomials,
Part I: the global topology of parameter space, Acta Math. {\bf 160} (1988)
143-206; Part II: patterns and parapatterns, Acta Math., to appear.

\ref [Bou] T. Bousch, Sur quelques probl\`emes de dynamique holomorphe,
Th\`ese, Univ. Paris-Sud Orsay 1992.

\ref [CCMM] F. Cohen, R. Cohen, B. Mann, R. J. Milgram, The topology of
rational functions and divisors of surfaces, Acta Math. {\bf 166} (1991)
163-221.

\ref [DH1] A. Douady and J. H. Hubbard, It\'eration des polyn\^omes
quadratiques complexes, CRAS Paris {\bf 294} (1982) 123-126.

\ref [DH2] A. Douady and J. H. Hubbard, \'Etude dynamique des polyn\^omes
complexes I \& II, Publ. Math. Orsay (1984-85).

\ref [DH3] A. Douady and J. H. Hubbard, On the dynamics of polynomial-like
mappings, Ann. Sci. Ec. Norm. Sup. (Paris) {\bf 18} (1985), 287-343.

\ref [DM] P. Doyle and C. McMullen, Solving the quintic by iteration, Acta
Math. {\bf 163} (1989) 151-180.


\ref [GK] L. Goldberg and L. Keen, The mapping class group of a generic
quadratic rational map and automorphisms of the 2-shift, Invent. Math.
{\bf 101} (1990) 335-372.

\ref [HP] F. v. Haeseler and H.-O. Peitgen, Newton's method and complex
dynamical systems, Acta Appl. Math. {\bf 13} (1988) 3-58.

\ref [HW] G. H. Hardy and E. M. Wright, An Introduction to the Theory of
Numbers, Clarendon Press 1938, 1945, 1954.

\ref [La] P. Lavaurs, Syst\`emes dynamiques holomorphes: explosion de points
p\'eriodiques para\-boliques, Th\`ese, Univ. Paris-Sud Orsay 1989.

\ref [Ly] M. Lyubich, An analysis of the stability of the dynamics of rational
functions, Funk. Anal. i. Pril. {\bf 42} 1984), 72-91;
Selecta Math. Sovietica {\bf 9} (1990) 69-90.

\ref [Mak] P. Makienko, in preparation.

\ref [Man] B. Mandelbrot, The Fractal Geometry of Nature, Freeman 1982.

\ref [Mc] C. McMullen, Automorphisms of rational maps, in
``Holomorphic Functions and Moduli'', edit. Drasin et al., Springer 1988.

\ref [M1] J. Milnor, On the 3-dimensional Brieskorn manifolds \[M(p,q,r)\], 
pp. 175-225 of ``Knots, Groups, and 3-Manifolds'', ed. Neuwirth, Annals Stud.
{\bf 84}, Princeton U. Press 1975.

\ref [M2] ------, Dynamics in one complex variable: Introductory lectures,
Stony Brook IMS preprint 1990/5.

\ref [M3] ------, Remarks on iterated cubic maps, Stony Brook IMS
preprint 1990/6; to appear in ``Experimental Mathematics''.

\ref [M4] ------, Hyperbolic components in spaces of polynomial maps,
Stony Brook IMS\break preprint 1992/3.

\ref [M5] ------, On cubic polynomials with periodic critical point, in
preparation.

\ref [MSS] R. Ma\~n\'e, P. Sad and D. Sullivan, On the dynamics of rational
maps, Ann. Sci. \'Ec. Norm. Sup. Paris (4) {\bf 16} (1983) 193-217.

\ref [Pe] C. Petersen, No elliptic limits, in preparation.

\ref [Pr1] F. Przytycki, Remarks on the simple connectedness of basins of sinks
for iteration of rational maps, Dyn. Sys. \& Erg. Th., Banach Center Pub.
{\bf 23} (1989) 229-235.

\ref [Pr2] ------, Polynomials in the hyperbolic components, in preparation.

\ref [R1] M. Rees, Ergodic rational maps with dense critical point forward
orbit, Erg. Th. \& Dy. Sy. {\bf 4} (1984) 311-322.

\ref [R2] ------, Positive measure sets of ergodic rational maps, Ann. Sci.
\'Ecole Norm. Sup. (4) {\bf 19} (1986) 383-407.

\ref [R3] ------, Components of degree two hyperbolic rational maps,
Invent. Math. {\bf 100} (1990), 357-382.

\ref [R4] ------, Combinatorial models illustrating variation of
dynamics in families of rational maps, Proc. Int. Congr. Math. Kyoto 1990.

\ref [R5] ------, A partial description of parameter space of rational
maps of degree two,\break Part I: Acta Math. {\bf 168} (1992) 11-87;
Part II: Stony Brook IMS preprint 1991/4.

\ref [Se] G. Segal, The topology of spaces of rational functions, Acta Math.
{\bf 143} (1979) 39-72.

\ref [Sh1] M. Shishikura, On the quasiconformal surgery of rational functions,
Ann. Sci. \'Ec. Norm. Sup. {\bf 20} (1987) 1-29.

\ref [Sh2] M. Shishikura, On a theorem of M. Rees for the matings of
polynomials, preprint, IHES 1990.

\ref [Su] D. Sullivan, Quasiconformal homeomorphisms and dynamics I, solution
of the Fatou-Julia problem on wandering domains, Ann. Math. {\bf 122}
(1985) 401-418.

\ref [Ta] Tan Lei, Accouplements des polyn\^omes complexes, Th\`ese, Orsay
1987; Mating of quadratic polynomials, to appear.

\ref [Th] W. Thurston, Three-dimensional Geometry and Topology, to appear.

\ref [vW] B. L. van der Waerden, Modern Algebra I (various editions).

\ref [Wh] G. T. Whyburn, Topological
characterization of the Sierpinski curve, Fundamenta Math. {\bf 45} (1958)
320-324.

\ref [W]~B. Wittner,  On the bifurcation loci of rational maps of
degree two, Thesis, Cornell University 1988.

\ref [Y] ~Yin Yongcheng, On the Julia sets of quadratic rational maps,
Complex Variables {\bf 18} (1992) 141-147.

\vfill

\end